\documentclass[10pt]{amsart}

\usepackage{graphicx}        
\usepackage{multicol}        
\usepackage[bottom]{footmisc}
\usepackage{amscd,amsmath,amssymb,amsfonts}
\usepackage[cmtip, all]{xy}

\usepackage{bbm}
\usepackage{tikz-cd}
\usepackage{enumerate}
\usepackage{mathrsfs}
\usepackage{amsthm}

\newtheorem{theorem}{Theorem}[section]
\newtheorem{proposition}[theorem]{Proposition}
\newtheorem{lemma}[theorem]{Lemma}

\newtheorem{corollary}[theorem]{Corollary}

\theoremstyle{definition}
\newtheorem{definition}[theorem]{Definition}
\theoremstyle{remark}
\newtheorem{remark}[theorem]{Remark}

\newtheorem{ex}[theorem]{Example}

\numberwithin{equation}{section}

\newcommand{\red}{{\rm red}}

\newcommand{\Pic}{{\rm Pic}}
\newcommand{\End}{{\rm End}}
\newcommand{\Hom}{{\rm Hom}}
\newcommand{\im}{{\rm im}}
\newcommand{\Spec}{{\rm Spec\,}}

\newcommand{\Char}{{\rm char}}

\newcommand{\0}{\emptyset}
\newcommand{\sA}{{\mathcal A}}
\newcommand{\sB}{{\mathcal B}}
\newcommand{\sC}{{\mathcal C}}

\newcommand{\sE}{{\mathcal E}}
\newcommand{\sF}{{\mathcal F}}
\newcommand{\sG}{{\mathcal G}}
\newcommand{\sH}{{\mathcal H}}
\newcommand{\sI}{{\mathcal I}}
\newcommand{\sJ}{{\mathcal J}}
\newcommand{\sK}{{\mathcal K}}
\newcommand{\sL}{{\mathcal L}}
\newcommand{\sM}{{\mathcal M}}
\newcommand{\sN}{{\mathcal N}}
\newcommand{\sO}{{\mathcal O}}
\newcommand{\sP}{{\mathcal P}}
\newcommand{\sQ}{{\mathcal Q}}

\newcommand{\sT}{{\mathcal T}}
\newcommand{\sU}{{\mathcal U}}
\newcommand{\sV}{{\mathcal V}}
\newcommand{\sW}{{\mathcal W}}

\newcommand{\A}{{\mathbb A}}

\newcommand{\F}{{\mathbb F}}
\newcommand{\G}{{\mathbb G}}

\renewcommand{\L}{{\mathbb L}}

\renewcommand{\P}{{\mathbb P}}
\newcommand{\Q}{{\mathbb Q}}
\newcommand{\R}{{\mathbb R}}

\newcommand{\U}{{\mathbb U}}

\newcommand{\X}{{\mathbb X}}

\newcommand{\Z}{{\mathbb Z}}

\renewcommand{\det}{\operatorname{det}}

\newcommand{\id}{{\operatorname{\rm Id}}}

\newcommand{\<}{\langle}
\renewcommand{\>}{\rangle}

\renewcommand{\dim}{{\operatorname{\rm dim}}}

\newcommand{\Coh}{\operatorname{Coh}}

\newcommand{\Hilb}{\operatorname{Hilb}}

\newcommand{\del}{\partial}

\renewcommand{\max}{{\operatorname{\rm max}}}

\newcommand{\Sym}{{\operatorname{Sym}}}

\newcommand{\Tor}{{\operatorname{\rm Tor}}}

\newcommand{\Gr}{{\operatorname{\rm Gr}}} 
\newcommand{\rnk}{{\operatorname{rank}}} 
\newcommand{\Bl}{\text{Bl}}

\newcommand{\KQ}{{\operatorname{KQ}}}

\newcommand{\GW}{{\operatorname{GW}}}

\newcommand{\sHom}{\mathcal{H}om}
 
\newcommand{\Aut}{{\operatorname{Aut}}}

\newcommand{\GL}{\operatorname{GL}}
\newcommand{\SL}{\operatorname{SL}}

\newcommand{\BSL}{\operatorname{BSL}}

\newcommand{\Deg}{\text{deg}}

\newcommand{\ind}[1]{}

\newcommand{\inp}[1]{}

\newcommand{\Ann}{{\operatorname{Ann}}}

\newcommand{\cm}{{\operatorname{cm}}}
\newcommand{\ncm}{{\operatorname{ncm}}}
\newcommand{\sEnd}{\mathcal{E}nd}

\newcommand{\sym}{{\operatorname{sym}}}
\newcommand{\wt}{\operatorname{wt}}

\newcommand{\floor}[1]{\lfloor{#1}\rfloor}

\makeatletter
 \@namedef{subjclassname@2020}{\textup{2020} Mathematics Subject Classification}
\makeatother

\begin{document}

\date{ \today}

\author[M.~Levine]{Marc~Levine}
\address{Marc~Levine, Universit\"at Duisburg-Essen,
Fakult\"at Mathematik, Campus Essen, 45117 Essen, Germany}
\email{marc.levine@uni-due.de}

\author[S.~Pauli]{Sabrina Pauli}
\address{Sabrina Pauli, Universit\"at Duisburg-Essen,
Fakult\"at Mathematik, Campus Essen, 45117 Essen, Germany}
\email{sabrina.pauli@uni-due.de}

\thanks{Both authors have been supported by the ERC programme QUADAG.  This paper is part of a project that has received funding from the European Research Council (ERC) under the European Union's Horizon 2020 research and innovation programme (grant agreement No. 832833).\\ 
\includegraphics[scale=0.08]{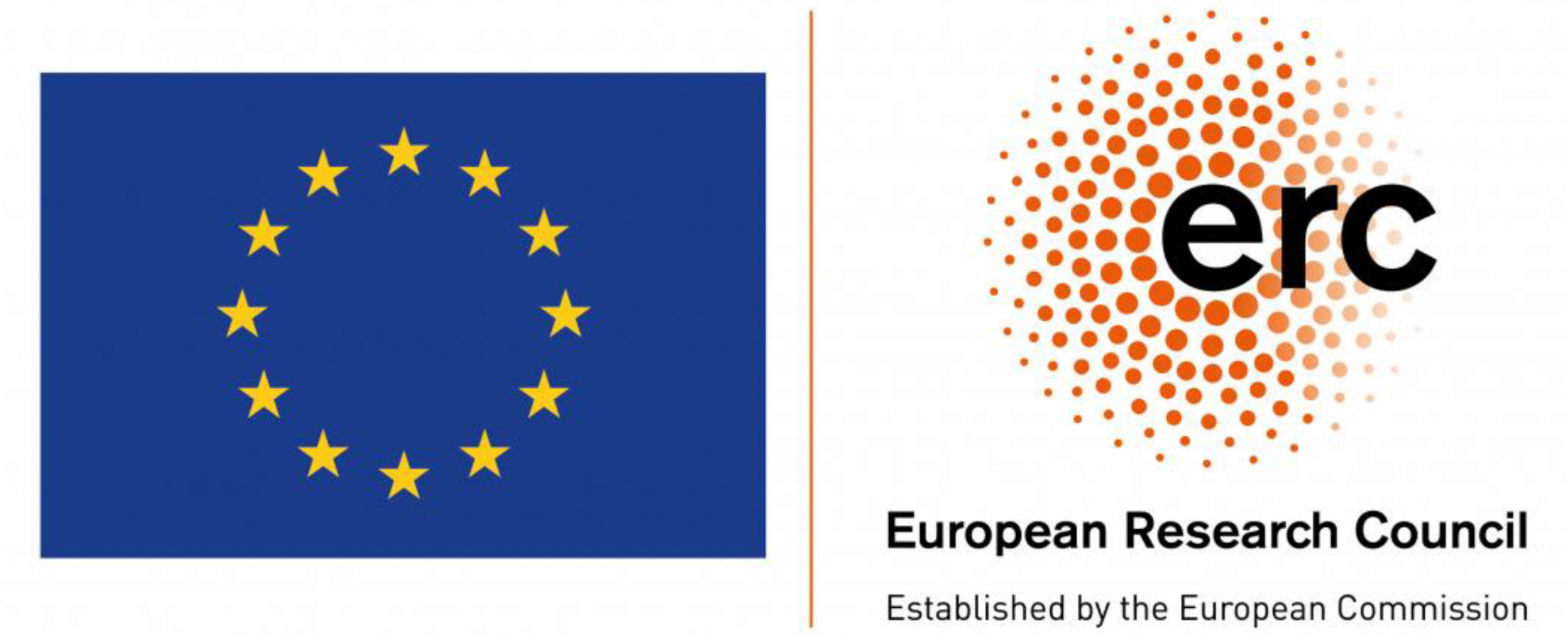}}

\title{Quadratic Counts of Twisted Cubics}


\setcounter{tocdepth}{1}

\begin{abstract}  
Using a quadratic version of the Bott residue theorem, we give a quadratic refinement of the count of twisted cubic curves on hypersurfaces and complete intersections in a projective space.  
\end{abstract}

\maketitle

\tableofcontents

\section{Introduction}  In \cite{ES} Ellingsrud and Str{\o}mme  compute the number of twisted cubic curves on some smooth projective hypersurfaces and complete intersections. Their arguments rely on understanding the closure $H_n$ in an appropriate Hilbert scheme of the locus of smooth twisted cubics in $\P^n$,  and some vector bundles on $H_n$ arising in the resolution of the ideal sheaf of the universal curve. They introduce the standard torus action on $\P^n$ and using Bott's residue formula arrive at their computation; note that the end result is not a closed formula valid in all cases, but rather gives an algorithm which they implement to give concrete numbers in some cases of interest, including the famous quintic threefold.

The main goal in this paper is refine the integer-valued computations of Ellingsrud and Str{\o}mme to elements in the Grothendieck-Witt ring $\GW(k)$, where $k$ is the field of definition of the given hypersurface or complete intersection. For this, we need to look a bit more closely at the resolutions described in \cite{ES} and \cite{EPS} in order to construct a so-called {\em relative orientation} of the vector bundle whose sections describe the twisted cubic curves contained in a given hypersurface (the section depending on the selected hypersurface). Once we have this information, we apply the quadratically refined version of Bott's residue formula, found in \cite[Theorem 9.5]{LevineAB}. This allows to proceed essentially as in \cite{ES}, although one does need  to extend the torus action to an action of the normalizer $N_{\SL_2}$ of the torus in $\SL_2$ in order to retain the quadratic information, and in addition one needs to keep careful track of the orientation condition, which boils down to keeping careful track of signs.

In addition to the numerical condition on the degrees of the complete intersection $X$ that is needed in order for one to expect to have only finitely many twisted cubics on $X$,  the requirement of a relative orientation imposes an additional congruence. For instance, in the case of a hypersurface of degree $m$ in $\P^n$, the numerical condition is $3m+1=4n$, and for the relative orientation one needs $n$ to be even and $m\equiv1\mod 4$ (see Proposition~\ref{prop:relOrientationObservations}).  Just as in \cite{ES}, we do not have a closed form formula for the quadratic count, but we have a number of examples of concrete computer-assisted computations. 

We should add that the localization methods of \cite{LevineAB} destroy the two-torsion information in $\GW(k)$, so {\it a priori} we retain only the information given by the signature function for each real embedding of the field of definition $k$ (as well as the rank information given by the numerical count of Ellingsrud and Str{\o}mme); this does give a complete computation for $k=\R$. However, we show that the {\em problem} being posed lives over  $\Z[1/6]$, so one does have an answer in the form of an element of $\GW(\Z[1/6])$, which yields the computation over an arbitrary base-field $k$ of characteristic $\neq 2,3$ by base-extension (see \S\ref{sec:Integrality}, especially Corollary~\ref{cor:Integrality}). Noting that $\GW(\Z)\to \GW(\R)$ is an isomorphism, and that the cokernel of $\GW(\Z)\to \GW(\Z[1/6])$ is killed by passage to $\GW(k)$ for $\Z[1/6]\subset k$ if $2$ and $3$ are squares in $k$, we see that the signature and rank yield a complete answer for all such fields (for example,  for  the fields $\F_p$ with $p\equiv \pm 1\mod 24$, i.e., for roughly 1/4 of these fields). See Remark~\ref{rem:Refinement} for details.

We computed the count in all relatively oriented cases where $n\le 12$, feeding our computations into a Macaulay2 program. Here are the computations in tabular form.  \\[10pt]
\scalebox{0.92}{
\begin{tabular}{|c|c|r|r|}
\hline
$n$&degree(s)&signature&rank\\
\hline
\hline
4&(5)&765&317206375\\
\hline
5&(3,3)&90&6424326\\
\hline
10&(13)&768328170191602020&794950563369917462703511361114326425387076\\
\hline
11&(3,11)&4407109540744680&31190844968321382445502880736987040916\\
\hline
11&(5,9)&313563865853700&163485878349332902738690353538800900\\
\hline
11&(7,7)&136498002303600&31226586782010349970656128100205356\\
\hline
12&(3,3,9)&43033957366680&3550223653760462519107147253925204\\
\hline
12&(3,5,7)&5860412510400&67944157218032107464152121768900\\
\hline
12&(5,5,5)&1833366298500&6807595425960514917741859812500\\
\hline
\end{tabular}}
\ \\[2pt]
Note that among the above cases, only for $n=4, m=(5)$, $n=5, m=(3,3)$ is the variety Calabi-Yau.
 
The signature is what is computed in this paper and gives a signed count of the real twisted cubics in a real hypersurface or complete intersection. Each curve is counted with its ``local quadratic degree'', which we have not attempted to discuss here; note that in any case, the signed count gives a lower bound for the number of real twisted cubics, assuming the local degrees are all $\pm1$, which we expect to be the case for the general hypersurface or complete intersection. The rank is the classical count of the twisted cubics in a hypersurface or complete intersection over an algebraically closed field (the cases $(4,(5))$ and $(5, (3,3))$ are given in \cite{ES} and the other cases were computed by transforming the formulas of \cite{ES} into a Macaulay2 program. This program, and the one computing the signatures, were written by the second author). Putting these together gives the quadratic count $Q$, 
\[
Q:=s+\frac{r-s}{2}\cdot H, 
\]
valid for a field $k$ of characteristic $\neq 2,3$ for which both $2$ and $3$ are squares in $k$, 
where $s$ is the signature, $r$ is the rank, and $H$ is the hyperbolic form $H(x,y)=x^2-y^2$. For a general $k$ (of characteristic $\neq 2,3$), 
\[
Q=s+\frac{r-s}{2}\cdot H+\epsilon_1(\<2\>-1)+\epsilon_2(\<3\>-1)+\epsilon_3(\<6\>-1)
\]
with the $\epsilon_j$  constants depending only on $n, m_*$, and with $\epsilon_1, \epsilon_2\in \{0,1\}$,
$\epsilon_3\in \{0,1,2,3\}$ (see Lemma~\ref{lem:Coker}).

The first named author thanks John Christian Ottem for suggesting the quadratic count of twisted cubics as an application of the quadratic equivariant localization method, and Raman Parimala for her help with understanding aspects of the Grothendieck-Witt ring of $\Z[1/d]$.

Unless explicitly mentioned otherwise, we work over a field $k$ of characteristic $\neq 2,3$.

\section{Eagon-Northcott complexes and flat families}

For use throughout the paper, we recall the simplest case of the Eagon-Northcott complex. Let $R$ be a commutative ring and let $A=(a_{ij})$ be a $3\times2$ matrix with the $a_{ij}\in R$. For $1\le i<j\le 3$, let $e_{ij}(A)$ denote the determinant of the $2\times2$ submatrix of $A$ with rows $i,j$ and let $(A)\subset R$ denote the ideal $(e_{12}(A), e_{13}(A), e_{23}(A))$. Let $B:R^3\to R$ be the map with matrix $(e_{23}(A), -e_{13}(A), e_{12}(A))$. We identify $B$ with $\bigwedge^2A^t$, giving us the Eagon-Northcott complex
\begin{equation}\label{eqn:EN1}
0\to R^2\xrightarrow{A} R^3\xrightarrow{\bigwedge^2A^t}R\to R/(A)\to 0
\end{equation}

\begin{theorem}[Eagon-Northcott]\label{thm:EN}  Suppose that $R$ is a Cohen-Macaulay ring and $(A)$ is not the unit ideal. \\[5pt]
1. Each minimal prime ideal of $(A)$ has height $\le 2$.\\[2pt]
2. The sequence \eqref{eqn:EN1} is exact if and only if each minimal prime ideal of $(A)$ has height $=2$.
\end{theorem}
This  follows from \cite[Theorem1, Theorem 2]{EN} in the special case $r=3, s=2$ and taking $R$ to be Cohen-Macaulay to identify the grade of an ideal with the minimum height of an associated prime.

\begin{definition} Let $R$ be a ring and consider the polynomial ring $R[X_0, X_1, X_2, X_3]$. Let 
\[
A=a_1X_1+a_2X_2+a_3X_3,\ B=b_2X_2+b_3X_3,\ C=c_2X_2+c_3X_3
\]
be linear forms in $R[X_0, X_1, X_2, X_3]$ and let  $q=AX_1^2+BX_1X_2+CX_2^2$. Let 
\[
\alpha_t(q)=\begin{pmatrix}tC&tB-X_0\\X_0&tA\\-X_1&X_2\end{pmatrix}
\]
We say that $q$ is {\em non-degenerate} if $(a_1, a_2, a_3, b_2, b_3, c_2, c_3)$ is the unit ideal in $R$, or equivalently, $q$ is non-zero modulo each maximal ideal in $R$.

Let $J(q)\subset R[t][X_0, X_1, X_2, X_3]$ be the ideal generated by the determinants of the $2\times2$ submatrices of $\alpha_t(q)$ and let $I(q)=(J(q), q)\subset R[t][X_0, X_1, X_2, X_3]$.  We let $\sJ(q)\subset \sI(q)\subset \sO_{\P^3_{R[t]}}$ denote the corresponding ideal sheaves

Let $C(q)\subset \P^3_{R[t]}$ be the subscheme defined by the  ideal sheaf
$\sI(q)$ and let $C'(q)\subset \P^3_{R[t]}$ be the subscheme defined by the  ideal sheaf
$\sJ(q)$. 
\end{definition}

\begin{lemma}\label{lem:TorsionFree} Suppose that $R$ is locally a UFD and $q$ is non-degenerate. Then  \\[5pt]
1. $C'(q)\cap \P^3_{R[t, t^{-1}]}=C(q)\cap \P^3_{R[t, t^{-1}]}$\\[2pt]
2. $C'(q)\subset \P^3_{R[t]}$ has pure codimension two, and has relative dimension one over $\Spec R[t,t^{-1}]$.\\[2pt]
3. The sheaf $\sO_{C(q)}$ is $t$-torsion free.
\end{lemma}

\begin{proof} The statement is local in $R$, so we may assume that $R$ is a local UFD.

 We have the exact sequence
\[
0\to \sI(q)/\sJ(q)\to \sO_{C'(q)}\to \sO_{C(q)}\to 0
\]
We note that $\sI(q)/\sJ(q)$ is generated by the image $[q]$ of $q$, that is, we have the surjection
\[
p_q:\sO_{\P^3_{R[t]}}(-3)\to \sI(q)/\sJ(q)
\]
sending $1$ to $[q]$.
We claim that the kernel of $p_q$ is the ideal generated by $(X_0, t)$. 

To see this, we have the identity
\[
tq=X_1e_{23}(\alpha_t(q))+X_2e_{13}(\alpha_t(q))\in J(q),
\]
so $t[q]=0$. This shows that $t\sI(q)\subset \sJ(q)$, so $(\sI/\sJ)/t(\sI/\sJ)=\sI/\sJ$ and also proves (1).  Moreover,  the kernel of $p_q$ is the same as the kernel of the induced map 
\[
\tilde{p}_q:\sO_{\P^3_{R[t]}}(-3)\to  \sO_{C'(q)}=\sO_{\P^3_{R[t]}}/\sJ(q).
\]

Let $\bar{\sJ}(q)$ be the image of $\sJ(q)$ in $\sO_{\P^3_{R[t]}}/(t)=\sO_{\P^3_R}$, that is
\[
\bar{\sJ}(q)=(X_0X_2, X_0X_1, X_0^2)\sO_{\P^3_R}.
\]
We have
\[
 \sO_{C'(q)}/t \sO_{C'(q)}=\sO_{\P^3_{R[t]}}/(t, \sJ(q))=\sO_{\P^3_R}/\bar{\sJ}(q).
 \]
 Letting $\bar{q}$ denote the image of $q$ in $\sO_{\P^3_R}/\bar{\sJ}(q)$, we claim that the annihilator of $\bar{q}$ in $\sO_{\P^3_R}$ is $(X_0)$. To see this, if we take a homogeneous $\alpha\in R[X_0,\ldots, X_3]$ with $\alpha\cdot q\in (X_0X_2, X_0X_1, X_0^2)$, then
 \[
 \alpha\cdot q=X_0(a\cdot X_0+b\cdot X_1+c\cdot X_2)
 \]
 for some $a,b,c\in  R[X_0,\ldots, X_3]$. Since $q$ is non-degenerate and $R$ is a UFD, $X_0$ divides $\alpha$, so $\alpha$ is in $(X_0)$. Conversely, a direct computation shows
 \[
 X_0\cdot q=(BX_1+CX_2)\cdot e_{23}(\alpha_t(q))-AX_1\cdot e_{13}(\alpha_t(q))\in J(q)
 \]
 This shows that $\Ann_{\sO_{\P^3_R}}(\bar{q})=(X_0)\sO_{\P^3_R}$ and in addition that 
 \[
 \Ann_{\sO_{\P^3_{R[t]}}}([q])\supset (X_0, t)\sO_{\P^3_{R[t]}}= \Ann_{\sO_{\P^3_{R[t]}}}(\bar{q})\supset
 \Ann_{\sO_{\P^3_{R[t]}}}([q])
 \]
so 
 $ \Ann_{\sO_{\P^3_{R[t]}}}([q])= (X_0, t)\sO_{\P^3_{R[t]}}$, as claimed.

 This also shows that the map
 \[
\sI(q)/\sJ(q)= \sI(q)/\sJ(q)\otimes_{R[t]}R\to \sO_{C'(q)}\otimes_{R[t]}R= \sO_{\P^3_R}/\bar{\sJ}(q)
\]
induced by the inclusion $\sI(q)/\sJ(q)\subset \sO_{C'(q)}$ is injective, since $[q]$ has the same annihilator in $\sO_{\P^3_{R[t]}}$ as does $\bar{q}$, and gives us the presentation of $\sI(q)/\sJ(q)$ via the Koszul complex
\begin{multline}\label{multline:Koszul}
0\to \sO_{\P^3_{R[t]}}(-4)\xrightarrow{\begin{pmatrix}X_0\\-t\end{pmatrix}}
\sO_{\P^3_{R[t]}}(-3)\oplus \sO_{\P^3_{R[t]}}(-4)\\\xrightarrow{\displaystyle(t, X_0)}  \sO_{\P^3_{R[t]}}(-3)\to \sI(q)/\sJ(q)\to 0
\end{multline}

We now prove (2).   Let $\mathfrak{m}$ be a maximal ideal of $R[t, t^{-1}]$. We have seen that $tq$ is a section of $\sJ(3)$. Moreover,  $e_{23}(\alpha_t(q))=X_0X_2+tAX_1$ is irreducible of degree two in $R[t, t^{-1}]/\mathfrak{m}[X_0, X_1, X_2, X_3]$ if  $(a_1, a_3)\not\subset\mathfrak{m}$, so assuming $(a_1, a_3)\not\subset\mathfrak{m}$, 
 $C'(q)\cap \P^3_{R[t, t^{-1}]/\mathfrak{m}}$ has codimension $\ge2$ in
$\P^3_{R[t, t^{-1}]/\mathfrak{m}}$. If $(a_1, a_3)\subset\mathfrak{m}$, then 
$e_{23}(\alpha_t(q))\equiv X_2(X_0-ta_2X_1)\mod\mathfrak{m}$ and since $q$ is non-degenerate the subscheme defined by $(X_0-ta_2X_1, tq)$ has codimension two in $\P^3_{R[t, t^{-1}]/\mathfrak{m}}$. For any $\mathfrak{m}$, the ideal $(X_2, e_{13}(\alpha_t(q)))\subset  R[t, t^{-1}]/\mathfrak{m}[X_0,\ldots, X_2]$ is $(X_2, (tb_3X_3-X_0)X_1)$, which also defines a codimension two closed subscheme of  $\P^3_{R[t, t^{-1}]/\mathfrak{m}}$, so $C'(q)\cap \P^3_{R[t, t^{-1}]}$ has codimension $\ge2$ in 
$\P^3_{R[t, t^{-1}]}$. Thus we have shown that $C'(q)$ has relative  dimension $\le 1$ over $\Spec R[t,t^{-1}]$.

 Working modulo $t$, the ideal $\bar{\sJ}(q)$ defines a closed subscheme of $\P^3_R$ of codimension $\ge1$, so $C'(q)$ has codimension $\ge2$ in $\P^3_{R[t]}$, and by Theorem~\ref{thm:EN}, $C'(q)$ has pure codimension two in  $\P^3_{R[t]}$. With what we have shown above, this implies that  $C'(q)$ has pure relative  dimension one over $\Spec R[t,t^{-1}]$.

Using Theorem~\ref{thm:EN} again, we see that the Eagon-Northcott complex
\[
0\to  \sO_{\P^3_{R[t]}}(-3)^2\xrightarrow{\alpha_t(q)}
\sO_{\P^3_{R[t]}}(-2)^3\xrightarrow{\bigwedge^2\alpha_t(q)^t}  \sO_{\P^3_{R[t]}} \to \sO_{C'(q)}\to0
\]
is exact. We map the Koszul complex \eqref{multline:Koszul} to the Eagon-Northcott complex over the inclusion $\sI(q)/\sJ(q)\to \sO_{\P^3_{R[t]}}$ by the 
maps
\begin{align*}
&\phi_0:\sO_{\P^3_{R[t]}}(-3)\to  \sO_{\P^3_{R[t]}}\\
&\phi_1:\sO_{\P^3_{R[t]}}(-3)\oplus \sO_{\P^3_{R[t]}}(-4)\to \sO_{\P^3_{R[t]}}(-2)^3\\
&\phi_2: \sO_{\P^3_{R[t]}}(-4)\to  \sO_{\P^3_{R[t]}}(-3)^2
\end{align*}
defined by  $\phi_0:=\times q$, 
\[
\phi_1:=\begin{pmatrix}X_1&BX_1+CX_2\\-X_2&AX_1\\0&0\end{pmatrix}
\]
and 
\[
\phi_2:=\begin{pmatrix}-X_2\\-X_1\end{pmatrix}.
\]
This induces the map $\bar\phi$ of complexes after applying $-\otimes_{R[t]}R$.

A direct computation using these resolutions shows that $\Tor^{R[t]}_1(\sI(q)/\sJ(q), R)$ is isomorphic to $\sO_{\P^3_R}(-3)/(X_0)\sO_{\P^3_R}(-4)$, with  generator
\[
v:=\begin{pmatrix}1\\0\end{pmatrix}:  \sO_{\P^3_R}(-3)\to \sO_{\P^3_{R}}(-3)\oplus \sO_{\P^3_{R}}(-4)
\]
and that $\bar{\phi}_1(v)$ is a generator for $\Tor^{R[t]}_1(\sO_{C'(q)}, R)$. Thus, in the long exact $\Tor$-sequence
\begin{multline*}
\ldots\to \Tor^{R[t]}_1(\sI(q)/\sJ(q), R)\to \Tor^{R[t]}_1(\sO_{C'(q)}, R)\to \Tor^{R[t]}_1(\sO_{C(q)}, R)\\\to \sI(q)/\sJ(q)\otimes_{R[t]}R\to \sO_{C'(q)}\otimes_{R[t]}R\to \sO_{C(q)}\otimes_{R[t]}R\to 0
\end{multline*}
the map $\Tor^{R[t]}_1(\sI(q)/\sJ(q), R)\to \Tor^{R[t]}_1(\sO_{C'(q)}, R)$ is surjective and the map 
$ \sI(q)/\sJ(q)\otimes_{R[t]}R=\sI(q)/\sJ(q)\to \sO_{C'(q)}\otimes_{R[t]}R$ is injective, so
\[
\Tor^{R[t]}_1(\sO_{C(q)}, R)=0.
\]  
Thus, the sheaf $\sO_{C(q)}$ is $t$-torsion free, proving (3). 
\end{proof}

\begin{proposition}  Suppose that $q$ is non-degenerate, then  $C(q)$ is flat of relative dimension one over $\Spec R[t]$.
\end{proposition}

\begin{proof} The statement is local in $R$, so we may assume that $R$ is local.  To prove the result for $R$,  suffices to prove the result for some ring $S$ and a non-degenerate $q_S=A_SX_1^2+B_XX_1X_2+C_SX_2^2\in S[X_1,X_2, X_3]$ such that there is a ring homomorphism $\psi:S\to R$ sending $q_S$ to $q$. As $q=AX_1^2+BX_1X_2+CX_2^2$ with $A=a_1X_1+a_2X_2+a_3X_3$, $B=b_2X_2+b_3X_3$, $C=c_2X_2+c_3X_3$, $a_i, b_i, c_i\in R$,  we may replace $R$ with a suitable localization $S$ of a polynomial ring over $\Z$,  and suitable element $q_S$. Thus, we may assume from the start that $R$ is a regular local ring, in particular, a UFD. 

Let $C_0(q)\subset \P^3_R$ denote the closed subscheme with ideal $\sI(q)_0:=(\alpha_0(q))+(q)\subset R[X_0, X_1, X_2, X_3]$ and let  $\bar{C}_0(q)$ be the closed subscheme with ideal sheaf $\sqrt{\sI(q)_0}=(X_0, q)\sO_{\P^3_R}$. As $q$ is non-degenerate, 
$\bar{C}_0(q)$ is a  complete intersection subscheme of $\P^3_R$, equi-dimensional of relative dimension one over $R$, and hence $\bar{C}_0(q)$ is flat over $R$. Note that the quotient $\sqrt{\sI(q)_0}/ \sI(q)_0$ is generated as $\sO_{\P^3_R}$-module by the image $[X_0]$ of $X_0$. For each $r\in R\setminus \{0\}$, $r[X_0]\neq0$, as $\sI(q)_0\subset (X_0, X_1, X_2)^2\sO_{\P^3_R}$ but 
$rX_0$ is not in $(X_0, X_1, X_2)^2\sO_{\P^3_R}$. Moreover, $(X_0, X_1, X_2)\cdot\sqrt{\sI(q)_0}/ \sI(q)_0=0$. This implies that $\sqrt{\sI(q)_0}/ \sI(q)_0=i_*(\sO_{\Spec R})$, where $i:\Spec R\to \sO_{\P^3_R}$ is the constant section with value $(0,0,0,1)$. From the exact sequence
\begin{equation}\label{eqn:ExactCubicSeq}
0\to \sqrt{\sI(q)_0}/ \sI(q)_0\to \sO_{C_0(q)}\to \sO_{\bar{C}_0(q)}\to 0
\end{equation}
we see that $C_0(q)$ is flat of relative dimension one over $R$. 
 
By Lemma~\ref{lem:TorsionFree}(1, 2) $C(q)\cap \P^3_{R[t,t^{-1}]}$ has relative dimension one over $R[t,t^{-1}]$, so $C(q)\cap \P^3_{R[t,t^{-1}]}$ has  codimension two in  $\P^3_{R[t,t^{-1}]}$, and this remains the case after  base-change by $R[t,t^{-1}]\to R'$ for any commutative noetherian ring $R'$. By Theorem~\ref{thm:EN}, this implies that the Eagon-Northcott complex
\[
0\to  \sO_{\P^3_{R[t]}}(-3)^2\xrightarrow{\alpha_t(q)}
\sO_{\P^3_{R[t]}}(-2)^3\xrightarrow{\bigwedge\alpha_t(q)^t}  \sO_{\P^3_{R[t]}} \to \sO_{C(q)}\to0
\]
is exact after base-change to $R'$ for any ring extension $R[t,t^{-1}]\to R'$, and thus $C(q)\cap \P^3_{R[t,t^{-1}]}$ is flat over $\Spec R[t,t^{-1}]$. 

Moreover, by Lemma~\ref{lem:TorsionFree}(3), we have
\[
\Tor_i^{R[t]}(\sO_{C(q)}, R[t]/(t))=0
\]
for all $i>0$, and thus for any prime ideal $\mathfrak{p}\subset R[t]$ containing $(t)$, we have
\[
\Tor_i^{R[t]}(\sO_{C(q)}, R[t]/\mathfrak{p})=\Tor_i^{R[t]/(t)}(\sO_{C_0(q)}, R[t]/\mathfrak{p})
\]
But as  $C_0(q)$ is flat over $R=R[t]/(t)$, it follows that
\[
\Tor_i^{R[t]}(\sO_{C(q)}, R[t]/\mathfrak{p})=0
\]
for $i>0$. Since we have already seen that $C(q)$ is flat over $R[t,t^{-1}]$, it follows that $C(q)$ is flat over $R[t]$. 
\end{proof} 

\section{Cubic rational curves in $\P^n$}
We fix a base-field $k$ of characteristic $\neq 2, 3$. Unless otherwise noted, all schemes are finite type separated $k$-schemes. We write  $\P^n$ for $\P^n_k$.

Let $K\supset k$ be an extension field and let $C\subset \P^n_K$ be a closed subscheme. The {\em linear span} of $C$ is the smallest linear subspace $L_C$ of $\P^n_K$ containing $C$. In fact, $L_C$ is always a $K$-subspace of $\P^n_K$, as  the ideal of $L_C$ is generated by the linear forms in the $K$-vector space $H^0(\P^n_K, \sI_C(1))$. 

\begin{lemma}\label{lem:linearspan} Let $\sO$ be a dvr with quotient field $K$ and residue field $\kappa$, and let $C\subset \P^n_\sO$ be a closed subscheme, flat over $\sO$. Let $L_{C_K}$ be the linear span of the generic fiber $C_K$ of $C$. \\[5pt]
1. There is a projective space bundle $L\subset \P^n_\sO$, smooth over $\sO$ with generic fiber $L_{C_K}$ and with $C\subset L$.\\[2pt]
2.  Suppose $L_{C_K}$  has dimension $m$ over $K$. Then there is an element $A\in \GL_{n+1}(\sO)$ such that $A\cdot C\subset \P^m_\sO\subset \P^n_\sO$, where $\P^m_\sO$ is the linear subspace defined by $X_{m+1}=\ldots=X_n=0$.
\end{lemma}

\begin{proof} (1) The linear space $L_{C_K}$ being defined over $K$, $L_{C_K}$ gives a $K$-point of a Grassmannian $\Gr(m+1, n+1)$, where $m$ is the dimension of $L_{C_K}$. Since $\Gr(m+1, n+1)$ is proper over $k$, this $K$-point extends to a morphism $f:\Spec \sO\to \Gr(m+1, n+1)$; pulling back the universal $\P^m$-sub-bundle of $\P^n\times\Gr(m+1, n+1)$ via $f$ gives us a $\P^m$-bundle $L\subset \P^n_\sO$ with $C_K\subset L_K$. Since $\P^n_\sO$ is separated over $\sO$, this implies that $C\subset L$.

For (2), the subspace $L$ is defined by an $m+1\times n+1$ matrix $M\in M_{m+1, n+1}(\sO)$, with the rows of $M$ giving $\sO$-points of $\P^n_\sO$ that span $L$. From the structure theory of matrices over the dvr $\sO$, there are invertible matrices $A\in \GL_{n+1}(\sO)$ and $B\in \GL_{m+1}(\sO)$ such that
\[
B\cdot M\cdot A=(I_{m+1}, 0_{m+1\times n-m}).
\]
\end{proof}

Let $\Hilb^{3m+1}(\P^n)$ be the Hilbert scheme parametrizing flat families of closed subschemes of $\P^n$ with Hilbert polynomial $3m+1$, and with universal curve $\sC^{3m+1, n}\to \Hilb^{3m+1}(\P^n)$. For $x\in \Hilb^{3m+1}(\P^n)$, we write $C_x\subset \P^n_{k(x)}$ for the corresponding curve.

In $\Hilb^{3m+1}(\P^n)$ we have the open subscheme $H_n^{sm}$ of smooth cubic rational curves, classically known as {\em twisted cubics}. Let  $H_n$ denote the closure   of $H_n^{sm}$ in  $\Hilb^{3m+1}(\P^n)$, and let $\sC_n\subset H_n\times \P^n$ be the corresponding universal curve, with projection $p:=p_1:\sC_n\to H_n$. For $x\in H_n$, we call the corresponding  curve $C_x\subset \P^n_{k(x)}$  a {\em rational cubic curve in $\P^n$}.  It is well-known that $H_n^{sm}$ is smooth and irreducible, so $H_n$ is an integral $k$-scheme. 

If we wish to work over another base-scheme, we write $H_n/S$ for the closure in $\Hilb^{3m+1}(\P^n_S)$ of the open subscheme  $H^{sm}_{n, S}\subset \Hilb^{3m+1}(\P^n_S)$ parametrizing flat families of twisted cubics.

\begin{lemma}\label{lem:3Space} 1. For each $x\in H_n$, $C_x$ has linear span of dimension 3. \\[2pt]
2. After a $k(x)$-linear change of coordinates, $C_x$ is one of two types:\\[5pt]
a. $C_x$ is Cohen-Macaulay. Passing to the algebraic closure $\overline{k(x)}$ of $k(x)$, there is a $3\times 2$ matrix  $(L_{ij})$ of linear forms  $L_{ij}\in \overline{k(x)}[X_0,\ldots, X_3]$ such that the ideal of $C_x\times_{k(x)}\overline{k(x)}$ is generated by the linear forms $X_4,\ldots, X_n$ and the determinants of the $2\times2$ submatrices of $(L_{ij})$. Moreover these three determinants are $\overline{k(x)}$-linearly independent elements of $\overline{k(x)}[X_0,\ldots, X_3]$.\\[2pt]
b. $C_x$ is not Cohen-Macaulay. There are linear forms $B, C\in \overline{k(x)}[X_2, X_3]$ and $A\in \overline{k(x)}[X_1, X_2, X_3]$, not all zero, such that the ideal of $C_x\times_{k(x)}\overline{k(x)}$ is generated by the linear forms $X_4,\ldots, X_n$, the quadratic terms $X_0^2, X_0X_1, X_0X_2$ and the cubic polynomial $q_C:=AX_1^2+BX_1X_2+CX_2^2$.\\[2pt]
3. Let $\sI_x\subset \sO_{\P^n_{k(x)}}$ denote the ideal sheaf of $C_x$. Then for $i>0$, 
 $H^i(\P^n_{k(x)}, \sI_x(m))=0$ for $m\ge0$ if $C_x$ is Cohen-Macaulay, and 
 $H^i(\P^n_{k(x)}, \sI_x(m))=0$ for $m\ge1$ if $C_x$ is not Cohen-Macaulay.
\end{lemma}

\begin{proof} Since a smooth cubic curve in a $\P^2$ has genus 1, the linear span of a smooth rational degree three curve $C\subset \P^n$ must have dimension at least three. The fact that $h^0(\P^1, \sO_{\P^1}(3))=4$   says that the kernel of the restriction map
\[
H^0(\P^n, \sO_{\P^n}(1))\to H^0(C, \sO_C(1))
\]
has codimension at most 4, so there are at least $n-3$ linearly independent linear forms vanishing on $C$. Thus the linear span of $C$ has dimension at most three, that is, the linear span of a  smooth rational degree three curve in $\P^n$ is a $\P^3$.

 As $C_x$ is the specialization over some dvr $\sO$ of a smooth rational cubic curve $C_K\subset \P^n_K$, it follows by Lemma~\ref{lem:linearspan} that, after a $k(x)$-linear change of coordinates, $C_x$ is contained in the linear subspace $\P^3\subset \P^n$ defined by $X_4=\ldots, X_n=0$ and $C_x$ is a $k(x)$ point of  the $H_3\subset H_n$ defined by this $\P^3$ in $\P^n$. This reduces the proof of (1) and (2) to the case $n=3$. 

We pass to the algebraic closure of $k(x)$ and write $C_x$ for $C_x\times_{k(x)}\overline{k(x)}$.  If $C_x$ is Cohen-Macaulay, the resolution of $\sI_{C_x}$ given in \cite[Lemma 1]{PS} is an Eagon-Northcott complex of the form
\[
0\to \sO_{\P^3}(-3)^2\xrightarrow{(L_{ij})}\sO_{\P^3}(-2)^3\xrightarrow{\bigwedge^2(L_{ij})^t}\sI\to 0
\]
with $(L_{ij})$ the matrix of linear forms mentioned in (2a) and with $\bigwedge^2(L_{ij})^t$ the $1\times 3$ matrix of determinants $(e_{23}(L_{ij}), -e_{13}(L_{ij}), e_{12}(L_{ij}))$. This takes care of the case (2a). \cite[Lemma 2]{PS} yields the assertion (2b), describing the non-Cohen-Macaulay curves. 

In the Cohen-Macaulay case, one shows using the Eagon-Northcott resolution  that $h^i(\P^3, \sI(m))=0$ for $m\ge0$, and $h^0(\P^3, \sI(1))=0$. Thus the linear span of $C_x$ is the entire $\P^3$. In the non-Cohen-Macaulay case, let $\sJ$ be the ideal $(X_0, q)\sO_{\P^3}$, so $\Spec \sO_{\P^3}/\sJ$ has linear span the plane $X_0=0$.  In the affine space $\A^3:=\P^3\setminus \{X_3=0\}$, $C_x$ has defining ideal $(x_0x_1, x_0x_2, x_0^2, q(x_0, x_1, x_2, 1))$, $x_i:=X_i/X_3$, so $X_0$ is not in $H^0(\P^3, \sI(1))$, hence $C_x$ has linear span $\P^3$. This proves (1).

For (3), we first consider the Cohen-Macaulay case. We may assume that $C_x$ is contained in the $\P^3$ defined by $X_4=\ldots=X_n=0$; let $\bar{\sI}$ denote the ideal sheaf of $C_x\subset \P^3$ and let $\sK=(X_4,\ldots, X_n)\sO_{\P^n}$. This gives the exact sequence
\[
0\to  \sK\to \sI\to i_*\bar\sI\to 0
\]
with $i:\P^3\to \P^n$ the inclusion. Using the Koszul resolution of $\sK$, one shows that $H^i(\P^n,\sK(m))=0$ for $i,m\ge0$ and the Eagon-Northcott resolution of $\bar\sI$ shows that $H^i(\P^n,\bar\sI(m))=0$ for $i,m\ge0$. In the non-Cohen-Macaulay case, a similar argument reduces us to the case $n=3$; changing notation, we may assume $k=k(x)$, $k$ is algebraically closed and $C_x$ has ideal $\sI$ as described in the proof of (2b). We have the exact sequence
\[
0\to \sI\to \sJ\to \sJ/\sI\to 0
\]
where  $\sJ=(X_0, q)\sO_{\P^3}$. Letting $i:\Spec k\to \P^3$ be the inclusion at the point $(0,0,0,1)$, we have $\sJ/\sI\cong i_*k(-1)$ with generator $X_0$. Thus the map
\[
H^0(\P^3, \sJ(1))\to H^0(\P^3, \sJ/\sI(1))
\]
is surjective. Clearly $H^i(\P^3, \sJ/\sI(m))=0$ for all $i>0$ and all $m$, and the Koszul resolution of $\sJ$ shows that $H^i(\P^3, \sJ(m))=0$ for all $i>0$ and all $m\ge0$. Thus $H^i(\P^3, \sI(m))=0$ for all $i>0$ and all $m\ge1$, completing the proof.
\end{proof}

\begin{proposition}\label{prop:DefSpace} $H_n$ is a smooth $k$-scheme of dimension $4n$ and the non-Cohen-Macaulay locus $H^\ncm_n$ is a smooth divisor. 
\end{proposition}

\begin{proof} For $n=3$, this follows from the main result \cite[\S 6, Theorem]{PS}. The proof of {\em loc. cit.} relies on the construction of the universal deformation space of the non-Cohen-Macaulay curve $C_0$ defined by the ideal $\sI:=(X_0X_1, X_0X_2, X_0^2, X_1^3)\sO_{\P^3}$. Let $P=k[X_0,\ldots, X_3]$ and let $I\subset P$ be the ideal  $(X_0X_1, X_0X_2, X_0^2, X_1^3)$. The argument is in two steps. First they explicitly compute the universal deformation space of $I$, considered as a homogeneous ideal in $P$. To do this, they construct an  explicit resolution 
\[
0\to F_3\to F_2\to F_1\to P\to P/I\to 0
\]
Using this resolution, they write down a $k$-basis of the 16-dimensional tangent space $T_{M_3',I}$, with the first 10 basis elements $\del/\del u_1,\ldots, \del/\del u_{10}$ corresponding to deformations that fix the flag $(X_0=X_1=X_2=0)\subset (X_0=0)$ and the remaining 6 ``trivial" vectors coming from the translations by $\GL_4$ acting on such flags. In the subspace spanned by the $\del/\del u_i$, they construct the versal deformation space of $I\subset P$ explicitly, and show this is a union of affine spaces $\A_k^6\cup \A^9_k$, intersecting along a linear subspace $\A^5_k\subset \A^6_k$,  with $\A_k^6\setminus\A_k^5$ is the locus of Cohen-Macaulay deformations. The universal deformation space is then the product of the versal space with the $\A^6$ corresponding to the 6 ``trivial'' deformations, embedded as a linear space in $\GL_4$ transversal to the isotropy group of the flag, and parametrizing the big affine cell in the corresponding flag variety.

The second step  relies on   \cite[\S 3, Comparison Theorem]{PS}. Using this result, they show that the map $M'_3=\A^{12}\cup \A^{15}\to \Hilb^{3m+1}(\P^3)$ induced by the universal family over $\A^{12}\cup \A^{15}$ is an open immersion. It is easy to show that the $\GL_4$ translates of the image covers   $H^\ncm_3$. This yields the main theorem \cite[\S 6, Theorem]{PS}, with $H_3$ corresponding to the component $\A^{12}$ and $H_3^\ncm$ corresponding to the intersection  $\A^{11}_k=\A^{12}_k\cap \A^{15}_k$ (and all their respective translates by $\GL_4$).

We consider the analogous situation for $C_0\subset \P^3\subset \P^n$, with ideal sheaf $\sI_n=(\sI, X_4,\ldots, X_n)\sO_{\P^n}\subset \sO_{\P^n}$ and homogeneous ideal $I_n:=(I, X_4,\ldots, X_n)\subset P_n:=k[X_0,\ldots, X_n]$. Here the $\P^3$ is defined by $X_4=\ldots=X_n=0$. We first compute the universal deformation space of $I_n\subset P_n$. Let $A_n:=P_n/I_n$. Note that $A_n\cong A_3$. We first consider the cotangent complex $\L_{A_n/P_n}$. We have the fundamental distinguished triangle (see \cite[Tag 08QR, (90.7.0.1), Proposition 90.7.4]{Stacks})
\[
\L_{P_3/P_n}\otimes_{P_3}^LA_3\to \L_{A_n/P_n}\to \L_{A_3/P_3}\to \L_{P_3/P_n}\otimes_{P_3}A_3[1],
\]
where we view $P_3$ as the quotient $P_n/(X_4,\ldots, X_n)$ . Since $P_3$ is defined by the complete intersection ideal $(X_4,\ldots, X_n)$, we have $\L_{P_3/P_n}=P_3(-1)^{n-3}[1]$ \cite[Tag 08SH, Lemma 90.14.2]{Stacks}, and the distinguished triangle yields the exact sequence
\begin{multline*}
0\to T_2(A_3/P_3, A_3)\to T_2(A_n/P_n,A_n)\to 0\\
\to T_1(A_3/P_3, A_3)\to T_1(A_n/P_n,A_n)\to
A_n(1)^{n-3}\to 0
\end{multline*}
We write $M_d$ for the degree $d$ graded component in a $\Z$-graded object $M$.

From the above sequence, we see that the obstruction space $T_2(A_n/P_n,A_n)_0$ for $I_n\subset P_n$  agrees with that for $I_3\subset P_3$. Since the curve $C_0$ has $\P^3$ as linear span, the surjection $P_n\to A_n$ sends the degree one component to the four dimensional linear span of $X_0,\ldots, X_3$ in $A_n$. We can split the surjection $T_1(A_n/P_n,A_n)_0\to [A_n]_1^{n-3}$ by realizing $[A_n]_1^{n-3}$ as the global sections of the normal bundle of $\P^3$ in $\P^n$. This is realized as the universal deformation space of $\P^3$ in $\P^n$ given by the standard affine cell in $\Gr(4, n+1)$, which in turn is given by an affine space $V\subset \GL_{n+1}$, $V=\A^{4(n-3)}$, with  $\GL_{n+1}$ acting  by the usual action on $\A^{n+1}$. The corresponding first order deformations of $I_n$ given by  the subspace of $\mathfrak{gl}_n$ corresponding to $V$  give a $k$-subspace of  $ T_1(A_n/P_n,A_n)_0$ mapping isomorphically to $[A_n(1)^{n-3}]_0=k^{4(n-3)}$. Since the two obstruction spaces agree, this implies that the universal deformation space of $I_n\subset P_n$ is $M'_3\times V$.

The map $f:M'_3\times V\to \Hilb^{3m+1}(\P^n)$ corresponding to the universal family thus sends $\A^{12}\times V$ to $H_n$ and $\A^{11}\times V$ to $H_n^\ncm$. The hypotheses for the Comparison Theorem in this case follow easily from the case $n=3$ by translating by $V$, so the map $f$ is an open immersion of $M'_3\times V$ onto an open neighborhood of $[C_0]\in \Hilb^{3m+1}(\P^n)$. The $\GL_{n+1}$ translates of $f(\A^{11}\times V)$ covers $H^\ncm_n$ and the  $\GL_{n+1}$ translates of $f(\A^{12}\times V)$ cover $H_n$, since each curve $C$ in $H_n$ is contained in some $\P^3$ in $\P^n$. Thus $H_n$ is smooth and $H_n^\ncm$ is a smooth divisor in $H_n$. 
\end{proof}

We return to our study of $H_n$ over a fixed field $k$ of characteristic $\neq 2,3$.

\begin{lemma} 1. There is a morphism $\Phi:H_n\to \Gr(4, n+1)$ with $\Phi(x)=\pi_{C_x}$ for all $x\in H_n$.\\[2pt]
2. $\Phi$ is a smooth $\GL_{n+1}$-equivariant morphism, making $H_n$ into a fiber bundle (for the \'etale topology) over 
$\Gr(4, n+1)$  with fiber $H_3$.
\end{lemma}

\begin{proof} (1) Let $\sC_n\subset H_n\times\P^n$ be the universal curve, with defining ideal sheaf $\sI_n\subset \sO_{H_n\times\P^n}$. By Lemma~\ref{lem:3Space}(3) and standard base-change theorems on coherent cohomology, for all $m\ge1$, $R^ip_{1*}(\sI_n(m))=0$ for $i>0$, $p_{1*}(\sI_n(m))$ is locally free, and the map
\[
p_{1*}(\sI_n(m))\otimes k(x)\to H^0(\P^3_{k(x)}, \sI_{C_x}(m))
\]
is an isomorphism for all $x\in H_n$.  In particular $p_{1*}(\sI_n(1))$ is a locally free subsheaf of $p_{1*}(\sO_{H_n\times\P^n})=\sO_{H_n}\otimes_kH^0(\P^n, \sO_{\P^n}(m))$  of rank $n+1-4$, with $p_{1*}(\sI_n(1))\otimes k(x)\to p_{1*}(\sO_{H_n\times\P^n})\otimes k(x)$ injective for all $x\in H_n$. Taking the dual surjection  defines the morphism $\Phi:H_n\to \Gr(4, n+1)$ with $\Phi(x)=\pi_{C_x}$ for all $x\in H_n$.

For (2), it is clear that $\Phi$ is $\GL_{n+1}$ equivariant, where $\GL_{n+1}$ acts on $H_n$ and $\Gr(4, n+1)$ through its standard linear action on $\A^{n+1}$. 

To see that $\Phi$ is a smooth morphism, we first check over the $H_n^\cm$, an open subscheme of $\Hilb^{3m+1}(\P^n)$. For $x\in H^\cm_n$, the tangent space is given by
\[
T_{H_n, x}=\Hom(\sI_{C_x}, \sO_{C_x}).
\]
We may suppose that $\pi_{C_x}$ is our standard $\P^3$ defined by $X_4=\ldots=X_n=0$; let $\bar{\sI}_{C_x}\subset \sO_{\P^3}$ be the ideal sheaf of $C_x\subset \P^3$. We thus have the exact sequence
\[
0\to \sO_{\P^n}(-1)^{n-3}\to \sI_{C_x}\to i_*\bar{\sI}_{C_x}\to 0
\]
where $i:\P^3\to \P^n$ is the inclusion.  Looking at the Eagon-Northcott resolution of $\bar{\sI}_{C_x}$, we see that $H^1(\P^n, i_*\bar{\sI}_{C_x}(m))=0$ for all $m\ge0$, giving the short exact sequence
\[
0\to \Hom_{\sO_{\P^3}}(\bar{\sI}_{C_x}, \sO_{C_x})\to 
 \Hom_{\sO_{\P^n}}(\bar{\sI}_{C_x}, \sO_{C_x})\to H^0(C_x, \sO_{C_x}(1))^{n-3}\to 0
\]
Identifying $\sHom(\sO_{\P^n}(-1)^{n-3}, \sO_{\P^3})$ with the normal sheaf $\sN_i$ of $i$ allows us to identify the term $H^0(C_x, \sO_{C_x}(1))^{n-3}$ with $H^0(C_x,  \sN_i\otimes \sO_{C_x})$. Also, since $\P^3$ is the linear span of $C_x$ and $H^1(C_x, \bar{\sI}_{C_x}(1))=0$, we see that the restriction map
\[
T_{\Gr(4, n+1), \P^3}=H^0(\P^3, \sN_i)\to H^0(C_x, \sN_i\otimes\sO_{C_x})
\]
is an isomorphism. This identifies the term $H^0(C_x, \sO_{C_x}(1))^{n-3}$ in the above exact sequence with $i_x^*\Phi^*T_{\Gr(4, n+1), \Phi(x)}$, where $i_x:x\to H_n$ is the inclusion. In other words, $d\Phi_x$ is surjective and hence $\Phi$ is smooth at $x$.

Note that $H_n^\ncm\subset H_n$ is a proper closed subscheme, stable under the $\GL_{n+1}$-action. Thus for $y\in \Gr(4, n+1)$,  $H_n^\cm\cap \Phi^{-1}(y)\subset H_n\cap \Phi^{-1}(y)$ is dense in $\Phi^{-1}(y)$. Moreover, $\Gr(4, ,n+1)$ and $H_n$ are both smooth, hence $H_n$ is Cohen-Macaulay and each maximal ideal $\mathfrak{m}_y\subset  \sO_{\Gr(4, n+1)}$ is generated by a regular sequence. Since $\Phi^{-1}(y)_{\red}$ is  isomorphic to $H_3\times_kk(y)$, the fiber $\Phi^{-1}(y)$ 
has the expected dimension, and thus each fiber $\Phi^{-1}(y)$ has no embedded components; this also implies that the map $\Phi$ is flat. As $\Phi^{-1}(y)$ contains a smooth open dense subscheme, and $\Phi^{-1}(y)^{\red}\cong H_3\times_kk(y)$ is smooth over $k(y)$,  it follows that $\Phi^{-1}(y)=\Phi^{-1}(y)^{\red}$ and is thus smooth over $k(y)$, hence $\Phi$ is a smooth morphism. 
\end{proof}

For $x\in H_n$ and $C=C_x$, we often write $\pi_C$ for the $\P^3\subset \P^n$ corresponding $\Phi(x)\in \Gr(4, n+1)$.

Let $\Pi\subset \Gr(4, n+1)\times\P^n$ be the universal $\P^3$-bundle, and let $\Pi_n:=H_n\times_{ \Gr(4, n+1)}\Pi$, with projection $\pi_1$ to $H_n$ and $p_2$ to $\Pi$; we let $\pi_2:\Pi_n\to \Gr(4, n+1)$ denote the composition $p_1\circ p_2$.  $\sC_n$ is thus a closed subscheme of $\Pi_n$. Via the second projection $p_2:\Pi_n\to \P^n$, we may form the twist of a coherent sheaf $\sF$ on $\Pi_n$ by $p_2^*\sO_{\P^n}(m)$, which we denote by $\sF(m)$. We similarly twist a sheaf $\sG$ on $\sC_n$ by $p_2^*\sO_{\P^n}(m)$, and denote this by $\sG(m)$.

We consider 
$ \Gr(4, n+1)$ as embedded in a projective space by the Pl\"ucker embedding, giving us the invertible sheaf $\sO_{\Gr(4,n+1)}(m)$ for $m\in \Z$; hopefully the two different notations for twisting will not cause confusion.

This gives us the diagram
\begin{equation}\label{eqn:MainDiag}
\xymatrix{
&\sC_n\ar@{^(->}[d]^i\ar[ddl]_p\ar[ddr]^{p_2}\\
&\Pi_n\ar[dl]^{\pi_1}\ar[d]^{p_2}\ar@/_10pt/[dd]_{\pi_2}\ar[dr]^{p_2}\\
H_n\ar[dr]_\Phi&\Pi\ar[d]^{p_1}\ar[r]_{p_2}&\P^n\\
&\Gr(4, n+1)
}
\end{equation}

 We let $\bar{\sI}_C\subset \sO_{\pi_C}$ denote the ideal sheaf of $C\subset \pi_C$,  and $\sI_C\subset \sO_{\P^n}$ the ideal sheaf of $C\subset \P^n$. Similarly, let $\bar{\sI}_{\sC_n}\subset \sO_{\Pi_n}$ and $\sI_{\sC_n}\subset \sO_{H_n\times\P^n}$ denote the respective ideal sheaves of $\sC_n$ in $\Pi_n$ and in $H_n\times\P^n$.

\section{Some sheaves and presentations}

 Let $\sE:=p_*(\bar\sI_{\sC_n}(2))$ and let 
 \[
 \beta:p^*\sE(-2)\to \bar\sI_{\sC_n}\subset \sO_{\Pi_n}
 \]
 be the canonical map. Let $\sF$ be the kernel of the multiplication map
 \[
 m:\sE\otimes p_*\sO_{\Pi_n}(1)\to p_*\sO_{\Pi_n}(3)
 \]
 giving the map
 \[
 \alpha:p^*\sF(-3)\to p^*\sE(-2).
 \]

\begin{lemma} 1.  For all $x\in H_n$, we have
\[
\dim_{k(x)} H^0(\pi_{C_x},  \bar{\sI}_{C_x}(2))=3
\]
and $H^i(\pi_{C_x}, \bar{\sI}_{C_x}(m))=0$ for $m\ge 1$, $i\ge1$. \\[2pt]
2.   $\sE$ is a rank three locally free sheaf on $H_n$.
\end{lemma} 

\begin{proof} (1) implies (2) by standard results on the higher direct images of flat sheaves.

To prove (1),  we may reduce to the case $n=3$, so for $x\in H^\cm_n$, $\sI_{C_x}$ is generated by the $2\times 2$ determinants of a $3\times 2$ matrix of linear forms $(L_{ij})$. Let $R$ be the graded ring $k(x)[X_0, X_1, X_2, X_3]$ and let $I\subset R$ be the homogeneous ideal generated by the $2\times 2$ determinants $e_{ij}(L_{ij})$.  We have the corresponding (graded) Eagon-Northcott complex  \eqref{eqn:EN1}
\[
0\to R(-3)^2\xrightarrow{(L_ij))}R(-2)^3\xrightarrow{\bigwedge^2(L_ij))^t} R\to R/I\to 0
\]
which is exact since $C_x$ has the expected codimension 2 in $\P^3$ (Theorem~\ref{thm:EN}). Moreover, Riemann-Roch for the cubic CM-curve $C_x$ implies
\[
\dim_{k(x)} H^0(C_x, \sO_{C_x}(m))=3m+1
\]
for all $m\ge0$. This agrees with the Hilbert polynomial of the free part of the Eagon-Northcott complex,  which implies that the map $H^0(\P^3, \sO_{\P^3}(m))\to H^0(C_x, \sO_{C_x}(m))$ is surjective for all $m$ and that $H^i(C_x, \sO_{C_x}(m))=0$ for all $m\ge0$, $i\ge1$. Via the exact sequence
\[
0\to \bar\sI_{C_x}\to \sO_{\P^3}\to \sO_{C_x}\to 0
\]
we see that 
\[
h^0(\P^3, \bar\sI_{C_x}(m))=\begin{pmatrix}m+3\\3\end{pmatrix}-3m-1
\]
and 
\[
H^i(\P^3, \bar\sI_{C_x}(m))=0
\]
for all $m\ge0$ and $i\ge1$. For $m=2$, this yields $h^0(\P^3, \bar\sI_{C_x}(2))=3$. 

 If $C_x$ is a non-CM curve, then we may consider $C_x$ as a curve in $\P^3$ and after a linear change of coordinates, $\bar{\sI}_{C_x}\subset \sO_{\P^3}$ has generators $X_0^2, X_0X_1, X_0X_2, q$, with $q$ of the form 
\[
q=AX_1^2+BX_1X_2+CX_2^2
\]
where $A,B, C$ are linear forms in $X_1, X_2, X_3$. Let $C^0_x$ be the cubic curve 
$V(X_0, q)$ and let $p:=(0,0,0,1)\in C^0_x$. Then we have the exact sequence
\[
0\to \bar\sI_{C_x}\to \bar\sI_{C^0_x}\to i_{p^*}(k(p ))\to 0
\]
and $C^0_x$ is cubic plane curve. Also, the section $X_0$ of $\bar\sI_{C^0_x}(1)$ generates the term $i_{p^*}(k(p ))$.  

The ideal  $\bar\sI_{C^0_x}$ of $C^0_x$ is a complete intersection of degrees $1$ and $3$, so we have the exact sequence
\[
0\to \sO_{\P^3}(-4)\to \sO_{\P^3}(-1)\oplus \sO_{\P^3}(-3)\to \bar\sI_{C^0_x}\to0
\]
Thus $H^i(\P^3, \bar\sI_{C^0_x}(m))=0$ for $m\ge1$, $i\ge1$. Since $H^0(\bar\sI_{C^0_x}(1))\to k(p )$ is surjective, it follows that $H^i(\P^3, \bar\sI_{C_x}(m))=0$ for all $m\ge1$ and
\[
h^0(\P^3,   \bar\sI_{C_x}(m))=h^0(\P^3,   \bar\sI_{C^0_x}(m))-1=
\frac{m^3+6m^2-7m+6}{6}-1
\]
For $m=2$, this yields $h^0(\P^3,   \bar\sI_{C_x}(2))=3$. 
 \end{proof}
 
\section{Local descriptions}
We have the locally free sheaf $\sE:=p_*(\bar\sI_{\sC_n}(2))$ on $H_n$, with
\[
\sE\otimes k(x)=H^0(\pi_{C_x}, \bar{\sI}_{C_x}(2))
\]
We similarly have the locally free sheaf $\sV:=p_{1*}\sO_\Pi(2)$ of rank $10=h^0(\P^3, \sO_{\P^3}(2))$ on $\Gr(4, n+1)$. This gives us the Grassmannian bundle
\[
r:\Gr(3, \sV)\to \Gr(4, n+1)
\]
The inclusions $H^0(\pi_{C_x}, \bar{\sI}_{C_x}(2))\hookrightarrow H^0(\pi_{C_x}, \sO_{\pi_{C_x}}(2))$ define a a vector-bundle inclusion $\sE\hookrightarrow \Phi^*\sV$, which induces the (projective) morphism
\[
q:H_n\to \Gr(3,\sV)
\]
over $\Gr(4, n+1)$. This enlarges our diagram \eqref{eqn:MainDiag} to the diagram
\[
\xymatrix{
&\sC_n\ar[dl]_p\ar@{^(->}^i[d]\ar[dr]^{p_2}\\
H_n\ar[dr]_\Phi\ar[d]_q &  \Pi_n\ar[d]_{\pi_2}\ar[l]_{\pi_1} \ar[r]^{p_2}&\P^n\\
 \Gr(3, \sV)\ar[r]_r&\Gr(4, n+1)
}
\]

Letting  $E_{3, \sV}$ be the universal bundle on $\Gr(3, \sV)$  we have a canonical isomorphism $q^*E_{3, \sV}\cong \sE$. Let $X\subset \Gr(3, \sV)$ denote the image of $q$ and let $F\subset X$ denote the image of $H_n^\ncm$. It is shown in \cite[\S1]{EPS} that the restriction of $q$ to $H_n^\cm$ defines an isomorphism $q^\cm:H_n^\cm\to X\setminus F$; we henceforth identify $H_n^\cm$ and $X\setminus F$ via $q$. 

We recall another construction of $X$, given in \cite{EPS} in the case $n=3$; the general case is just formed by replacing the point $\Gr(4, 4)$ with our parameter space $\Gr(4, n+1)$.  Consider the sheaf $\sM:=\sHom(\sO_{\Gr(4, n+1)}^2, \sO_{\Gr(4, n+1)}^3\otimes \sV)$ and the corresponding vector bundle $M$ over  $\Gr(4, n+1)$, with fiber over $y\in  \Gr(4, n+1)$ the $3\times 2$ matrices of linear forms on $\Pi_y\cong \P^3$.  For $(L_{ij})$ in a fiber $M_y$, we have the corresponding subspace  of $\sV\otimes k(x)$ spanned by the determinants of the $2\times 2$-submatrices of $(L_{ij})$, $\<e_{12}( L_{ij}), e_{13}(L_{ij}), e_{23}(L_{ij})\>$. We let $\sU\subset M$ be the    subscheme parametrizing those $(L_{ij})$ such that   $\<e_{12}( L_{ij}), e_{13}(L_{ij}), e_{23}(L_{ij})\>$ has the maximal rank three. Sending $(L_{ij})\in \sU$ to $\<e_{12}( L_{ij}), e_{13}(L_{ij}), e_{23}(L_{ij})\>$ gives a morphism
\[
\psi: \sU\to \Gr(3, \sV).
\]
This gives us the diagram
\begin{equation}\label{eqn:MainDiag2}
\xymatrix{
&&\sC_n\ar[dl]_p\ar@{^(->}[d]^i\ar[dr]^{p_2}\\
&H_n\ar[dr]^-\Phi\ar[d]_q &  \Pi_n\ar[d]_{\pi_2} \ar[l]_{\pi_1}\ar[r]^{p_2}&\P^n\\
&X\ar@{^(->}[d]_{i_X}\ar[r]^-{r_X}&\Gr(4, n+1)\\
  \sU\ar[r]_-{\psi}&\Gr(3, \sV)\ar[ru]_r
}
\end{equation}
with $r_X$ the restriction of $r$.

\begin{theorem} 1. $\sU$ is an open subscheme of $M$.\\[2pt]
2. We have the the diagonal inclusion $\G_m\to \GL_3\times \GL_2$, sending $t$ to the pair of diagonal matrices with entries $t$; let $\Gamma$ be the quotient group-scheme  over $\Gr(4, n+1)$,
\[
\Gamma:=(\G_m\backslash\GL_3\times \GL_2)/\Gr(4, n+1).
\]
 Let  $\Gamma$ act  on $\sU$ by 
\[
(g, h)\cdot (L_{ij}):=g\circ (L_{ij})\circ h^{-1}.
\]
Then $\sU$ admits a smooth projective geometric quotient  $\Gamma\backslash \sU$, with quotient map $Q:\sU\to \Gamma\backslash \sU$ making $\sU$ a  $\Gamma$-torsor over $\Gamma\backslash \sU$.\\[2pt]
3. The map $\psi:\sU\to \Gr(3, \sV)$ has image $X$ and descends to an isomorphism
\[
\bar\psi:\Gamma\backslash \sU\to X.
\]
\end{theorem}

\begin{proof} This result is proven in \cite[Proposition 1, Lemma 2]{EPS} for $k$ an algebraically closed field of characteristic zero; the case of $k$ an arbitrary characteristic zero field follows by faithfully flat descent. 

It is not hard to show that for all $k$,  the $\Gamma$-action is free and proper. The arguments 
of {\it loc.\,cit.}, with Mumford's GIT extended by the work of Seshadri \cite{S1, S2} give the proof in positive characteristic.
\end{proof}

This gives us the diagram
\begin{equation}\label{eqn:MainDiag3}
\xymatrix{
&&\sC_n\ar[dl]_p\ar@{^(->}[d]^i\ar[dr]^{p_2}\\
&H_n\ar[dr]^-\Phi\ar[d]_q &  \Pi_n\ar[d]_{\pi_2} \ar[l]_{\pi_1}\ar[r]^{p_2}&\P^n\\
\Gamma\backslash \sU\ar@{=}[r]&X\ar@{^(->}[d]_{i_X}\ar[r]^-{r_X}&\Gr(4, n+1)\\
 \sU\ar[r]_-{\psi}\ar[u]^Q&\Gr(3, \sV)\ar[ru]_r
}
\end{equation}
 
Via this description of $X$ as $\Gamma\backslash \sU$, we may define  locally free sheaves $\sE_X, \sF_X$ on $X$ of rank three and two, respectively, as follows. We have the action of $\GL_3\times \GL_2$ on  $\sO_\sU^3$ and $\sO_\sU^2$ over the action on $\sU$ by letting a pair $(g, h)$ act on 
$\sO_\sU^3$  by
\[
(g,h)(v) :=\frac{\det h}{\det g}\cdot g(v), 
\]
and $\sO_\sU^2$ by
\[
(g,h)(w) :=\frac{\det h}{\det g}\cdot h(w).
\]
This descends to an action of $\Gamma$ on both $\sO_\sU^3$ and $\sO_\sU^2$ over $\sU$, giving the  locally free sheaves $\sE_X$ and $\sF_X$ on $X$ by descent.

Let $\Pi_\sU:=\sU\times_{\Gr(4, n+1)}\Pi$, with projections $q_1:\Pi_\sU\to \sU$, $q_2:\Pi_\sU\to \Pi_n$. Via $q_2$, we may twist a coherent sheaf $\sG$ on $\Pi_\sU$, setting $\sG(m):=\sG\otimes q_2^*\sO_{\Pi_n}(m)$.  Let 
\[
\sL:\sO_{\Pi_\sU}(-3)^2\to \sO_{\Pi_\sU}(-2)^3
\]
be the map with value $(L_{ij})$ at $(L_{ij})\in \sU$ and let
$\bigwedge^2\sL^t:\sO_{\Pi_\sU}(-2)^3\to \sO_{\Pi_\sU}$ be the map with value $(e_{23}(L_{ij}), -e_{13}(L_{ij}), e_{12}(L_{ij}))$ at $(L_{ij})\in \sU$. Here as before $e_{ij}(L_{ij})$ is the determinant of the $2\times 2$ minor of $(L_{ij})$ with rows $i,j$. Let $\sJ\subset \sO_{\Pi_\sU}$ be the image of $\bigwedge^2\sL^t$, giving us the Eagon-Northcott complex
\begin{equation}\label{eqn:EagonNorthcott1}
0\to \sO_{\Pi_\sU}^2(-3)\xrightarrow{\sL} \sO_{\Pi_\sU}(-2)^3\xrightarrow{\bigwedge^2\sL^t}\sO_{\Pi_\sU}\to \sO_{\Pi_\sU}/\sJ\to0
\end{equation}
on $\Pi_\sU$. We let $\sU^\cm:=(\bar\psi\circ Q)^{-1}(X\setminus F)$ and $\sU^\ncm:=(\bar\psi\circ Q)^{-1}(F)$.

\begin{lemma} The complex \eqref{eqn:EagonNorthcott1} is exact.
\end{lemma}

\begin{proof}  For $x:=(L_{ij})\in \sU^\cm$, the ideal $\sJ_x$ defines a CM-curve, hence has height 2 in $\sO_{x\times_{\Gr(4, n+1)}\Pi_n}\cong \sO_{\P^3}$. For $x\in \sU^\ncm$, the ideal $\sJ_x$ defines a (non-reduced) plane in  $x\times_{\Gr(4, n+1)}\Pi_n\cong \P^3$, hence has height 1. But since $F$ is a proper closed subset of the irreducible $k$-scheme $X$, this implies that the subscheme of $\Pi_\sU$ defined by $\sJ$ has codimension 2. By the theorem of Eagon-Northcott (Theorem~\ref{thm:EN}), this implies that the complex \eqref{eqn:EagonNorthcott1} is exact.
\end{proof}

\begin{lemma} Let $\Pi_X:=X\times_{\Gr(4, n+1)}\Pi$ with projections $q_{X1}:\Pi_X\to X$, $q_{X2}:\Pi_X\to \Pi$.\\[5pt]
1. The complex \eqref{eqn:EagonNorthcott1} descends to a complex
\[
0\to q_{X1}^*\sF_X(-3)\xrightarrow{\alpha_X} q_{X1}^*\sE_X(-2)^3\xrightarrow{\bigwedge^2\alpha_X^t}\sO_{\Pi_X}\to \sO_{\Pi_X}/\sJ_X\to0
\]
on $\Pi_X$. Here $\sJ_X\subset \sO_{\Pi_X}$ is defined as the image of $\bigwedge^2\alpha_X^t$. Moreover, this complex is exact.\\[2pt]
2. The map
\[
\beta_X:\sE_X\to q_{1X*}(\sO_{\Pi_X}(2))=r_X^*\sV
\]
induced by $\bigwedge^2\alpha_X^t$ defines an isomorphism of $\sE_X$ with $r_X^*(E_{3, \sV})$. 
\end{lemma}

\begin{proof}
 We use the standard basis $e_1, e_2, e_3$ for $ \sO_\sU^3$, and the standard basis  
$f_1, f_2$ for $\sO_\sU^2$. We use the basis $e_2\wedge e_3, -e_1\wedge e_3, e_1\wedge e_2$ for $\bigwedge^2\sO_\sU^3$ and the basis $f_1\wedge f_2$ for  $\bigwedge^2\sO_\sU^2$. With respect to these bases, 
\[
\bigwedge^2\sL:\bigwedge^2(\sO_{\Pi_\sU}^2(-3))\to \bigwedge^2(\sO_{\Pi_\sU}(-2))
\]
 has matrix
\[
\begin{pmatrix}e_{23}(L_{ij})\\-e_{13}(L_{ij})\\e_{12}(L_{ij})\end{pmatrix}
\]
at $(L_{ij})\in \sU$.  Taking the transpose and twisting gives us the map
\[
\bigwedge^2\sL^t:   (\bigwedge^2\sO_{\Pi_\sU}^3)^\vee(-2)\to   \det^{-1}\sO_{\Pi_\sU}^2
\]
with matrix $(e_{23}(L_{ij}),-e_{13}(L_{ij}),e_{12}(L_{ij}))$ at a point $(L_{ij})$,  in our choice of bases. 

We have the canonical isomorphism
\[
\sO_{\Pi_\sU}^3\cong (\bigwedge^2\sO_{\Pi_\sU}^3)^\vee\otimes\det \sO_{\Pi_\sU}^3
\]
so $\bigwedge^2\sL^t$ induces the map
\[
\tilde\bigwedge^2\sL^t:   \sO_{\Pi_\sU}^3(-2)\to    \det^{-1}\sO_{\Pi_\sU}^2\otimes \otimes\det \sO_{\Pi_\sU}^3.
\]

The group $\Gamma$ acts on $\det^{-1}\sO_{\Pi_\sU}^2\otimes \otimes\det \sO_{\Pi_\sU}^3$ by
\[
(g,h)\mapsto (\det h)^{-1}(\det g)^2(\det h)^{-2}\cdot \det g\cdot (\det g)^{-3}\cdot (\det h)^3=1
\]
and as we have constructed the map $\tilde\bigwedge^2\sL^t$ by natural operations from the map $\sL$,  $\tilde\bigwedge^2\sL^t$ descends to a map of locally free sheaves on $X$
\[
\bigwedge^2\alpha^t:q_{X1}^*\sE_X(-2)\to  \sO_{\Pi_X}
\]
Our ideal $\sJ_X$ is by definition the image of $\bigwedge^2\alpha_X^t$. 

Thus, descent by the action of $\Gamma$ on the complex \eqref{eqn:EagonNorthcott1} does give a complex of the asserted form on $\Pi_X$. Since the complex \eqref{eqn:EagonNorthcott1} is exact, so is the complex formed by descent. This proves (1).

For (2), we have the canonical isomorphism
\[
q_{X1*}(\sO_{\Pi_X}(2))=q_{X1*}(q_{2X}^*\sO_{\Pi}(2))=r_X^*p_{1*}*\sO_{\Pi}(2)=r_X^*\sV,
\]
justifying the identity $q_{1X*}(\sO_{\Pi_X}(2))=r_X^*\sV$ in the statement of (2). By construction if $(L_{ij})$ has image $x\in X$, then the image of $\sE_X\otimes k(x)$ in $H^0(\sO_{\Pi_{r(x)}}(2))=\sV\otimes k(x)$ is the subspace generated by the  $2\times 2$ minors of $(L_{ij})$. By definition of the  map $q$,  this is exactly the subspace  $E_{2, \sV}\otimes k(x)\subset \sV\otimes k(x)$, which proves (2).
\end{proof}

\begin{lemma}\label{lem:Blowup} The map $q:H_n\to X$ induces an isomorphism of $H_n$ with the blow-up $\mu_F:\Bl_FX\to X$.
\end{lemma}

\begin{proof}  The isomorphism $q^\cm:H^\cm_n\to X\setminus F$ induces a birational map $\phi:H_n\dashrightarrow \Bl_FX$ over $X$; to show that $\phi$ is an isomorphism, it suffices to show that after pullback by $\sU\to X$, the induced map $\phi_\sU:H_n\times_X\sU\dashrightarrow \Bl_{F_\sU}\sU$ is an isomorphism. We reduce to the case $n=3$ by fibering over $\Gr(4,n+1)$. Let $D\subset \Bl_FX$ denote the exceptional divisor.

In addition to the action of $\Gamma$, we have the action of $\GL_4$ on $\sU$ via its action on the linear forms on $\P^3$, which thus commutes with the $\Gamma$-action. The matrices in $F_\sU$ 
correspond closed subschemes with ideals of the form $(L_0^2, L_0L_1, L_0L_2)$, with $L_0, L_1, L_2$ linearly independent linear forms; from this, one sees that $F_\sU$ is a single $\GL_4\times\Gamma$ orbit. We take for base-point the matrix
\[
\alpha_0:=\begin{pmatrix} 0&-X_0\\X_0&0\\-X_1&X_2\end{pmatrix}
\]
Acting by the Lie algebra of $\GL_4\times \Gamma$, we find a basis for the normal bundle of $F_\sU$ in $\sU$ at $\alpha_0$ given by matrices of the form
\[
v_{A,B,C}:=\begin{pmatrix} C&B\\0&A\\0&0\end{pmatrix}
\]
with $C, B$ linear forms in $X_2, X_3$, $A$ a linear form in $X_1, X_2, X_3$. A   path $\gamma$ in $\Bl_{F_\sU}\sU$ lying over $\alpha_0$ and intersecting the exceptional divisor transversely thus corresponds to a path in $\sU$ of the form 
\[
\alpha_t:=\begin{pmatrix} tC&tB-X_0\\X_0&tA\\-X_1&X_2\end{pmatrix} +O(t^2).
\]
giving for $t\neq0$ the flat family of ideals $\sI_t$ with generators 
\[
e_{12}(t)=t^2AC-tBX_0+X_0^2,\ e_{13}(t)=tCX_2+tBX_1-X_0X_1,\ e_{23}(t)=X_0X_2+X_1tA
\]
We have the element
\[
tq:=X_2e_{13}(t)+X_1e_{23}(t)=t(AX_1^2+BX_1X_2+CX_2^2)
\]
so for $t\neq 0$, the element $q$ is in $\sI_t$. 

One can check by a direct computation that, if $f$ is a local section of $\sO_{\P^3_S}$ such that $tf$ is a local section of $\sI$, then $f$ is a local section of $\sI$, that is, $\sO_{\P^3_S}/\sI$ is flat over $k[t]_{(t)}$.  Clearly $\sI\otimes k(t)\subset \sO_{\P^3_{k(t)}}$ defines a $CM$-curve, and $\sI\otimes k\subset \sO_{\P^3}$  gives the subscheme $C_0$ defined by  the ideal $\sI_0$ with generators $(X_0^2, X_0X_1, X_0X_2, q)$.

We can compute the corresponding tangent vector at $[C_0]\in H_3$ as the element of $\Hom(\sI_{C_0}, \sO_{C_0})$ as induced by $-1$ times the map $\sI\to  \sO_{\P^3_S}/\sI_0\otimes_kk[t]/t^2=t\sO_{C_0}\oplus \sO_{C_0}$. This maps $\sI$ to the summand $t\sO_{C_0}$ and sends $t\sI$ to zero. Explicitly,  the map $\sI_{C_0}\to \sO_{C_0}$ is given by
\[
x_0^2\mapsto Bx_0,\ x_0x_1\mapsto Cx_2+Bx_1, x_0x_2\mapsto Ax_1,\ q\mapsto 0.
\]
From our conditions on $A, B, C$, it is easy to see that this map is zero as a map to the normal bundle if and only if $A=B=C=0$. 
This shows that $q_\sU^{-1}(F_\sU)\subset H_3\times_XF_\sU$ is a reduced Cartier divisor, and hence $q^{-1}(F)\subset H_3$ is also a reduced Cartier divisor.

As set-theoretically $q^{-1}(F)=H_3^\ncm$, the universal property of the blow-up implies that the morphism $q:H_3\to X$ lifts to a well-defined morphism $\tilde{q}:H_3\to \Bl_FX$, with $\tilde{q}^*(D)=H_3^\ncm$. The points $x$ of $H_3$ lying over $\alpha_0$ are uniquely determined by a choice of degree three generator for $\sI_{C_x}$, and this degree three generator is of the form $AX_1^2+BX_1X_2+CX_2^2$, with $A, B, C$ linear forms in $X_1, X_2, X_3$. We may take $B, C$ to be linear forms in $X_2, X_3$ and $A$ to be a linear form in $X_1, X_2, X_3$ without changing the resulting ideal, and then the choice of $A, B$ and $C$ is unique up to a (single) scalar. Thus,   the morphism $\tilde{q}_\sU$ defines a bijection from the projective space on   $F_\U\otimes k(\alpha_0)$  with the fiber of $q_\sU$ over $\alpha_0$, and hence $\tilde{q}$ gives a bijection of $H_3^\ncm$ with $D$. Since the map $\tilde{q}$ is an isomorphism over $X\setminus F$, $\tilde{q}$ is a birational bijection of smooth proper $k$-schemes, hence an isomorphism by Zariski's main theorem.
\end{proof}

Define locally free sheaves $\sF_n$, $\sE_n$ on $H_n$ by $\sF_n:=q^*\sF_X$, $\sE_n:=q^*\sE_X$.
Note that $\Pi_n=H_n\times_X\Pi_X$, so pulling back the complex
\[
0\to q_{X1}^*\sF_X(-3)\xrightarrow{\alpha_X} q_{X1}^*\sE_X(-2)\xrightarrow{\bigwedge^2\alpha_X^t}\sO_{\Pi_X} 
\]
to $\Pi_n$ via $q\times\id$ gives a complex
\[
0\to \pi_1^*\sF_n(-3)\xrightarrow{\alpha_n} \pi_1^*\sE_n(-2)\xrightarrow{\bigwedge^2\alpha_n^t}\sO_{\Pi_n} 
\]
on $\Pi_n$. Let $\sJ_n\subset \sO_{\Pi_n}$ be the image of $\bigwedge^2\alpha_n^t$, giving the closed subscheme $\sC_n'$  of $\Pi_n$ and the complex
\begin{equation}\label{eqn:EagonNorthcott3}
0\to \pi_1^*\sF_n(-3)\xrightarrow{\alpha_n} \pi_1^*\sE_n(-2)\xrightarrow{\bigwedge^2\alpha_n^t}\sO_{\Pi_n} \to \sO_{\sC_n'}\to0
\end{equation}

Let $\Pi^\ncm_n\to H_n^\ncm$ be the pullback of $\Pi_n$ to $H_n^\ncm$ let $\bar{\sJ}_n$ be the image of $\sJ_n$ in  $\sO_{\Pi^\ncm_n}$ and let $\sP_n\subset \Pi^\ncm_n$ be the closed subscheme defined by $\sqrt{\bar{\sJ}_n}$. Then $\sP_n\to H_n^\ncm$ is a $\P^2$-bundle, linearly embedded in the $\P^3$-bundle $\Pi^\ncm_n\to H_n^\ncm$. Let $\iota:\sP_n\to \Pi_n$ denote the inclusion.

\begin{lemma}\label{lem:ENCohomology} 1. The complex \eqref{eqn:EagonNorthcott3} is exact.\\[2pt]
2. We have $\sJ_n\subset \sI_{\sC_n}$ and the quotient $\sI_{\sC_n}/\sJ_n$ is of the form $\iota_*(\sL)$ for a certain invertible sheaf $\sL$ on $\sP_n$.   \\[2pt]
3. For each $x\in H_n^\ncm$, the invertible sheaf $\sL\otimes k(x)$ on $\sP_{n, x}\cong \P^2_{k(x)}$ is isomorphic to $\sO_{\sP_{n, x}}(-3)$.
\end{lemma}

\begin{proof} We have already proven that,  after restriction to $H^\cm$,  (1) and (2) hold and $\sJ_n= \sI_{\sC_n}$. The statement is also local on $H_n$ in the flat topology, so we need only prove the lemma after  pulling back via  $\Bl_{F_\sU}\sU\to H_n$ and restricting to a neighborhood of the exceptional divisor $D\subset \Bl_{F_\sU}\sU$.

As $H_n\to \Gr(4, n+1)$ is $\GL_{n+1}$ equivariant, we can restrict to a fiber over the standard $\P^3\subset \P^{n+1}$ defined by $x_4=\ldots=x_n=0$, so we reduce to the case $n=3$. 

For $n=3$,  $F_\sU$ is a homogeneous space for the action of $ \Gamma$ on $\sU$, the cohomology of the pullbacks of \eqref{eqn:EagonNorthcott3} and the quotient $\sI_{\sC_n}/\sJ_n$ are both flat over $F_\sU$, so we can restrict our attention to a fiber $D_x$ of $D\to F_\sU$ over a single point $x\in F_\sU(k)$. Similarly,  it suffices to prove that we can find a flat morphism $\pi:S\to H_3$ with $S$ a smooth $k$-scheme such that
\\[5pt]
1.   The the closed subscheme $V:=S\times_{H_3}D$ of $S$ is a smooth Cartier divisor on $S$, the projection $V\to D$ has  image $D_x$ and the induced map $V\to D_x$ is smooth.\\[2pt]
2.   The cartesian square
\[
\xymatrix{
V\ar[r]\ar[d]&S\ar[d]^\pi\\
D\ar@{^(->}[r]&H_3
}
\]
is transverse in the category of smooth $k$-schemes,  that is,  the induced map on the normal bundles $N_V\to i_{D_x}^*N_D$ is an isomorphism\\[2pt]
3.  The statement of the lemma holds after pulling back by $\pi:S\to H_3$.
\\[5pt]

We take $x$ to be the matrix
\[
\begin{pmatrix}0&-x_0\\x_0&0\\x_1&x_2\end{pmatrix}.
\]
The set of points $y$ lying over $x$ are parametrized by the non-zero cubic forms
\[
q:=Ax_1^2+Bx_1x_2+Cx_2^2
\]
with $A, B, C$ as in the proof of Lemma~\ref{lem:Blowup}. We take as parameter space $V$ the complement of 0 in the affine space of dimension 7 that parametrizes the cubic forms $q$, that is
\[
A=a_1x_1+a_2x_2+a_3x_3,\ B=b_2x_2+b_3x_3,\ C=c_2x_2+c_3x_3
\]
and $V=\Spec k[a_1, a_2, a_3, b_2, b_3,  c_2, c_3]\setminus\{0\}$.  Our fiber in $D$ is thus $\P(V)$, so we have the canonical smooth surjective  map $\pi_0:V\to D$

We take as family of  curves  over $V\times \Spec k[t]_{(t)}$
\[
\alpha_t(q):=
\begin{pmatrix}tC&tB-x_0\\x_0&tA\\x_1&x_2\end{pmatrix}.
\]

Let $S=V\times \Spec k[t]_{(t)}$. The corresponding flat family of subschemes $C\subset \P^3_{S}$ is given by the sheaf of ideals
$\sI$ generated by $(e_{12}(\alpha_t(q)), e_{13}(\alpha_t(q)), e_{23}(\alpha_t(q)), q)$; the corresponding ideal sheaf $\sJ$ is   generated by  $(e_{12}(\alpha_t(q)), e_{13}(\alpha_t(q)), e_{23}(\alpha_t(q)))$. Since $C':=\Spec \sO_{\P^3_S}/\sJ$ has codimension 2 in $\P^3_S$, the Eagon-Northcott complex
\[
 0\to\sO_{\P^3_S}(-3)^2\xrightarrow{\alpha_t(q)}\sO_{\P^3_S}(-2)^3\xrightarrow{\bigwedge^2\alpha_t(q)^t}
 \sO_{\P^3_S}\to \sO_{C'}\to0
 \]
 is exact. Putting  $\sO_{\P^3_S}$ in degree zero,  the cohomology of 
 \[
0\to  \sO_{\P^3_S}(-3)^2\xrightarrow{\alpha_t(q)}\sO_{\P^3_S}(-2)^3\xrightarrow{\bigwedge^2\alpha_t(q)^t}
 \sO_{\P^3_S}\to \sO_C\to0
 \]
 is equal to $\sI/\sJ$, supported in degree zero.
 
 As a $\sO_{\P^3_S}$ module, $\sI/\sJ$ is generated by the image of $q$. Since 
 \[
 tq=x_2e_{13}(\alpha_t(q))+x_1e_{23}(\alpha_t(q))
 \]
we have $t\sI\subset \sJ$,  $t\cdot \sI/\sJ=0$ and
 \[
 \sI/\sJ=\sI/(\sJ+t\sI)\cong (\sI/t\sI)/(\sJ/t\sI).
 \]

 We write $0$ for $\Spec k[t]_{(t)}/(t)\subset \Spec k[t]_{(t)}$ and let $C_0\subset \P^3_k$ be fiber of $C$ over $V\times0\hookrightarrow S$.  Since $C$ is flat over $S$, the sequence
 \[
 0\to \sI/t\sI\to \sO_{\P^3_S}/(t)\to \sO_{C_0}\to 0
 \]
 is exact. Thus $\sJ/t\sI\subset \sI/t\sI\subset \sO_{\P^3_S}/(t)$  is given by the ideals  
 \[
 (x_0^2, x_0x_1, x_0x_2)\sO_{\P^3_S}/(t)\subset (x_0^2, x_0x_1, x_0x_2, q)\sO_{\P^3_S}/(t).
 \]
Thus as $\sO_{\P^3_S}$-module, $\sI/\sJ$ is the $\sO_{\P^3_S}/(t)=\sO_{\P^3_V}$-module corresponding to the graded module $(x_0^2, x_0x_1, x_0x_2, q)/ (x_0^2, x_0x_1, x_0x_2)$ over the sheaf of polynomial rings $R:=\sO_V[x_0, x_1, x_2, x_3]$. This is a cyclic module with generator the image $[q]$ of $q$, so we need only compute the annihilator of $[q]$.

Clearly $x_0q$ is in $\sJ$, so $x_0\in \Ann_R[q]$. On the other hand, if $P\cdot[q]=0$ for some $P\in R$, we have $\alpha, \beta, \gamma\in R$ with
\[
Pq=\alpha\cdot x_0^2+\beta\cdot x_0x_1+\gamma\cdot x_0x_2
\]
Since   $q$ does not involve $x_0$ we have $P=x_0P'$ for some $P'\in R$. Thus $\Ann_R[q]=x_0R$. 

The plane $\iota_x:\sP_{3,x}\hookrightarrow \P^3_S$ is the subscheme of $\P^3_S$ defined by $(t,x_0)$.  As $q$ has degree three, this implies
\[
\sI/\sJ\cong \iota_{x*}\sO_{\sP_{3,x}}(-3)
\]
which completes the proof.
\end{proof}

The complex \eqref{eqn:EagonNorthcott3} and inclusion $\sI_{\sJ_n}\subset \sI_{\sC_n}$ gives us the complex
\begin{equation}\label{eqn:EagonNorthcott4}
0\to p^*\sF_n(-3)\xrightarrow{\alpha_n} p^*\sE_n(-2)\xrightarrow{\bigwedge^2\alpha_n^t}\sO_{\Pi_n} \to \sO_{\sC_n}\to0
\end{equation}
with $\sO_{\Pi_n}$ in degree zero

For later use, we note the following consequence of  Lemma~\ref{lem:ENCohomology}.
\begin{proposition}\label{prop:MainResolution} Retaining the notation of Lemma~\ref{lem:ENCohomology}, the cohomology sheaves $\sH^i$ of   the complex \eqref{eqn:EagonNorthcott4} satisfy\\[5pt]
1. $\sH^i=0$ for $i\neq0$\\[2pt]
2. There is a canonical isomorphism $\sH^0\cong \iota_*(\sL)$.
\end{proposition}

\section{Determinants and canonical bases}\label{sec:Det}

Recall the morphisms $r_X:X\to \Gr(4, n+1)$ and $\Phi:H_n\to \Gr(4, n+1)$. 
\begin{lemma}\label{lem:TangentDet} 1. We have the following presentation of the relative tangent sheaf $\sT_{r_X}$:
\begin{equation}\label{eqn:TangentPresentation}
0\to \sO_X\to \sEnd(\sE_X)\oplus \sEnd(\sF_X)\to \sHom(\sF_X, \sE_X)\otimes q_{X1*}\sO_{\Pi_X}(1)\to \sT_{r_X}\to 0
\end{equation}
2. Let $j:H_n^\cm\to H_n$ be the inclusion. We have the presentation of $j^*\sT_\Phi$
\[ 
0\to \sO_{H_n^\cm}\to j^*\sEnd(\sE)\oplus j^*\sEnd(\sF)\to j^*\sHom(\sF, \sE)\otimes \pi_{1*}\sO_{\Pi_n}(1)\to j^*\sT_\Phi\to 0
\]
3. We have a canonical isomorphisms
\[
\det \sT_{\Phi}\cong  (\det\sF^{\otimes -3}\otimes\det\sE^{\otimes 2})^{\otimes 4}\otimes 
\det(\pi_{1*}\sO_{\Pi_n}(1)))^{\otimes 6}\otimes \sO_{H_n}(-6\cdot H_n^\ncm) 
\]
and
\[
\det \sT_{H_n}\cong\det \sT_{\Phi}\otimes\Phi^*\sO_{\Gr(4,n+1)}(n+1)
\]
\end{lemma}

\begin{proof}  Let $Q:\sU\to X$ denote the quotient map, $r_\sU:\sU\to \Gr(4, n+1)$ the structure map.  We have the vector bundle $M\to \Gr(4, n+1)$ associated to the locally free sheaf $\sM=\sHom(\sO_\Gr^2, \sO_\Gr^3)\otimes p_{1*}\sO_\Pi(1)$, and $\sU$ is an open subscheme of $M$. Thus $\sT_{r_\sU}=r_\sU^*\sM$.

 A section $m=(m_{ij})$ of $r_\sU^*\sM=\sHom(\sO_\sU^2, \sO_\sU^3)\otimes r_\sU^* p_{1*}\sO_\Pi(1)$  is a matrix of linear forms on the fibers of the  $\P^3$-bundle $\Pi_\sU\to \Gr(4, n+1)$,  with a matrix $m$ giving the tangent vector at $L\in \sU\subset \sHom(\sO_\sU^2, \sO_\sU^3)\otimes q_{\sU1*}\sO_{\Pi_\sU}$ defined by $L+\epsilon m$. 

Since $X$ is the quotient of $\sU$ by $\Gamma$, the pullback $Q^*\sT_{r_X}$  is the quotient of $\sM$ by the action of the Lie algebra of $\Gamma$, and $\sT_{r_X}$ itself is thus given by descent via the natural action of $\Gamma$ on this quotient. Explicitly, the Lie algebra of $\GL_3\times\GL_2$ acts on $\sM$ at a point $L$ of $\sU$ by 
\[
(a,b)\cdot m=m+a\cdot L-L\cdot b
\]
The map $\G_m\to \GL_3\times\GL_2$ induces the map of Lie algebras $\lambda\mapsto (\lambda\cdot \id, \lambda\cdot \id)$ giving the presentation of $Q^*\sT_{r_X}$
\[
0\to \sO_\sU\to \sEnd(\sO_\sU^3)\oplus \sEnd(\sO_\sU^2)\to 
\sHom(\sO_\sU^2, \sO_\sU^3)\otimes r_\sU^* p_{1*}\sO_\Pi(1)\to Q^*\sT_{r_X}\to0.
\]
Applying $\Gamma$-descent  yields the asserted presentation of $\sT_{r_X}$.

(2) follows from (1) and the fact that $q$  restricts to an isomorphism $q^\cm:H_n^\cm\to X\setminus F$.

For (3), we have
\[
\rnk \sF_X=2,\ \rnk \sE_X=3,\ \rnk q_{X1*}\sO_{\Pi_X}(1)=\dim_kH^0(\P^3, \sO_{\P^3}(1))=4
\]
so (1) gives us the canonical isomorphism
\[
\det\sT_{r_X}\cong  (\det\sF_X^{\otimes -3}\otimes\det\sE_X^{\otimes 2})^{\otimes 4}\otimes 
\det(q_{X1*}\sO_{\Pi_X}(1))^{\otimes 6}
\]
$H_n$ is isomorphic to the blow-up of $X$ along the codimension 7 closed subscheme $F$ with exceptional divisor $H_n^\ncm$. Since $F$ is smooth over $\Gr(4, n+1)$, the formula for the relative canonical class of a blow-up yields our formula for $\det \sT_{\Phi}$. For the final formula, we have the exact sequence
\[
0\to  \sT_{\Phi}\to \sT_{H_n}\to \Phi^*\sT_{\Gr(4, n+1)}\to 0
\]
so the canonical isomorphism 
$\det\sT_{\Gr(4, n+1)}\cong \sO_{\Gr(4, n+1)}(n+1)$
 yields our formula for $\det\sT_{H_n}$.
\end{proof}

Let $\sE_{m,n}:=p_*\sO_{\sC_n}(m)$.   

\begin{lemma}\label{lem:EnmPresentation} For $m\ge0$, $\sE_{m,n}$ is a locally free sheaf on $H_n$ of rank $3m+1$ and $R^ip_*\sO_{\sC_n}(m)=0$ for $i>0$. 
\end{lemma}

\begin{proof}  We have $R^ip_*\sO_{\Pi}(j)=0$ for $j>-3$; applying this to the resolution  \eqref{eqn:TangentPresentation} gives the exact sequence in (2). This also shows that $j^*\sE_{m,n}$ is a locally free sheaf on $H^\cm_n$ of rank $3m+1$ and that $R^ip_*\sO_{\sC_n}(m)=0$ for $i>0$. For $x\in H^\ncm$, the plane curve $\bar{C}_x$ is a plane cubic, showing that $H^0(\bar{C}_x, \sO_{\bar{C}_x}(m))$ has dimension $3m$ and $H^i(\bar{C}_x, \sO_{\bar{C}_x}(m))=0$ for $i>0$. This implies $H^0(C_x, \sO_{\bar{C}_x}(m))$ has dimension $3m+1$ and $H^i(C_x, \sO_{\bar{C}_x}(m))=0$ for $i>0$. As $\sO_{\sC_n}(m)$ is flat over $H_n$, this completes the proof.
\end{proof}

We make a brief detour to define a canonical relative orientation on $\Sym^mW$, for $W$ a given $k$-vector space. Suppose $W$ has dimension $r$ over $k$ and let $w_1,\ldots, w_r$ be a basis. We order the basis $\Sym^mW$  of monomials of degree $m$ in the $w_j$ lexicographically, using the order $w_1>w_2>\ldots>w_r$ on the $w_j$. To fix this notion, $w_1^m$ will come first in this order.

Let $\sym^m_r$ denote the dimension of $\Sym^mW$ and let $e_1(w_*),\ldots, e_{\sym^m_r}(w_*)$ denote the lexicographically ordered basis of $\Sym^mW$.

\begin{lemma}\label{lem:DetSymIso}  Let
\[
N_{m,r}:= \binom{m+r-1}{r}
\]
Then as a rank one $\Aut_k(W)$ representation there is a unique isomorphism
 \[
 \phi_{W,m}:\det(W)^{\otimes N_{m,r}}\xrightarrow{\sim}\det\Sym^mW
 \]
 satisfying
 \[
 \phi_{W,m}((w_1\wedge\ldots w_r)^{\otimes N_{m,r}})=e_1(w_*)\wedge\ldots\wedge e_{\sym^m_r}(w_*)
 \]
 for each basis $w_1,\ldots, w_r$ of $W$.
 \end{lemma}
 
 \begin{proof} As $\det\Sym^mW$ is a rank one representation of $\Aut_k(W)$, there is a unique integer $N$ such that $\det\Sym^mW$ is isomorphic to the representation $\det(W)^{\otimes N}$, and one can determine $N$ by examining the action of $t\cdot \id_W$. On  $\det(W)^{\otimes N}$, 
 $t\cdot \id_W$ acts by $t^{rN}$, while on $\det\Sym^mW$, $t\cdot \id_W$ acts by $t^{m\cdot\sym^m_r}$. Thus
 \[
 N=\frac{m}{r}\cdot \binom{m+r-1}{m-1}=
 \frac{m}{r}\cdot \binom{m+r-1}{r}=  \binom{m+r-1}{r-1}=N_{m,r}
 \]
 
 To see that $\phi:=\phi_{W,m}$ is a well-defined $\Aut_kW$-equivariant isomorphism, we temporarily write $\phi$ as $\phi_{w_*}$, to indicate the possible dependence on the choice of ordered basis $w_*:=(w_1,\ldots, w_r)$. Then we need only show that
\[
\phi_{\alpha(w_*)}((\alpha(w_1)\wedge\ldots\wedge\alpha(w_r))^{\otimes N_{m,r}}))=\det\Sym^m(\alpha)(\phi_{w_*}((w_1\wedge\ldots\wedge w_r)^{\otimes N_{m,r}}))
\]
 for $\alpha\in \Aut_k(W)$ of the form
  \\[2pt]
 1. An elementary matrix $e_{ij}(\lambda)$ acting on $w_1,\ldots, w_r$ by $w_\ell':=e_{ij}(\lambda)(w_\ell)$ with $w_\ell'=w_\ell$ for $\ell\neq i$ and $w_i'=w_i+\lambda\cdot w_j$. \\
2. A diagonal automorphism $(t_1,\ldots, t_r)$ with $w_j':=(t_1,\ldots, t_r)(w_j):=t_jw_j$, \\
 3. A permutation $\sigma\in \Sigma_r$ acting by $w_j':=\sigma(w_j):=w_{\sigma(j)}$\\[2pt]
 \\[2pt] 
 For (1),  $\Sym^m(e_{ij}(\lambda))$ sends $e_\ell(w_*)$ to $e_\ell(w_*')$. This map is represented by a uni-potent matrix with respect to the basis $e_*(w_*)$, and thus
 \begin{align*}
 \phi_{w_*'}((w'_1\wedge\ldots\wedge w'_r)^{\otimes N_{mr}}))&=
 e_1(w_*')\wedge \ldots\wedge e_{\sym^m_r}(w_*')\\
 &=e_1(w_*)\wedge \ldots\wedge e_{\sym^m_r}(w_*)\\
 &=\Sym^m(e_{ij}(\lambda))(\phi_{w_*}((w_1\wedge\ldots\wedge w_r)^{\otimes N_{mr}}))
 \end{align*}
 For (2), this follows by a similar argument,  since the action of $(t_1,\ldots, t_r)$ sends the ordered basis $w_1,\ldots, w_r$ to the ordered basis $t_1w_1,\ldots, t_rw_r$. 
 
 For (3), the permutation $\sigma$ sends the $j$th ordered monomial $e_j(w_1,\ldots, w_r)$, with respect to the lexicographical order for $w_1>\ldots> w_r$,  to the $j$th ordered monomial $e_j(w_{\sigma(1)},\ldots, w_{\sigma(r)})$ with respect to the lexicographical order for $w_{\sigma(1)}>\ldots> w_{\sigma(r)}$, in other words,
 \[
e_j(w'_*)=\Sym^m(\sigma)(e_j(w_*))
\]
and thus
\[
\phi_{w_*'}((w'_1\wedge\ldots\wedge w'_r)^{\otimes N_{mr}}))=
\Sym^m(\sigma)(\phi_{W,m}((w_1\wedge\ldots\wedge w_r)^{\otimes N_{mr}}))
\]
 \end{proof}
  
\begin{proposition}\label{prop:DetSymIso} Let $\sW$ be a locally free coherent sheaf  of rank $r$ on a $k$-scheme $X$.  Then there is a unique isomorphism of invertible sheaves  $\phi_{\sW,m}:(\det \sW)^{\otimes N_{m,r}}\xrightarrow{\sim} \det\Sym^m\sW$ such that, if $w_1,\ldots, w_r$ is a framing for $\sW$ over some open subscheme $U$ of $X$, then
\[
\phi_{\sW,m}((w_1\wedge\ldots\wedge w_r)^{\otimes N_{m,r}})=e_1(w_*)\wedge\ldots e_{\sym^m_r}(w_*)
\]
where $e_j(w_*)$ is the $j$th monomial section of $\Sym^m\sW$ over $U$, with respect to the lexicographical order on the monomials in the $w_j$ induced by the order $w_1>\ldots>w_r$.\end{proposition}
\begin{proof} This is an immediate consequence of Lemma~\ref{lem:DetSymIso}.
\end{proof}

We will be using a modification of the lexicographical orientation, which we now describe. As we will only be using this for a rank four bundle, we restrict to $r=4$. 
\begin{definition}\label{def:CanonicalOrientation} 1. Given an ordered basis $x_1,\ldots, x_4$ for $W$ we call a monomial 
$x_1^{m_1}x_2^{m_2}x_3^{m_3}x_4^{m_4}$ {\em positive} if $m_1>m_2$ or if $m_1=m_2$ and $m_3>m_4$, {\em negative} if 
$m_1<m_2$ or if $m_1=m_2$ and $m_3<m_4$, and  {\em neutral} if $m_1=m_2$ and $m_3=m_4$.\\[2pt]
2. If $M:=x_1^{m_1}x_2^{m_2}x_3^{m_3}x_4^{m_4}$ is a positive monomial, define $M^*:=
x_1^{m_2}x_2^{m_1}x_3^{m_4}x_4^{m_3}$.\\[2pt]
3. Note that $M\mapsto M^*$ defines a bijection of the positive monomials with the negative monomials and is the identity on the neutral monomials. The {\em canonical order} on monomials in $\Sym^mW$ is defined as follows: If $M_1$ and $M_2$ are positive or neutral monomials, set $M_1>M_2$ if $M_1$ is larger than $M_2$ in the lexicographical order.   If $M_1>M_2$ are positive monomials and $M_3$ is neutral, we set $M_1>M_1^*>M_2>M_2^*>M_3$. \\[2pt]
4. Let $\sigma_m(x_*)$ be the permutation of the set of the monomials of $\Sym^mW$ that transforms the lexicographical order to the canonical order, and let $\epsilon_m$ be the sign of $\sigma_m(x_*)$. \\[2pt]
5. Let $\sW$ be a rank four invertible sheaf on some $k$-scheme $X$. The {\em canonical orientation} on $\Sym^m\sW$ is the isomorphism
\[
\epsilon_m\phi_{\sW,m}:\det\sW^{\otimes N_{m,4}}\xrightarrow{\sim} \det\Sym^m\sW
\]
\end{definition}
\begin{remark} 1. Definition~\ref{def:CanonicalOrientation}(5) is justifed by the fact that $\epsilon_m$ is independent of the choice of ordered basis $x_1,\ldots, x_4$ for $W$. One can verify this along the lines of the proof of Lemma~\ref{lem:DetSymIso}.\\[2pt]
2. Fix $m$ and choose positive integers $a_1, a_2$ with $a_1> m\cdot a_2$, give, $x_1$ weight $a_1$, $x_2$ weight $-a_1$, $x_3$ weight $a_1$ and $x_4$ weight $-a_4$. Then a monomial
$M:=x_1^{m_1}x_2^{m_2}x_3^{m_3}x_4^{m_4}$ is positive (resp. negative, resp. neutral) if and only if  $\wt(M)>0$, resp. $\wt(M)<0$, resp. $\wt(M)=0$, that is, if and only if $a_1(m_1 -m_2)+a_2(m_3-m_4)$ is positive, resp. negative, resp. zero.  Moreover, if $a_1>(2m-1)a_2$, then ordering the positive monomials by weight is the same as the lexicographical order and  the canonical order on the non-neutral monomials is the order by the absolute value of the weight, with the positive monomial of weight $w$ immediately preceeding the negative monomial of weight $-w$. \end{remark}

\begin{ex}  We illustrate with the example of degree 2. The lexicographical order is
\[
x_1^2>x_1x_2>x_1x_3>x_1x_4>x_2^2>x_2x_3>x_2x_4>x_3^2>x_3x_4>x_4^2
\]
The positive monomials are
\[
x_1^2, x_1x_3, x_1x_4, x_3^2
\]
The negative monomials are
\[
x_2^2, x_2x_3 x_2x_4, x_4^2
\]
the neutral monomials are
\[
x_1x_2, x_3x_4
\]
The canonical order is
\[
x_1^2>x_2^2>x_1x_3>x_2x_4> x_1x_4>x_2x_3> x_3^2>x_4^2>x_1x_2>x_3x_4
\]

If we give $x_1$ weight $3$, $x_2$ weight $-3$, $x_3$ weight 1 and $x_4$ weight -1, the weights of the monomials in the canonical order are
\[
6, -6, 4, -4, 2,-2, 2, -2, 0,0
\]
If we give $x_1$ weight $4$, $x_2$ weight $-4$, $x_3$ weight 1 and $x_4$ weight -1, the weights of the monomials in the canonical order are
\[
8, -8, 5, -5, 3,-3, 2, -2, 0,0
\]
\end{ex}

Let $V$ be a fixed $k$-vector space of dimension $n+1$ and let $E_{4,n+1}\to \Gr(4, n+1)$ be the tautological rank four subsheaf of $\sO_{\Gr(4, n+1)}^{n+1}$, giving the exact sequence
\[
0\to E_{4,n+1}\to \sO_{\Gr(4, n+1)}\otimes V\to \sQ_{n-3, n+1}\to 0
\]
We have the similarly defined sequence on $\P^n$, 
\[
0\to \sO_{\P^n}(-1)\to \sO_{\P^n}\otimes V\to \sQ_{n,n+1}\to 0
\]

\begin{lemma}\label{lem:OPiDet} Let $N_m=N_{m,4}$.  Then for all $m\ge0$ we have canonical isomorphisms 
\[
\phi_{m,n}:\sO_{\Gr(4,n+1)}(N_m) \xrightarrow{\sim}\det p_{1*}(\sO_\Pi(m))
\]
\[
\psi_n:\sO_{\Gr(4,n+1)}(1)\xrightarrow{\sim} \det E_{4,n+1}^\vee.
\]
with 
\[
\phi_{m,n}\circ\psi_n^{\otimes N_m}=\epsilon_m\cdot \phi_{m, E_{4,n+1}^\vee}
\]
\end{lemma}

\begin{proof} The Pl\"ucker embedding is defined by the global sections of $\det E^\vee_{4,n+1}$, which gives the isomorphism $\psi_n$. We then define $\phi_{m,n}:=
\epsilon_m\cdot \phi_{m, E_{4,n+1}^\vee}\circ \psi_n^{\otimes N_m}$.
\end{proof}

\begin{remark} Given a framing $x_1,\ldots, x_4$ of $E_{4, n+1}^\vee$ over some open subset $U$ of $\Gr(4, n+1)$, this gives a corresponding isomorphism $\sO_{\Gr(4, n+1)}(1)\cong \sO_{\Gr(4, n+1)}$ over $U$ and if $s_1,\ldots, s_{\sym^m_4}$ is the canonical framing of $ p_{1*}(\sO_\Pi(m))$ with respect to $x_1,\ldots, x_4$, then 
\[
\phi_{m,n}(1)=s_1\wedge\ldots\wedge s_{\sym^m_4}
\]
\end{remark}

We have the $\P^2$-bundle $\rho_n:\sP_n\to H^\ncm$ and the invertible sheaf $\sL$ on $\sP_n$ given by Lemma~\ref{lem:ENCohomology}.  

We twist  the complex \eqref{eqn:EagonNorthcott4} by $\sO_{\Pi_n}(m)$. This gives us  the complex on $\sC_n$
\begin{equation}\label{eqn:ResOnCn}
0\to p^*\sF_n(m-3)\xrightarrow{\alpha_n} p^*\sE_n(m-2)\xrightarrow{\bigwedge^2\alpha_n^t}\sO_{\Pi_n}(m) \to \sO_{\sC_n}(m)\to0
\end{equation}
with $\sO_{\Pi_n}(m)$ put in degree zero.

\begin{lemma}\label{lem:MainRes} Take $m\ge1$. If we apply $p_*$ to the complex \eqref{eqn:ResOnCn}, we obtain the complex
\begin{equation}\label{eqn:ResOnHn}
0\to \sF_n\otimes p_*[\sO_{\Pi_n}(m-3)]\xrightarrow{\alpha_n} \sE_n\otimes  p_*[\sO_{\Pi_n}(m-2)]\xrightarrow{\bigwedge^2\alpha_n^t}p_*\sO_{\Pi_n}(m) \to \sE_{m,n}\to0.
\end{equation}
Moreover, this complex has cohomology sheaves $\sH^i=0$ for $i\neq0$ and 
$\sH^0\cong i_{H^\ncm_n*}(\rho_{n*}(\sL(m)))$.
\end{lemma}

\begin{proof} This follows from Proposition~\ref{prop:MainResolution}, together with the fact that $\sO_{\P^3}(m)$, $\sO_{\P^3}(m-3)$,  $\sO_{\P^3}(m-2)$ $\sO_{\sC_n}(m)$ and  $\sO_{\P^2}(m-3)$ have vanishing higher cohomology for $m\ge1$. 
\end{proof}

\begin{proposition}\label{prop:DetCMSectionBundle} Let 
\[
M_m:=-N_m+3N_{m-2}-2N_{m-3}=\frac{(3m-1)m}{2}.
\]
Let
\[
K_m:=\binom{m-1}{2}.
\]
Then for all $m\ge1$, the isomorphisms $\phi_{j,n}$ and the complex of Lemma~\ref{lem:MainRes} give rise to a canonical isomorphism
\[
\rho_{m,n}:\det\sE_{m,n}\xrightarrow{\sim} \Phi^*\sO_{\Gr(4,n+1)}(M_m)\otimes \det\sE^{\otimes-\binom{m+1}{3}}\otimes\det\sF^{\otimes\binom{m}{3}}\otimes\sO_{H_n}(K_mH_n^\ncm).
\]
\end{proposition}

\begin{proof}  Since
\[
N_m=N_{m,4}=\binom{m+3}{4}
\]
we have 
\[
M_m=\binom{m+3}{4}-3\binom{m+1}{4}+2\binom{m}{4}=\frac{(3m-1)m}{2}
\]
as claimed.

The complex of Lemma~\ref{lem:MainRes} gives us the canonical isomorphism
\begin{align*}
\det\sE_{m,n}&=\det p_*\sO_{\sC_n}(m)\\
&\cong\det (\sF_n \otimes p_*[\sO_{\Pi_n}(m-3)])\otimes  \det(\sE_n^3 \otimes  p_*[\sO_{\Pi_n}(m-2)])^{-1}\\&\hskip120pt\otimes \det p_*[\sO_{\Pi_n}(m)] \otimes \det i_{H^\ncm_{n*}}(\rho_{n*}(\sL(m))).
\end{align*}
 Recalling that the restriction of $\sL$ to a fiber of $\rho_n$ is canonically isomorphic to $\sO_{\P^2}(-3)$, we find that, for $x\in H^\ncm_{n*}$, $\rho_{n*}(\sL(m))\otimes k(x)$ has dimension 
\[
h^0(\P^2, \sO_{\P^2}(m-3))=\binom{m-1}{2}=K_m
\]
so 
\[
\det   i_{H^\ncm_{n*}}(\rho_{n*}(\sL(m)))=\sO_{H_n}(K_m\cdot H^\ncm).
\]
Our isomorphism $\phi_{j,n}:\sO_{\Gr(4, n+1)}(N_j)\xrightarrow{\sim} \det p_{1*}\sO_{\Pi}(j)$ of Lemma~\ref{lem:OPiDet}, together with the fact
\[
\dim_kH^0(\P^3, \sO_{\P^3}(j))=\binom{ j+3}{3}
\]
and the formula for  $\det(\sA\otimes\sB)$ gives us the canonical isomorphisms
\[
\det \sF_n \otimes p_*[\sO_{\Pi_n}(m-3)]=\det\sF^{\otimes\binom{m}{3}}\otimes \sO_{H_n}(2N_{m-3})
\]
\[
\det \sE_n(m-2)^3 \otimes  p_*[\sO_{\Pi_n}(m-2)]=\det\sE^{\otimes\binom{m+1}{3}}\otimes \sO_{H_n}(3N_{m-2})
\]
and 
\[
\det p_*\sO_{\Pi_n}(m)=\sO_{H_n}(N_{m})
\]
Fitting this together yields our formula for $\det\sE_{m,n}$. 
\end{proof}

\begin{definition} We say that a locally free  sheaf $\sP$  on  a smooth $k$-scheme $Y$ is {\em relatively oriented} if $\det\sP\otimes\det\sT_Y$ is a square in $\Pic(Y)$. A choice of an invertible sheaf $\sL$ on $Y$ and an isomorphism
$\rho:\det\sP\otimes\det\sT_Y\xrightarrow{\sim}\sL^{\otimes 2}$  is a {\em relative orientation} for $\sP$
\end{definition}

The group $\Gamma$ has a  free rank one character group, with generator  the character $(g, h)\mapsto (\det g)^2\cdot (\det h)^{-3}$ of $\GL_3\times \GL_2$. Let $\gamma\in \Pic(X)$ be the class of the invertible sheaf on $X=\Gamma\backslash\sU$ induced by this character of $\Gamma$.

\begin{lemma} \label{lem:Pic} Suppose $n\ge 4$. \\[5pt]
1. $\Pic(X)$ is free abelian group of rank 2 with generators $\gamma$ and $r_X^*\sO_{\Gr(4, n+1)}(1)$.\\[2pt]
2. $\Pic(H_n)=\Z^3$ with generators $\Phi^*\sO_{\Gr(4, n+1)}(1)$, $q^*\gamma$ and $\sO_{H_n}(H_n^{\ncm})$.\\[5pt]
If $n=3$, $\Pic(X)=\Z$ with generator $\gamma$ and $\Pic(H_n)=\Z^2$ with generators   $q^*\gamma$ and $\sO_{H_n}(H_n^{\ncm})$.
\end{lemma}

\begin{proof} We first consider the case $n\ge 4$. (2) follows from (1) and the formula for the Picard group of the blowup.

For (1), the scheme $\sU$ is an open subscheme of the vector bundle $M$ over $\Gr(4, n+1)$. Moreover, in a fiber over $x\in \Gr(4, n+1)$, it follows from \cite[Lemma 2]{EPS} that $M_x\setminus\sU_x$ has codimension $\ge2$ in $M_x$. Thus the projection $r_\sU:\sU\to \Gr(4, n+1)$ induces an isomorphism $\Pic(\Gr(4, n+1))\to \Pic(\sU)$. This implies (1) by standard descent theory.

The case $n=3$ is similar, except that $\Gr(4, n+1)$ is replaced with the point $\Gr(4,4)$, so $\Pic(\sU)=\{0\}$.
\end{proof}

\begin{remark}\label{rem:Cong1} It follows directly from our construction of $\sE_X$ and $\sF_X$ that
\[
\det\sE_X=\gamma=\det\sF_X.
\]
\end{remark}

For positive integers $m_1,\ldots, m_r,n$, define
\begin{equation}\label{eqn:BundleNotation}
\sE_{m_1,\ldots,m_r;n}:=\oplus_{i=1}^r\sE_{m_i,n}
\end{equation}

\begin{proposition}
\label{prop:relOrientationObservations}
Assume $\rnk(\sE_{m_1,\ldots,m_r;n})=\dim H_n$ and that the $m_i$ are all odd.
Then $\sE_{m_1,\ldots,m_r;n}$ is relatively oriented if and only if:
\begin{enumerate}
    \item The number of $m_i$'s such that $m_i\equiv-1\mod 4$ is even.
    \item When $n$ is odd, the number of direct summands $r$ in $\sE_{m_1,\ldots,m_r;n}$ is even.
    \item When $n$ is even, the number of direct summands $r$ in $\sE_{m_1,\ldots,m_r;n}$ is odd.
\end{enumerate}
\end{proposition}

\begin{proof}
Since $\rnk(\sE_{m_1,\ldots,m_r;n})=\sum_{i=1}^r(3m_i+1)$ and $\dim H_n=4n$, the condition
$\rnk(\sE_{m_1,\ldots,m_r;n})=\dim H_n$ implies $n\ge4$. Then  $\Pic(H_n)=\Z^3$ with generators $\sO_{H_n}(H^\ncm_n)$,  $\gamma:=\det \sF=\det\sE$, and $\sO_{\Gr(4,n+1)}(1)$. 
 
 Since we are assuming all the $m_i$ are odd, $\det\sE$ and $\det\sF$ appear in $\det\sE_{m_i,n}$ to an even power. $\sO_{H_n}(H^\ncm_n)$ occurs to an even power in $\det\sT_{H_n}$, and is even in $\det\sE_{m_i,n}$ if and only if $m_i\equiv 1\mod 4$ (for odd $m_i$). This proves (1). 

$\sO_{\Gr(4, n+1)}(1)$ occurs in $\det\sE_{m_i,n}$ to the power $M_{m_i}$, which is odd  if $m_i\cong 1\mod 4$ and even if $m_i\equiv -1\mod 4$

 $\sO_{\Gr(4, n+1)}(1)$ occurs in $\det\sT_{H_n}$ to the power $n+1$, so if $n$ is odd, there must be an even number of $m_i\equiv 1\mod 4$, and if $n$ is even there must be an odd number of $m_i\equiv 1\mod 4$. This proves (2) and (3).
 
The statement about the choice of orientation follows directly from the construction of the $\rho_{m_i,n}$.
  \end{proof}

For $n\ge4$, $m\ge2$,  $\sE_{m,n}$ is relatively oriented if an only if $n$ is even and  and $m\equiv 1\mod 4$. Assuming this to be the case, the isomorphisms $\rho_{m,n}$ of Proposition~\ref{prop:DetCMSectionBundle} and that of Lemma~\ref{lem:TangentDet} give rise to a relative orientation for $\sE_{m,n}$.

\begin{proof}  In all cases, Lemma~\ref{lem:TangentDet}  gives an isomorphism of $\det\sT_\Phi$ with a square and an isomorphism, $\det\sT_{H_n}\cong \Phi^*\sO_{\Gr(4, n+1)}(n+1)$ modulo a square.  

In case $m$ is odd, $\det\sE$ and $\det\sF$ appear in $\det\sE_{m,n}$  even powers. If $m\equiv1]\mod 4$ and  $n$ is even,  $M_m$  and  $n+1$ are both odd and $K_m$ is even. Thus  $\sO_{\Gr(4, n+1)}(1)$ appears to an odd power  in both formulas and the twist by $\sO_{H_n}(H^\ncm_n)$ is even in both formulas. 

If $m$ is not congruent to 1 modulo 4,   $\sO_{H_n}(H^\ncm_n)$ occurs to an even power in $\det\sT_{H_n}$ but to an odd power in $\det\sE_{m,n}$.  As  
$\sO_{H_n}(H^\ncm_n)$,  $\Phi^*\sO_{\Gr(4, n+1)}(1)$ and $\det \sF=\det\sE$ are independent generators of $\Pic(H_n)=\Z^3$,  $\sE_{m,n}$ is not relatively oriented in this case.

If $m\equiv 1\mod 4$ and $n$ is odd, then $\det\sE$ and $\det\sF$ appear in $\det\sE_{m,n}$  to even powers, but $K_m$ is odd and $n+1$ is even, so 
$\sE_{m,n}$ is not relatively oriented, by the same reasoning as in the previous case.
\end{proof}

\begin{remark}\label{rem:Cong2} If one of the $m_j$ is even, then $\sE_{m_j,n}$ has odd rank $3m_j+1$, so the Euler class $e(\sE_{m_j,n})\in H^{3m_j+1}(H_n, \sW(\det^{-1}\sE_{m_j,n}))$. By the multiplicativity of the Euler class, $e(\oplus_j\sE_{m_j,n})$ vanishes as well, so we can ignore discussing orientation conditions in this case.
\end{remark}

\begin{ex} 1. We consider the case $r=1$ $m_1=5$, $n=4$. Then on $H_4^\cm$, 
\[
\det\sE_{5,4}=\Phi^*\sO_{\Gr(4,5)}(35)\otimes \det\sE^{\otimes-20}\otimes\det\sF^{\otimes 10}
\]
and
\[
\det \sT_{H_4^\cm}=(\det\sF^{\otimes -3}\otimes\det\sE^{\otimes 2})^{\otimes 4}\otimes 
\det(p^\cm_*\sO_\Pi(1)))^{\otimes 6}\otimes \Phi^*\sO_{\Gr(4,5)}(5)
\]
Thus the vector bundle $\sE_{4,5}$ is relatively oriented on $H_4^\cm$ and has rank $3\cdot 5+1=16$, equal to the dimension $4\cdot 4$ of $H_4$.  \\[5pt]
2. For $r=1$, $m=9$, $n=7$, the condition $3m+1=4n$ is satisfied.   On $H_7^\cm$, 
\[
\det\sE_{9,7}=\Phi^*\sO_{\Gr(4,8)}(M_9)\otimes \det\sE^{\otimes-120}\otimes\det\sF^{\otimes 84}
\]
with $M_9=(3\cdot 9-1)\cdot 9/2=13\cdot 9$ odd, 
and
\[
\det \sT_{H_7^\cm}=(\det\sF^{\otimes -3}\otimes\det\sE^{\otimes 2})^{\otimes 4}\otimes 
\det(p^\cm_*\sO_\Pi(1)))^{\otimes 6}\otimes \Phi^*\sO_{\Gr(4,8)}(8)
\]
so $\sE_{9,7}$ is not relatively oriented. 
\end{ex}

\begin{remark}\label{rem:GEquivarianceProperty}The group scheme $\GL_{n+1}$ acts naturally on $\P^n$, and this action extends to give a canonical $\GL_{n+1}$-linearization on all ``natural'' sheaves on $\P^n$, such as $\sO_{\P^n}(m)$ and the tangent sheaf $\sT_{\P^n}$. The $\GL_{n+1}$-action on $\P^n$ induces a $\GL_{n+1}$-action on the diagrams \eqref{eqn:MainDiag}, \eqref{eqn:MainDiag}, and \eqref{eqn:MainDiag}, where $\GL_{+1}$ acts on $\sU$ via its action on linear forms $H^0(\P^n, \sO_{\P^n}(1))$. This extends to a $\GL_{n+1}$ linearization of all ``natural'' sheaves on the various schemes in these diagrams, such as $\sT_{H_n}$, $\sT_{\Gr(4,n+1)}$, $\sE_X$, $\sF_X$, and so on. As the divisor $H^\ncm_n$ is fixed by  the $\GL_{n+1}$-action, we also an induced $\GL_{n+1}$-action on $H^\ncm_n$ and $H^\cm_n$, and a linearization of $\sO_{H_n}(d\cdot H^\ncm_n)$. 

The  resolution of $j^*\sT_\Phi$ carries a natural $\GL_{n+1}$ linearization, so the description of $\det\sT_{H_n}$  given in  Lemma~\ref{lem:TangentDet}  is $\GL_{n+1}$-equivariant. Similarly,  the complex \eqref{eqn:ResOnHn} carries a natural $\GL_{n+1}$ linearization and the description of 
$\det\sE_{m,n}$ in Proposition~\ref{prop:DetCMSectionBundle} is also $\GL_{n+1}$-equivariant.
\end{remark}

\section{Extending to $\Z[1/6]$}\label{sec:Integrality}

We   extend the picture to the closure  of the locus of twisted cubics over $\Z[1/6]$, $H_n/\Z[1/6]$. Our first task is to show that $H_n/\Z[1/6]$ is smooth over $\Z[1/6]$ and that the non-Cohen-Macaulay locus $H_n^\ncm/\Z[1/6]$ is a smooth divisor. The issue here is that in general, closure does not commute with base-change, so we need to check that the fiber over some $\F_p$ of the closure of $H^{sm}_{n, \Z[1/6]}$ in the Hilbert scheme is the same as $H_n/\F_p$.

For this, we first recall the notion and basic properties of a {\em Nagata ring}. All our rings are assumed commutative. A  noetherian domain $R$ with quotient field $K$ satisfies the property $N$-2 if for every finite field extension $K\subset L$, the integral closure of $R$ in $L$ is a finite $R$-module. A noetherian ring $R$ is a Nagata ring if for each prime ideal $P$, the domain $R/P$ satisfies $N$-$2$. An arbitrary localization of a Nagata ring is a Nagata ring and if $(f_1,\ldots, f_n)$ generate the unit ideal in $R$ and $R[1/f_i]$ is a Nagata ring for each $i$, then $R$ is a Nagata ring (\cite[Tag 032E, Lemma 10.162.6, Lemma 10.162.7]{Stacks}. A a scheme $X$ is a {\em Nagata scheme} if each point $x\in X$ admits an affine open neighborhood $U$ such that $\sO_X(U)$ is a Nagata ring (see \cite[Tag 033R, Definition 28.13.1]{Stacks}). $X$ is a Nagata scheme if and only if $X$ admits an  open cover by Nagata schemes, if and only if for each affine open $U\subset X$, $\sO_X(U)$ is a Nagata ring \cite[Tag 033R, Lemma 28.13.6]{Stacks}. In addition, a fundamental result of Nagata states that if $R$ is a Nagata ring and $A$ is a finitely generated $R$-algebra, then $A$ is also a Nagata ring (see \cite[Theorem 72]{Matsumura} or \cite[Tag 032E, Proposition 10.162.15]{Stacks}). Thus, if $X$ is a Nagata scheme and $f:Y\to X$ is a morphism, locally of finite type, then $Y$ is also a Nagata scheme. Finally, Nagata has shown that a Noetherian complete local ring is a Nagata ring (\cite[Tag 032E, Lemma 10.162.8]{Stacks}, \cite[Corollary 2, Chap. 12]{Matsumura}). Besides the Noetherian complete local rings,  a Dedekind domain with characteristic zero quotient field (such as $\Z$) is a Nagata ring. 

We present a version of a lemma of Hironaka \cite[Lemma 4]{Hironaka}.

\begin{lemma} Let $\sO$ be a discrete valuation ring with  residue field $k$ and quotient field $K$, and let $X\to \Spec \sO$ be an integral, finite-type $\sO$ scheme, flat over $\sO$. Let $\sP\subset \sO_X$ be the ideal sheaf of the reduced special fiber  $(X\times_\sO k)_\red$. Suppose that $\sO$ is a Nagata ring and  that  both the generic fiber $X\times_\sO K$ and the reduced special fiber  $(X\times_\sO k)_\red$ are normal and integral.  Suppose in addition that, for $x\in X$ the  generic point of $X\times_\sO k$, we have $\sP_x=(\pi)\sO_{X,x}$.  Then $X$ is normal.
\end{lemma}

\begin{proof} Our proof is taken directly from  Hironaka {\em loc. cit.} with minor changes. We note that the assumption that $X\times_\sO K$ is integral follows from our assumption that $X$ is integral. Choose a parameter $\pi$ for $\sO$.

The conclusion is local on $X$, so we can replace $X$ with $\Spec A$, with $A$ a Noetherian local $\sO$-algebra essentially of finite type over $\sO$, such that\\[5pt]
i. $(\pi)A$ has a unique minimal prime $P$ and   $(\pi)A_P=PA_P$\\[2pt]
ii. $A/P$ is normal\\[2pt]
iii. $A[1/\pi]$ is normal. \\[5pt]
Since $\sO$ is a Nagata ring, so is $A$. Since $A_P$ is a Noetherian local domain with principal maximal ideal $(\pi)A_P$, $A_P$ is a discrete valuation ring, hence integrally closed.   Let $A'=A_P\cap A[1/\pi]$, the intersection taking place in the quotient field  $F$ of $A$. By (iii), $A'$ is equal to the integral closure of $A$ in  $F$. 

 Let $P'$ be a minimal prime of $(\pi)A'$. Then $P'\cap A=P$. Since $A_P$ is integrally closed,
 $A'_P=A_P$,  so $P'$ is the unique prime of $A'$ lying over $P$ and $A'_{P'}=A_P$. By (i), we have $(\pi)A'_{P'}=P'A'_{P'}$ and $P'$ is the unique minimal prime of $(\pi)A'$.  Since $A'$ is normal, $(\pi)A'$ has no embedded components \cite[Theorem 38]{Matsumura}. Thus   $P'=(\pi)A'$ and $A'/(\pi)=A'/P'$. But  after localizing at $P$ we have 
\[
(A/P)_P=A_P/PA_P=A'_{P'}/P'A_{P'}=(A'/P')_{P'}
\]
so the domains $A/P$ and $A'/P'$ have the same quotient fields. By (ii), we have 
\[
A/P= A'/P'=A'/(\pi)A'.
\]

Since $A$ is a Nagata ring, $A'$ is a finite $A$-module. But $A'[1/\pi]=A[1/\pi]$ by (iii),  so there is an integer $e\ge0$ with $\pi^eA'\subset A$. We take $e$ minimal and want to show that $e=0$, that is, $A=A'$, which will complete the proof.

If $e>0$, then $A'/(\pi)A'=A/P$ implies $A'= A+\pi\cdot A'$. Then 
\[
\pi^{e-1}\cdot A'= \pi^{e-1}A+\pi^eA'\subset A, 
\]
contrary to the minimality of $e$. 
\end{proof}

\begin{corollary} $H_n/\Z[1/6]$ is a smooth projective scheme over $\Z[1/6]$, and for each prime $p>3$, the base-change $(H_n/\Z[1/6])\times_{\Z[1/6]}\F_p$ is isomorphic to $H_n/\F_p$. \end{corollary}

\begin{proof} Let  $p>3$ be a prime and let $\sO:=\Z_{(p )}$, with quotient field  $\Q$ and residue field $k=\F_p$; $\sO$ is a Nagata ring.  It suffices to show that $H_n/\sO$ is smooth over $\sO$. $H_n/\sO$ is the closure of $H_n/\Q$ in $\Hilb^{3m+1}(\P^n_\sO)$, so $H_n/\sO$ is integral and is flat and projective over $\sO$. For each $x\in H_n/\Q$, the corresponding curve $C_x$ is geometrically connected, so the same holds for all $x\in H_n/\sO$. By the classification given in the proofs of \cite[Lemma 1, Lemma 2, Lemma 4]{PS} and our determination of the universal deformation space in Proposition~\ref{prop:DefSpace}, there 
 is a neighborhood $U$ of $H_n/k$ in $\Hilb^{3m+1}(\P^n_k)$ such that, for $x$ in  $U\setminus H_n/k$,  the curve $C_x$ is one of the following two types:\\[5pt]
a. $C_x$ is not connected.\\[2pt]
b. $C_x^{\red}$ is a plane cubic and the kernel of $\sO_{C_x}\to \sO_{C_x^{\red}}$ is non-zero and supported at a smooth point $p$ of $C_x^{\red}$.\\[5pt]
Thus, if  $H_n/\sO\times_\sO k\setminus H_n/k\neq\0$, there is a point $x\in H_n/\sO\times_\sO k\setminus H_n/k$  with corresponding curve $C_x$ of type (b).  Since $H_n/\sO$ is the closure of the locus of twisted cubics, there is a discrete valuation ring $\sO'$ and  closed subscheme $\sC\subset \P^n_{\sO'}$ defining a flat family over $\sO'$ with generic fiber $\sC_\eta$ a  smooth 
twisted cubic and special fiber $\sC_0$ our curve $C_x$ of type (b), with its corresponding smooth point $p$. We may assume that $\sO'$ is complete, hence a Nagata ring. But then there is a neighborhood $V$ of $p$ in $\sC$ that is an integral scheme, flat over $\sO'$ with smooth generic fiber and with $V_0\setminus \{p\}$ and $V_0^{\red}$ both smooth and integral; here $V_0$ is the special fiber of $V\to \Spec \sO'$. By Hironaka's lemma, $V$ is a normal scheme. But as $V_0\setminus \{p\}$ is reduced, and there are no embedded components in $V_0$ (\cite[Theorem 38]{Matsumura} again)  the special fiber $V_0$ is reduced, contrary to the assumption that $\sC_0=C_x$ has an embedded component at $p$.

Thus, set-theoretically, the fiber $H_n/\sO\times_\sO k$ is equal to $H_n/k$. The dense open subscheme $H^{sm}_{n, \sO}$  of twisted cubics is smooth over $\sO$ with special fiber $H^{sm}_{n, k}$ dense in $H_n/k$ and $H_n/\Q$ is smooth over $\Q$. Applying Hironaka's lemma again, it follows that $H_n/\sO$ is normal. As above, this implies that $H_n/\sO\times_\sO k=H_n/k$ and so $H_n/\sO$ is smooth over  $\sO$.

\end{proof}

\begin{remark} In \cite{HSS}, Heinrich, Skjelnes and Stevens extend  the smoothness result of Piene-Schlessinger to fields of characteristic $2$ and $3$, thus $H_3/k$ is a smooth $k$-scheme for all  fields $k$. However, we do not have a nice classification of the curves in $\Hilb^{3m+1}(\P^3_k)\setminus H_3$, needed for the argument above.  We also do not know if $H_n/k$ is smooth for fields of characteristic $2,3$ and $n>3$.
\end{remark}

For an arbitrary base-scheme $S$, we let $H_n^\cm/S\subset H_n/S$ denote the subscheme of $H_n/S$ parametrizing flat families of Cohen-Macaulay curves. 

\begin{lemma}\label{lem:HCM} 1. $H_n^\cm/\Z[1/6]$ is an open subscheme of $H_n/\Z[1/6]$.\\[2pt]
2.  For $f:T\to S$ a morphism of schemes we have $H_n^\cm/S\times_ST = H_n^\cm/T$, where we view both schemes as subschemes of $\Hilb^{3m+1}(\P^n_T)$ via the canonical isomorphism $\Hilb^{3m+1}(\P^n_T)\cong \Hilb^{3m+1}(\P^n_S)\times_ST$.
\end{lemma}

\begin{proof} (1) Let $x\in H_n^\cm/\Z[1/6]$ be a geometric point and let $C_x\subset \P^n_{k(x)}$ be the corresponding curve. It suffices to show that an open neighborhood of $x$ in $\Hilb^{3m+1}(\P^n_{\Z[1/6]})$ is contained in $H_n^\cm/\Z[1/6]$. 

Since $C_x$ is a degree three curve, we have the following possibilities for the 1-cycle $|C_x|$ associated to the scheme $C_x$; in all cases $C_x$ is connected and is not contained in a plane.\\[5pt]
a. $C_x$ is irreducible and  $|C_z|=1\cdot C_x^\red$\\[2pt]
b. $C_x^\red$ is a union of a line $\ell$ and a smooth conic $Q$, and $|C_z|=1\cdot \ell+1\cdot Q$.\\[2pt]
c1. $C_x^\red$ is a union of three distinct lines $\ell_i$ and  $|C_z|=\sum_{i=1}^31\cdot \ell_i$. The three lines do not pass though a common point.\\[2pt]
c2. $C_x^\red$ is a union of three distinct lines $\ell_i$ and  $|C_z|=\sum_{i=1}^31\cdot \ell_i$. The three lines   pass though a common point.\\[2pt]
d. $C_x^\red$ is a union of two distinct lines $\ell_i$ and  $|C_z|=1\cdot \ell_1+2\cdot \ell_2$\\[2pt]
e. $C_x^\red$ is a   line $\ell$ and  $|C_z|=\ 3\cdot \ell$.\\[5pt]
In cases (a), $C_x$ is a smooth cubic rational curve. In cases (b) and (c1), $C_x$ is a local complete intersection. In case (c2), $C_x$ is projectively equivalent to the union of three coordinate axes through $(0,0,0,1)$, defined by the ideal $(X_0X_1, X_0X_2, X_1X_2)$. In case (d), $C_x$ is projectively equivalent to the curve defined by the ideal $(X_0^2,X_2)\cap (X_0,X_1)=(X_0^2, X_0X_2, X_1X_2)$. In case (e), $C_x$ is projectively equivalent to the curve defined by the ideal $(X_0^2,X_0X_1, X_1^2)$.

In the local complete intersection cases (a), (b) and (c1),  the cotangent complex $\L_{C_x/\P^3_{k(x)}}$ is isomorphic in $D(\Coh_{\sO_{C_x}})$ to $\sI_{C_x}/\sI_{C_x}^2[1]$.

In cases (c2), (d) and (e), one can make an explicit computation of $\tau_{\ge-2}\L_{C_x/\P^3_{k(x)}}$. In fact, if a Cohen-Macaulay curve $C_x$ has an Eagon-Northcott complex with   matrix $(L_{ij})$, let $(L_{ij}^*)$ be the matrix
\[
\begin{pmatrix} L_{12}&L_{22}&L_{32}\\-L_{11}&-L_{21}&-L_{31}\end{pmatrix}
\]
This gives the map $\sO_{\P^3_{k(x)}}(-4)^3\xrightarrow{(L_{ij}^*)}\sO_{\P^3_{k(x)}}(-3)^2$; let 
$\overline{\sO}_{\P^3_{k(x)}}(-3)^2$ denote the cokernel of $(L_{ij}^*)$. Then 
\[
\tau_{\ge-2}\L_{C_x/\P^3_{k(x)}}\cong [\overline{\sO}_{\P^3_{k(x)}}(-3)^2\xrightarrow{(L_{ij})}
\sO_{C_x}(-2)^3]
\]
with $\sO_{C_x}(-2)^3$ in degree -1. To see this, the complex described above is the Lichtenbaum-Schlessinger cotangent 
complex associated to the Eagon-Northcott resolution of $\sO_{C_x}$, and the Lichtenbaum-Schlessinger cotangent complex computes the truncation $\tau_{\ge-2}$ of the true cotangent complex. 

Choosing a suitable matrix $(L_{ij})$ in the cases (c2), (d) and (e), one can show using the explicit complex above that $\sH^{-2}(\L_{C_x/\P^3_{k(x)}})=0$. Using the distinguished triangle
\[
\L_{\P^3_{k(x)}/\P^n_{\Z[1/6]}}\otimes\sO_{C_x}\to \L_{C_x/\P^n_{\Z[1/6]}}\to \L_{C_x/\P^3_{k(x)}}\xrightarrow{+1}
\]
and the fact that $\L_{\P^3_{k(x)}/\P^n_{\Z[1/6]}}\otimes\sO_{C_x}$ is concentrated in degree -1, we see that $\sH^2( \L_{C_x/\P^n_{k(x)}})=0$.  Thus, in all cases of a Cohen-Macaulay cubic rational curve $C_x$, the obstruction space for the deformations of $C_x$ in $\P^n_{\Z[1/6]}$  is $H^1(C_x, \sN_{C_x})$, with $\sN_{C_x}=\sHom(\sI_{C_x}, \sO_{C_x})$ the normal sheaf of $C_x$ in $\P^n_{\Z[1/6]}$. 

In \cite[Lemma 5]{PS} it is shown that $H^1(C_x, \bar{\sN}_{C_x})=0$, where $\bar{\sN}_{C_x}$ the normal sheaf in $\P^3_{k(x)}$. We have the exact sequence
\[
0\to \sO_{C_x}(1)^{n-3}\oplus \sO_{C_x}^\epsilon \to \sN_{C_x}\to \bar{\sN}_{C_x}\to 0
\]
with $\epsilon=1$ if $k(x)$ has positive characteristic and $\epsilon=0$ if $k(x)$ has   characteristic zero. Thus   $H^1(C_x, {\sN}_{C_x})=0$ and the deformations of $C_x$ are unobstructed.

By computing $H^0(C_x, \sN_{C_x})$, one can similarly show that the first order deformations are generated by first order deformations of the matrix $(L_{ij})$ in the Eagon-Northcott complex for $C_x$, which shows that all the nearby deformations of $C_x$ are in $H^\cm_n/\Z[1/6]$. This shows that $H^\cm_n/\Z[1/6]$ is open in $H_n/\Z[1/6]$.

Let $\pi:S\to \Spec \Z[1/6]$ be a $\Z[1/6]$-scheme. As the Eagon-Northcott complex remains exact after any base-change that preserves the codimension of a Cohen-Macaulay curve $C_x$ in a $\P^3_{k(x)}\subset \P^n_{k(x)}$,  it follows that $H^\cm_n/\Z[1/6]\times_{\Z[1/6]}S$ is contained in $H^\cm_n/S$. Conversely, for each geometric point $x$ of 
  $H^\cm_n/S$, the corresponding  geometric point of $\Hilb^{3m+1}(\P^n_{\Z[1/6]})$ is in $H^\cm_n/\Z[1/6]$, so $H^\cm_n/\Z[1/6]\times_{\Z[1/6]}S=H^\cm_n/S$. This proves (2).
\end{proof}

Working over $\Z[1/6]$, use a parallel notation to the case of the base-field $k$. We have the Grassmann variety $\Gr_{\Z[1/6]}(4, n+1)$, the morphism $\Phi:H_n/\Z[1/6]\to \Gr_{\Z[1/6]}(4, n+1)$,   the universal 3-plane bundle $\Pi_{\Z[1/6]}\to \Gr_{\Z[1/6]}(4, n+1)$, the locally free sheaf $\sV_{\Z[1/6]}:=p_{1*}\sO_{\Pi_{\Z[1/6]}}(2)$, the Grassmann bundle $r:\Gr(3,\sV_{\Z[1/6]})\to \Gr_{\Z[1/6]}(4, n+1)$ and the morphism $q:H_n/\Z[1/6]\to \Gr(3,\sV_{\Z[1/6]})$. This gives us the (reduced) image subscheme $X_{\Z[1/6]}:=q(H_n/\Z[1/6])$ and the closed subscheme $F_{\Z[1/6]}:=q(H_n^\ncm/\Z[1/6])\subset X_{\Z[1/6]}$. Let $\sU_{\Z[1/6]}$ be the $\Z[1/6]$-scheme of $3\times2$ matrices $L:=(L_{ij})$ of linear forms  on $\Pi_{\Z[1/6]}$ such that $\<e_{12}(L), e_{13}(L), e_{23}(L)\>$ has rank 3 in $\sV_{\Z[1/6]}$. We let $\Gamma_{\Z[1/6]}$ be the group-scheme over $\Z[1/6]$ $\GL_3\times\GL_2/\G_m$, with the action $(g,h)\cdot L=g\cdot L\cdot h^{-1}$ on $\sU_{\Z[1/6]}$.

\begin{lemma}\label{lem:IntegralStructure} 1. $F_{\Z[1/6]}\subset X_{\Z[1/6]}$ are smooth over $\Z[1/6]$\\[2pt]
2. The map $q:H_n/\Z[1/6]\to X_{\Z[1/6]}$ lifts to an isomorphism of $H_n/\Z[1/6]$ with the blow-up $\Bl_{F_{\Z[1/6]}}X_{\Z[1/6]}$.\\[2pt]
3. For a $\Z[1/6]$-algebra $R$, sending $L\in \sU_{\Z[1/6]}(R)$ to $\<e_{12}(L), e_{13}(L), e_{23}(L)\>\subset \sV_{\Z[1/6]}(R)$ defines a smooth morphism  $\psi:\sU_{\Z[1/6]}\to X$, making $\sU_{\Z[1/6]}\to X_{\Z[1/6]}$ a $\Gamma_{\Z[1/6]}$-torsor over $X$ (via the given action of $\Gamma$ on $\sU$).
\end{lemma}

\begin{proof} For (3), we consider the action map $\rho: \Gamma_{\Z[1/6]}\times_{\Z[1/6]}\sU_{\Z[1/6]}\to \sU_{\Z[1/6]}\times_{\Z[1/6]}\sU_{\Z[1/6]}$. As $\rho$ is an isomorphism after taking $-\times_{\Z[1/6]}\F_p$ for all $p>3$, it follows that $\rho$ is an isomorphism. This proves (3) and shows that $X_{\Z[1/6]}$ is smooth over $\Z[1/6]$. Similarly, since $H_n^\ncm/\Z[1/6]$ is smooth over $\Z[1/6]$ and for each $p$ $H_n^\ncm/\F_p\to F_{\F_p}$ is a projective space bundle, it follows that $F_{\Z[1/6]}$ is smooth over $\Z[1/6]$.

Since $H_n/\F_p\cong \Bl_{F_{\F_p}}X_{\F_p}$ for $p>3$, it follows that the closed subscheme $q^{-1}(F_{\Z[1/6]})\subset H_n/\Z[1/6]$ is flat over $\Z[1/6]$, as each local section of the ideal of the fiber over $\F_p$ lifts to a local section of the ideal over $\Z[1/6]$. Thus $q^{-1}(F_{\Z[1/6]})$ is a Cartier divisor on $H_n/\Z[1/6]$, and hence the map $q:H_n/\Z[1/6]\to X_{\Z[1/6]}$ lifts to a map $\tilde{q}:H_n/\Z[1/6]\to \Bl_{F_{\Z[1/6]}}X_{\Z[1/6]}$. As both source and target are smooth over $\Z[1/6]$, and $\tilde{q}$ is an isomorphism on each fiber over $\F_p$, $\tilde{q}$ is an isomorphism.
\end{proof}
 
 We thus may define the sheaves $\sF_{X_{\Z[1/6]}}$, $\sE_{X_{\Z[1/6]}}$ on $X_{\Z[1/6]}$ by $\Gamma_{\Z[1/6]}$-descent, as used to define $\sF_{X_k}$, $\sE_{X_k}$; we let $\sF_{\Z[1/6]}$, $\sE_{\Z[1/6]}$ be the respective pull-backs of  $\sF_{X_{\Z[1/6]}}$, $\sE_{X_{\Z[1/6]}}$ to $H_n/\Z[1/6]$. We have the universal curve $p:\sC_n/\Z[1/6]\to H_n/\Z[1/6]$ and the
  locally free sheaf $\sE_{m, n, \Z[1/6]}:=p_*(\sO_{C_n/\Z[1/6]}(m))$.    
 \begin{theorem}\label{thm:IntegralExtension} 1. We have a canonical isomorphism
 \begin{multline*}
 \det \sE_{n,m,\Z[1/6]}\cong  \Phi^*\sO_{\Gr_{\Z[1/6]}(4,n+1)}(M_{m,n})\otimes \det\sE_{\Z[1/6]}^{\otimes-\scriptscriptstyle\begin{pmatrix}m+1\\3\end{pmatrix}}\\\otimes\det\sF_{\Z[1/6]}^{\otimes\scriptscriptstyle\begin{pmatrix}m\\3\end{pmatrix}}\otimes\sO_{H_n/\Z[1/6]}(K_mH_n^\ncm/\Z[1/6]).
\end{multline*}
2. We have a canonical isomorphism
 \begin{multline*}
\det\sT_{H_n/\Z[1/6]}\cong (\det\sF_{\Z[1/6]}^{\otimes -3}\otimes\det\sE_{\Z[1/6]}^{\otimes 2})^{\otimes 4}\otimes 
\det(\pi_{1*}\sO_{\Pi_n}(1)))^{\otimes 6}\\\otimes \sO_{H_n/\Z[1/6]}(-6\cdot H_n^\ncm/\Z[1/6]) \otimes\Phi^*\sO_{\Gr_{\Z[1/6]}(4,n+1)}(n+1).
\end{multline*}
\end{theorem}

\begin{proof} Using Lemma~\ref{lem:IntegralStructure},  we can repeat the constructions used in \S\ref{sec:Det} to yield the computation of the determinants. 
\end{proof}

\begin{remark}  Just as in the case of a field, the construction of $\sE_{X_{\Z[1/6]}}$ and  
$\sF_{X_{\Z[1/6]}}$ by $\Gamma_{\Z[1/6]}$-descent gives a canonical isomorphism $\det\sE_{\Z[1/6]}\cong \det\sF_{\Z[1/6]}$. Thus the congruences in the $m_j, n$ needed for  $\oplus_j \sE_{n,m_j,\Z[1/6]}$ to be relatively oriented are the same as in the case of a field (see Remark~\ref{rem:Cong2}).
\end{remark}

\begin{corollary}\label{cor:Integrality} Let  $m_1,\ldots, m_r, n$ be integers with $\sum_j3m_j+1 = 4n$, $n\ge4$, $m_j\ge2$. Suppose that the congruences in the $m_j$ and $n$ as given by   Remark~\ref{rem:Cong2} are satisfied. Then
 there is an element $\Deg_{\Z[1/6]}e(\oplus_j\sE_{n,m_j, \Z[1/6]})\in \GW(\Z[1/6])$ such that for each field $k\supset \Z[1/6]$, the push-forward of the Euler class $e(\oplus_j\sE_{n,m_j,k})$ is a well-defined element of $\GW(k)$ and is equal to the image of $\Deg_{\Z[1/6]}e(\oplus_j\sE_{n,m_j, \Z[1/6]})\in \GW(\Z[1/6])$  under the base-change map $\GW(\Z[1/6])\to \GW(k)$. Here 
 $e(\oplus_j\sE_{n,m_j,k})$ is the Euler class in hermitian $K$-theory.
\end{corollary}

\begin{proof} We use hermitian $K$-theory to define the Euler classes; this is a well-defined construction over $\Z[1/2]$. The condition $\sum_j3m_j+1 = 4n$ says that $\oplus_j\sE_{n,m_j, k}$ has rank $\sum_j3m_j+1$ equal to the dimension $4n$ of $H_n/k$, and the congruences being satisfied implies that 
$\oplus_j\sE_{n,m_j, k}$ is relatively oriented. Thus the Euler class $e(\oplus_j\sE_{n,m_j,k})$ lives in the twisted hermitian $K$-theory $KQ^{8n, 4n}(H_n/k, \omega_{H_n/k})$, and as $H_n$ is smooth and proper of dimension $4n$ over $k$, we have the well-defined push-forward map
\[
p_{H_n*}: KQ^{8n, 4n}(H_n/k, \omega_{H_n/k})\to KQ^{0, 0}(k)=\GW(k)
\]
Similarly, we have the Euler class $e(\oplus_j\sE_{n,m_j, \Z[1/6]})\in  KQ^{8n, 4n}(H_n/\Z[1/6], \omega_{H_n/\Z[1/6]})$ and the push-forward
\[
p_{H_n/\Z[1/6]*}: KQ^{8n, 4n}(H_n/ \Z[1/6], \omega_{H_n/ \Z[1/6]})\to KQ^{0, 0}(\Z[1/6])=\GW(\Z[1/6])
\]
Defining 
\[
\Deg_{\Z[1/6]}e(\oplus_j\sE_{n,m_j, \Z[1/6]}):=p_{H_n/\Z[1/6]*}(e(\oplus_j\sE_{n,m_j, \Z[1/6]}))
\]
and using the commutativity  of ${}_*$ and ${}^*$ in cartesian squares gives the result.
\end{proof} 

\begin{lemma}\label{lem:Coker} 1. Let $d\in \Z$ be a square-free even integer and let $k$ be a field of characteristic prime to $d$ (or of characteristic zero), giving the ring homomorphism $\Z[1/d]\to k$. For $x\in \GW(\Z[1/d])$, the image of $x$ in $\GW(k)$ is determined by the rank $r$ and signature $s$ of $x$, in other words
\[
x=s+\frac{r-s}{2}\cdot H\in \GW(k)
\]
where $H$ is the hyperbolic form.\\[2pt]
2. For a unit $u$ in a commutative ring $R$, let $\<u\>$ denote the rank one quadratic form $x\mapsto ux^2$. Then the cokernel of the map $\GW(\Z)\to \GW(\Z[1/6])$ is additively generated by the forms $\<2\>-1$, $\<3\>-1$ and $\<6\>-1$.
\end{lemma}

\begin{proof} We are indebted to R. Parimala for the statement and proof of (1); as (1) in case $d=6$ follows from (2), and we will only be using (1) in case $d=6$, we give an independent proof of (2), inspired by Parimala's argument.

Since natural map $\GW(\Z)\to \GW(\R)$ is an isomorphism \cite[Chap. 5, Theorem 4.1]{Scharlau}, and an element of $\GW(\R)$ is determined by its rank and signature (this is Sylvester's theorem of inertia), we need to show that an element $x\in \GW(\Z[1/6])$ of rank zero and signature zero is in the subgroup generated by $x_2:=\<2\>-1$, $x_3:=\<3\>-1$ and $x_6:=\<6\>-1$..  Let $x_1=0$. 

The discriminant $d_x$ of $x$ is of the form $\pm2^a3^b$ mod squares (in $\Q^\times$). The assumption that $x$ has zero rank and signature easily implies that $d_x>0$, so $d_x$ is in $\{1, 2,3,6\}$. Subtracting $x_{d_x}$, we may assume that $x$ has trivial discriminant.  

The element $x_6-x_2-x_3$ is the Pfister form $\<\<3,2\>\>$, hence has trivial discriminant and has Hilbert symbol $\{3,2\}=-1$ at $p=3$. Adding  $x_6-x_2-x_3$  to $x$ if necessary, we may assume that  $x$ has trivial Hilbert symbol +1 at $p=3$. Since by assumption $x$ has trivial Hilbert symbol at all $p>3$, including $p=\infty$, the Hilbert reciprocity theorem \cite[Chap. 5, Theorem 5.1]{Scharlau} implies that  $x$ has trivial Hilbert symbol at $p=2$. Thus $x$ has trivial rank, signature, and discriminant, as well as trivial Hilbert symbol at all $p$, so $x=0$ in $\GW(\Q_p)$ for all $p\le \infty$ by \cite[Chap. 5, Theorem 4.2]{Scharlau}. Thus $x=0$ in $\GW(\Q)$ by the Hasse-Minkowski theorem. Since the map $\GW(\Z[1/6])\to \GW(\Q)$ is injective \cite[Chap. 4, Corollary 3.3]{HusMilnor}, this says that $x=0$.  
\end{proof}

\begin{remark}\label{rem:Refinement}  The methods of \cite{LevineAB} use Euler classes with values in the twisted cohomology of the sheaf of Witt groups, rather than in hermitian $K$-theory. However, in the case of a locally free, relatively oriented sheaf $\sV$ on a smooth projective $k$-scheme $p:X\to \Spec k$ (with $\Char k\neq 2$) and with $\rnk(\sV)=\dim_kX$, the push-forward of the Euler class $e^{\KQ}(\sV)$ in hermitian $K$-theory and the push-forward of the Euler class $e^\sW(\sV)$ in Witt cohomology both have the same image in the Witt group $W(k)$. Indeed, Witt theory $W(-)$ is gotten from hermitian $K$-theory by inverting the algebraic Hopf map $\eta$, and by \cite[Proposition 3.4, Theorem 9.1]{LevineAspects}, the Witt-theory Euler class $e^W(\sV)$ and the Witt cohomology Euler class $e^\sW(\sV)$ are unchanged if one twists $\sV$ by a line bundle. Thus, we may assume that $\sV$ admits a section $s$ with zero locus $Z$ of dimension zero. In this case, both $p_*(e^W(\sV))$ and $p_*(e^\sW(\sV))$ are given by the (finite) push-forward of the respective Euler classes with support in $Z$. But as the Witt groups of $Z$ agree with the 0th cohomology of the Witt sheaf on $Z$, the Euler classes $e^W_Z(\sV)$ and $e^\sW_Z(\sV)$ are equal in this common group. Thus the conclusion of Corollary~\ref{cor:Integrality} may be applied to the quadratic degree of the Euler class
$e^\sW(\oplus_j\sE_{n,m_j,k})$ computed using Witt-sheaf cohomology, in other words, there is an element $d_{n,m_*}\in \GW(\Z[1/6])$ such that $p_*(e^\sW(\oplus_{j=0}^r\sE_{n,m_j,k})))$ is the image of $d_{n,m_*}$ in $W(k)$, as long as $k$ of characteristic $\neq 2,3$.  In addition, there are integers $a,b,c$ (depending on $n, m_1,\ldots, m_r$) such that 
\[
p_*(e^\KQ(\oplus_{j=0}^r\sE_{n,m_j,k})))=s+\frac{r-s}{2}\cdot H+a(\<2\>-1)+b(\<3\>-1)+c(\<6\>-1)
\]
where $s$ is the signature of $d_{n,m_*}$, $r$ is the rank of $d_{n,m_*}$ and $H$ is the hyperbolic form.
\end{remark}

\section{Counting twisted cubics}

Let $\sE_{m_1,\ldots,m_r;n}:=\bigoplus_{i=1}^r\sE_{m_i}\rightarrow H_n$. By \cite[$\S3$]{ES} a hypersurface in $\mathbb{P}^n$ defined by a degree $m$ homogeneous polynomial defines a section of the bundle $\sE_m\rightarrow H_n$ and the zeros of the section correspond to the twisted cubics contained in the hypersurface. So if $\operatorname{rank}\sE_{m_1,\ldots,m_r;n}=\operatorname{dim}H_n$, a general section of $\sE_{m_1,\ldots,m_r}$ has finitely many zeros and thus there are finitely many twisted cubics on a general complete intersection of multidegree $(m_1,\ldots,m_r)$. Classically, one can compute this finite number of twisted cubics as the degree of the top chern class (or Euler class) of $\sE_{m_1,\ldots,m_r;n}$. We aim to quadratically refine this count by computing the degree of $e^\sW(\sE_{m_1,\ldots,m_r;n})$ valued in $\operatorname{W}(k)$. 

Let $N=N_{\SL_2}$ be the normalizer of the diagonal torus in $\SL_2$.
We use the quadratically refined version of Bott's residue formula \cite[Theorem 9.5]{LevineAB} which applies to schemes with an $N$-action with finitely many fixed points. 
We recall a description of the irreducible $N$-representations and their Euler classes from \cite[\S6, Theorem 7.1]{LevineWitt}, and then define the action of $N$ on $H_n$ and find all the fixed points.
From now on we use the notation $e$ for the Euler class $e^{\sW}$ in Witt sheaf cohomology and $e^N$ for the equivariant Euler class.

\subsection{Irreducible $N$-representations and their equivariant Euler classes}

$N$ is generated by $t:=\begin{pmatrix}t & 0\\0 & t^{-1}  \end{pmatrix}$ and $\sigma:=\begin{pmatrix}0 & 1\\ -1 & 0\end{pmatrix}$.
For an integer $a>0$ let $\rho_a:N\rightarrow \GL_2$ the irreducible $N$-representation that sends $t$ to $\begin{pmatrix}t^a& 0\\
0& t^{-a}\end{pmatrix}$ and $\sigma$ to $\begin{pmatrix}
0 & 1\\
(-1)^a& 0
\end{pmatrix}$
and let $\rho^-_a$ be the $2$-dimensional $N$-representation that sends $t$ to $\begin{pmatrix}t^a& 0\\
0& t^{-a}\end{pmatrix}$ and $\sigma$ to $-\begin{pmatrix}
0 & 1\\
(-1)^a& 0
\end{pmatrix}$.
Together with the trivial representation which we call $\rho_0$ and the $1$-dimensional representation $\rho_0^-$ which sends $t$ to $1$ and $\sigma$ to $-1$ in $\G_m=\GL_1(k)$, these are all the irreducible $N$-representations.
The corresponding vector bundles on $BN$ are denoted by $\widetilde{\mathcal{O}}(a)$ and $\widetilde{\mathcal{O}}^-(a)$.
Let $p:N\rightarrow \SL_2$ be the inclusion and let $e\in H^2(\BSL_2, \sW)$ be the Euler class of the bundle on $B\SL_2$ defined by the fundamental representation of $\SL_2$ on $\A^2$ in $H^2(B\SL_2,\mathcal{W})$. Let $\gamma$ be the invertible sheaf on $BN$ defined by the representation $\rho_0^-$ and let $\tilde{e}\in H^2(BN, \sW(\gamma))$ be the Euler class of $\tilde{\sO}(2)$.

The equivariant Euler classes of $\widetilde{\mathcal{O}}^{\pm1}(a)$ are computed in \cite[Theorem 7.1]{LevineWitt}. 
\begin{theorem}
\label{thm:EulerClasses}
For $a$ odd, define 
\[\epsilon(a):=\begin{cases} 1 & \text{for $a\equiv 1\mod 4$}\\
-1 & \text{for $a\equiv 3\mod 4$.}\end{cases}\]
Then 
\[e(\widetilde{\mathcal{O}}(a))=\begin{cases} \epsilon(a)\cdot a\cdot p^*e\in H^2(BN,\mathcal{W})& \text{for $a$ odd}\\
\frac{a}{2}\cdot \widetilde{e}\in H^2(BN,\mathcal{W}(\gamma)) & \text{for $a\equiv2\mod 4$}\\
-\frac{a}{2}\cdot \widetilde{e}\in H^2(BN,\mathcal{W}(\gamma))& \text{for $a\equiv0\mod 4$.}
\end{cases}\]
Moreover, $e(\tilde{\sO}^-(a))=-e(\tilde{\sO}(a))$ and $\tilde{e}^2=4e^2$, this latter identity 
taking place in $H^4(BN,\mathcal{W})$. \end{theorem}

\begin{remark}\label{rem:OrientedBasis} The Euler class $e(V)$ (in Witt sheaf cohomology) of a rank $r$ locally free sheaf $V$ on some smooth scheme $X$ lives canonically in $H^r(X, \sW(\det^{-1}V))$. Our computation for $e(\tilde{\sO}(a))$ thus relies on a particular choice of isomorphism $\det\tilde{\sO}(a)\cong \sO_{BN}$ for $a>0$ odd and $\det\tilde{\sO}(a)\cong \gamma$ for $a>0$ even. These are given by taking as basis $e_1, e_2$ for the representation $\rho_a$ as described by the matrix formula above, that is $\rho_a(t)(e_1)=t^ae_1$, $\rho_a(t)(e_2)=t^{-a}e_2$, $\rho_a(\sigma)(e_1)=(-1)^ae_2$, $\rho_a(\sigma)(e_2)=e_1$, and using $e_1\wedge e_2$ as basis for $\det\rho_a$. We call  such as basis for $\rho_a$ {\em oriented}; concretely, two oriented bases $(e_1, e_2)$ and $(e_1', e_2')$ for $\rho_a$  are related by $e'_i=\lambda e_i$ for some $\lambda\in k^\times$, $i=1,2$. Thus the generator $e_1\wedge e_2$ for $\det\rho_a$ is well-determined up to a square in $k^\times$, so 
our computation of $e(\tilde{\sO}(a))$ is thus valid for any choice of oriented basis. 

The identity  $e(\tilde{\sO}^-(a))=-e(\tilde{\sO}(a))$ uses a basis $e_1, e_2$ for 
$\tilde{\sO}^-(a)$, with the action of $\rho_a(t)$ as before and with $\rho_a(\sigma)(e_1)=(-1)^{a+1}e_2$, $\rho_a(\sigma)(e_2)=-e_1$. The identity $\widetilde{e}^2=4p^*e^2$  uses
 the basis $(e_1\wedge e_2)^2$ for $(\det\rho_2)^{\otimes 2}$ and similarly for $(\det\rho_1)^{\otimes 2}$.
\end{remark}

\subsection{Action and fixed points}
For a smooth $k$-scheme $Y$ with $N$-action and an $N$-linearized vector bundle $V$ of rank $r$ on $Y$, we have the equivariant Euler class $e^N(V)$ in the $N$-equivariant Witt sheaf cohomology $H^{2r}_N(Y, \sW(\det^{-1}(V))$ \cite{LevineAB}.

\begin{lemma}
Let $a_1,\ldots,a_r$ be positive odd integers. Then 
\[N\rightarrow N^{r}\text{, } t\mapsto (t^{a_1},\ldots, t^{a_r})\text{, }\sigma\mapsto (\sigma,\ldots,\sigma)\]
defines a group homomorphism.
\end{lemma}
\begin{proof}
We have the relation $\sigma\cdot \sigma=-1$, and $-1$ is sent to $(-1,\ldots,-1)$ which equals $(\sigma^2,\ldots,\sigma^2)$. Also the map described above respects the relation  $\sigma\cdot t=t^{-1}\sigma$, so gives a well-defined homomorphism of group schemes.
\end{proof}

Let $s=\floor{\frac{n+1}{2}}$ and choose $a_1,\ldots,a_s$ odd positive, pairwise distinct integers. 
$N$ acts on $\P^n$ via $N\rightarrow N^{\floor{\frac{n+1}{2}}}\hookrightarrow \GL_{n+1}$ where the first map equals the map in the lemma above. This action extends to an action on $H_n$.

A point $x\in H_n$ is fixed by the action if and only if the corresponding curve $C_x$ is mapped to itself via the action of $N$ on $\P^n$. Every twisted cubic spans a unique 3-plane, so this 3-plane needs to be fixed by the action for the corresponding point to be a fixed point. Thus, the possible fixed 3-planes are those subspaces $\Pi_{ij}$, $1\le i<j\le s$ defined by  by $x_\ell=0$, $\ell\neq 2i-2, 2i-1, 2j-2, 2j-1$. We will concentrate on the case $\P^3=\Pi_{12}$ defined by  $x_4=\ldots=x_n=0$; the remaining cases will follow by applying the suitable shuffle permutation to the variables. 

\begin{lemma}\label{lemma:FixCurve}
The action of $N$ on $H_3$ has 6 fixed points $y_1,\ldots,y_6$ with defining ideal of the corresponding curve\\[5pt]
$y_1:=(x_0x_2,x_0x_1,x_1x_3)$,\ $y_2:=(x_0x_2,x_2x_3,x_1x_3)$,\
 $y_3:=(x_0x_3,x_0x_1,x_1x_2)$, \\[2pt]
 $y_4:=(x_0x_3,x_2x_3,x_1x_2)$,\
$y_5:=(x_0^2,x_0x_1,x_1^2)$,\
 $y_6:=(x_2^2,x_2x_3,x_3^2)$.
\end{lemma}
\begin{proof}
The points on $H_n$ that are fixed by the action of $t$ are exactly the points fixed by the classical torus action studied by Ellingsrud-Str{\o}mme. Thus the set of fixed points of the $N$-action is a subset of the set of fixed points in \cite[Proposition 3.8]{ES}. One checks that the six listed fixed points are exactly the ones which are fixed by the action of $\sigma$.
\end{proof}

\begin{corollary}
The action of $N$ on $H_n$ has $\binom{\floor{\frac{n+1}{2}}}{2}\cdot 6$ fixed points, each of which are in $H_n(k)$.
\end{corollary}
\begin{proof}
There are $\binom{\floor{\frac{n+1}{2}}}{2}$ fixed $3$-planes in $\P^n$ each containing $6$ twisted cubics which are fixed by the action. The fixed $3$-planes are each defined by subset of the variables $X_0,\ldots, X_n$, hence are all defined over $\Z$. Similarly, the fixed curves $C_x$ described in Lemma~\ref{lemma:FixCurve} are all defined over $\Z$, so all these $N$-fixed points are in $H_n(k)$.
\end{proof}

Assume $\sE_{m_1,\ldots,m_r; n}:=\sE_{m_1,n}\oplus\ldots\oplus \sE_{m_r,n}\rightarrow H_n$ is relatively oriented and that $\sum_{i=1}^r(3m_i+1)=4n$. We also assume that all the $m_j$ are odd, following Remark~\ref{rem:Cong2}(2). We will use the canonical relative orientation, discussed   in Proposition~\ref{prop:relOrientationObservations}

In this case the equivariant Euler class $e^N(\sE_{m_1,\ldots,m_r; n})$ of $\sE_{m_1,\ldots,m_r; n}$ lives in $H_N^{4n}(H_n,\mathcal{W}(\det^{-1}\sE_{m_1,\ldots,m_r; n}))$ which is isomorphic to
$H_N^{4n}(H_n,\mathcal{W}(\omega_{H_n/k}))$ via the relative orientation. Let $\pi_{H_n}:H_n\rightarrow \Spec k$ be the structure map. Then we can push $e^N(\sE_{m_1,\ldots,m_r; n})$ forward  and get $(\pi_{H_n})_*e^N(\sE_{m_1,\ldots,m_r; n})\in H^0(BN,\mathcal{W})\cong \operatorname{W}(k)[e]$.
Applying the quadratic version of Bott's residue formula \cite[Theorem 9.5]{LevineAB} we get
\begin{equation}
\label{eq: Bott formula for twisted cubics0}
(\pi_{H_n})_*e^N(\sE_{m_1,\ldots,m_r; n})=\sum_{\text{fixed $3$-planes}}\left(\sum_{i=1}^6\frac{e^N(\sE_{m_1,\ldots,m_r; n}(y_i))}{e^N( \sT_{H_n}(y_i))}\right)\in \operatorname{W}(k)[e, 1/Me]
\end{equation}
for a suitable integer $M>0$. 

Of course, we are really interested in the ``quadratic degree'' 
\[
(\pi_{H_n})_*e(\sE_{m_1,\ldots,m_r; n})\in W(k)
\]
of our non-equivariant class $e(\sE_{m_1,\ldots,m_r; n})\in H^{4n}(H_n,\mathcal{W}(\omega_{H_n/k}))$.  This is the image of $(\pi_{H_n})_*e^N(\sE_{m_1,\ldots,m_r; n})$ under the reduction map $W(k)[e]\to W(k)$ sending $e$ to zero.  Both $e^N(\sE_{m_1,\ldots,m_r; n}(y_i))$ and $e^N( \sT_{H_n}(y_i))$ are elements of $H^d(BN, \sW(\gamma^{d_{y_i}}))$ for some integer $d_{y_i}$ depending on $y_i$, and where $d=\dim_kH_n$. Thus their ratio is a well-defined element of  $W(k)[1/M]\subset W(k)[e. 1/Me]$, and \eqref{eq: Bott formula for twisted cubics0} gives us the identity
\begin{equation}
\label{eq: Bott formula for twisted cubics}
(\pi_{H_n})_*e(\sE_{m_1,\ldots,m_r; n})=\sum_{\text{fixed $3$-planes}}\left(\sum_{i=1}^6\frac{e^N(\sE_{m_1,\ldots,m_r; n}(y_i))}{e^N( \sT_{H_n}(y_i))}\right)\in \operatorname{W}(k)[1/Me]
\end{equation}

\section{Relatively oriented bases and Euler class computation}
In order to compute \eqref{eq: Bott formula for twisted cubics}, we need to compute the equivariant Euler classes $e^N(\sE_{m,n}(y))$ and $e^N(\sT_{H_n}(y))$ for each fixed point $y$. We use the following strategy. We write the $N$-representations $\sE_{m,n}(y)$ and $\sT_{H_n}(y)$ as a direct sum of irreducible $N$-representations. Theorem \ref{thm:EulerClasses} gives us the equivariant Euler class of each irreducible summand and by the Whitney sum formula, the equivariant Euler classes of $\sE_{m,n}(y)$ and $\sT_{H_n}(y)$ are the respective products of the equivariant Euler classes of their irreducible summands. 

In order to get the correct choice of relatively oriented bases for  $\sT_{H_n}(y)$ and $\sE_{m,n}(y)$, we first choose bases for  $\sE_n(y)$ and $\sF_n(y)$. We have seen above that $y$ is in $H_n^\cm$. We will then use the resolutions \eqref{eqn:TangentPresentation} and \eqref{eqn:ResOnHn}
 for $\sT_{\Phi}(y)$ and $\sE_{m,n}(y)$ and choose our bases for $\sT_{\Phi}(y)$ and $\sE_{m,n}(y)$ so that they are compatible with respect to the resolutions and the given choice of bases for  $\sE_n(y)$ and $\sF_n(y)$. We will then use a choice of basis for $\sT_{\Gr(4, n+1)}(\Phi(y))$ to give a basis for $\sT_{H_n}(y)$. 
 
In case our chosen bases for $\sT_{\Phi}(y)$ or $\sE_{m,n}(y)$ are not oriented in the sense of Remark~\ref{rem:OrientedBasis}, we will keep track of the factors needed to make the basis oriented in this sense (these will always be $\pm1$) and correct our computation by the product of all these factors. 

We have the Euler sequence on $\Gr(4, n+1)$
\[
0\to \End(E_4)\to E_4^\vee\otimes \sO^{n+1}_{\Gr(4, n+1)}\to \sT_{\Gr(4, n+1)}\to 0
\]
and the canonical identification $\sO_{\Pi_n(\Phi(y))}(1)=E_4^\vee(y)$. Using this, our choice of coordinates for the projective space  $\Pi_n(\Phi(y))$ gives an oriented basis for  $ \sT_{\Gr(4, n+1)}(\Phi(y))$, and the identity
\[
\det \sT_{\Gr(4, n+1)}(\Phi(y))=(\det E_4^\vee)^{\otimes n+1}=\sO_{\Gr(4, n+1)}(n+1).
\]

Note that $\det\sE_n(y)$ and $\det\sF_n(y)$ both appear to even powers in $\det\sT_{H_n}(y)$ and
since we will always take $m$ odd, $\det\sE_n(y)$ appears to an even power in $\det\sE_{m,n}$.
For $m\equiv1\mod4$, the same holds for $\det\sF_n(y)$, but if $m\equiv-1\mod4$, 
 $\det\sF_n(y)$ occurs to an odd power in  $\det\sE_{m,n}$. However, if we apply  our method to each of $\sE_{m_1,n},\ldots, \sE_{m_r, n}$ and  $\oplus_{j=1}^r\sE_{m_j, n}$ is relatively oriented, then the number of $m_i\equiv-1\mod 4$ is even, so in all relatively oriented cases, our choice of basis for $\sE_n(y_j)$ and $\sF_n(y_j)$ at each of the fixed points $y_j$ will not play a role. At the individual $y_j$, we can even make an arbitrary change of basis for $\sE_n(y_j)$ in our computations, and a change of basis with determinant one for $\sF_n(y_j)$, without affecting the global result. 
 
 Similarly, the relative orientation condition implies that the powers of $\sO_{\Gr(4, n+1)}(1)$ occurring in $\det\oplus_{j=1}^r\sE_{m_j, n}$ and in $\det\sT_{H_n}$ are either both odd or both even, so our choice of basis for $\sT_{\Gr(4, n+1)}(\Phi(y))$ as above will give us a  bases for $\sT_{H_n}(y)$  and $\oplus_{j=1}^r\sE_{m_j, n}(y)$ that are relatively oriented with respect to our choice of global relative orientation.

Let $E=\sE_n(y)$ and $F=\sF_n(y)$ be the fibers  at a fixed point $y$. As all our fixed points are in the Cohen-Macaulay locus, the restriction of the complex  \eqref{eqn:ResOnHn} gives an equivariant resolution
\begin{multline}
    \label{eq:resolution bundle}
    0\rightarrow F\otimes H^0(\mathcal{O}_{\P^3}(d-3))\xrightarrow{L(y)}  E\otimes H^0(\mathcal{O}_{\P^3}(d-2))\\\xrightarrow{\bigwedge^2L(y)^t} H^0(\mathcal{O}_{\P^3}(d))\xrightarrow{\pi_{m,n}} \sE_{m,n}(y)\rightarrow 0
\end{multline}
of $\sE_{m,n}(y)$ (see Remark~\ref{rem:GEquivarianceProperty}). Similarly, the restriction of the complex \eqref{eqn:TangentPresentation} gives an equivariant resolution
\begin{multline}
    \label{eq:resulution tangent space}
    0\rightarrow k\xrightarrow{(\id_F,-\id_E)} \End(F)\oplus \End(E)\\\xrightarrow{(L(y)_*, L(y)^*)}  \Hom(F,E)\otimes H^0(\mathcal{O}_{\P^3}(1))\xrightarrow{\pi} \sT_{H_3}(y)\rightarrow 0
\end{multline}
of $ \sT_{H_3}(y)$.  The matrix $L(y)$ and bases for $E$ and $F$ are not uniquely determined by $y$; we fix our choices of $L(y)$ and the corresponding bases $e_1, e_2, e_3$ of $E$ and $f_1, f_2$ of $F$ as follows
\[
L(y_1)=\begin{pmatrix}-x_1&0\\x_2&x_3\\0&-x_0\end{pmatrix},\
L(y_2)=\begin{pmatrix}-x_3&0\\x_0&x_1\\0&-x_2\end{pmatrix},\ 
L(y_3)=\begin{pmatrix}0&-x_1\\x_2&x_3\\-x_0&0\end{pmatrix},
\]
\[
L(y_4)=\begin{pmatrix}0&-x_3\\x_0&x_1\\-x_2&0\end{pmatrix},\
L(y_5)=\begin{pmatrix}-x_1&0\\x_0&x_1\\0&-x_0\end{pmatrix},\
L(y_6)=\begin{pmatrix}-x_3&0\\x_2&x_3\\0&-x_2\end{pmatrix},
\]
in other words, $L(y_1)(f_1)=-x_1e_1+x_2e_2$, etc.

We now proceed to implement the strategy outlined above.  We have our action of $N$ on $\P^n$ and $H_n$ arising from our choice of odd integers $a_1>a_2,\ldots>a_s>0$, $s=\floor{\frac{n+1}{2}}$. For the computation involving the six $N$-fixed curves lying in $\Pi_{ij}$, we use the oriented basis $x_{2i-2}, x_{2i-1}, x_{2j-2}, x_{2j-1}$ with weights $a_i, -a_i, a_j, -a_j$ for $H^0(\Pi_{ij}, \sO(1))$. We fix an integer $M>0$ and assume that $a_i>M\cdot a_{i+1}$ for $i=1,\ldots, s$; for application to the study of the sheaf $\oplus_{i=1}^r\sE_{m_i,n}$, we will take $M=\max_i m_i$. 

In what follows, we will look at the case of $\Pi_{12}=\P^3$; the general case follows by symmetry. We let $V=H^0(\P^3, \sO_{\P^3}(1))$ with standard basis $x_0, x_1, x_2, x_3$.

\begin{definition}
Let $X$ be a $k$-vector space with an $N$ action. Take $x\in X\setminus\{0\}$  such that the line $k\cdot x$ is stable under the $\G_m$-action. Then $t\cdot x=t^wx$ for some integer $w\in \Z$. We say that $x$ is a weight vector with weight  $\wt(x)=w$.
\end{definition}

We use the canonical order defined in Definition \ref{def:CanonicalOrientation} on the monomials in $\Sym^rV$.  We assume that $r$ is at most our chosen integer $M$.
Recall from Definition \ref{def:CanonicalOrientation} that we call a monomial $g=x_0^{m_1}x_1^{m_2}x_3^{m_2}x_3^{m_4}$ 
\begin{itemize}
\item \emph{positive} if $m_1>m_2$ or $m_1=m_2$ and $m_3>m_4$. Note that is is exactly the case when $\wt(g)>0$.
\item \emph{negative} if $m_1<m_2$ or $m_1=m_1$ and $m_3<m_4$. Note that is is exactly the case when $\wt(g)<0$,
\item \emph{neutral} $m_1=m_2$ and $m_3=m_4$ and note that is is exactly the case when $\wt(g)=0$.
\end{itemize}
Also recall that we set $g^*=x_0^{m_2}x_1^{m_1}x_2^{m_4}x_3^{m_3}$.

For our chosen $N$-action, the induced action on $V$ makes each monomial $m_\alpha$ in $x_0,\ldots,  x_3$ a weight vector, with $\sigma(m_\alpha)=\pm m_\alpha^*$. In particular, each $k$-vector subspace $k\cdot m_\alpha+k\cdot m_\alpha^*$ of $k[x_0,\ldots, x_3]$ is an irreducible $N$-representation, and has dimension two exactly when $m_\alpha^*\neq m_\alpha$.

The $N$-representation $E:=H^0(\P^3, \sI_{C_y}(2))$ is naturally a sub-representation of $\Sym^2V$, and thus has a basis of weight vectors consisting of degree two monomials in $x_0,\ldots, x_ 3$.

By our choice of bases for $E$ and $F$, and the matrices $L(y_j)$, we see that $f_1, f_2$, $e_1, e_2, e_3$ are weight vectors with   $\wt(f_1)>0$, $\wt(f_2)<0$, $\wt(e_1)>0$, $\wt(e_2)=0$ and $\wt(e_3)<0$.   We call the bases $(e_1, e_2, e_3)$ for $E$, $(f_1, f_2)$ for $F$ and $x_0,\ldots, x_3$ for $V$ {\em standard bases}. Similarly, we call  $x_4,\ldots, x_n$ the standard basis for $H^0(\P^n,. \sO_{\P^n}(1))/V$. We make an analogous definition of the matrices $L(y_j)$ and standard bases in the case of the other fixed subspaces $\Pi_{ij}$ by acting by the evident shuffle permutation of the variables. We extend the $*$-action on monomials to $E$ and $F$ by setting $e_1^*=e_3$, $e_3^*=e_1$, $e_2^*=e_2$, $f_1^*=f_2$ and $f_2^*=f_1$.

\begin{ex} For $y_1=(x_0x_2, x_0x_1, x_1x_3)$ we have $e_1\mapsto -x_0x_2$, $e_2\mapsto -x_0x_1$, $e_3\mapsto -x_1x_3$. The weights are
\[
\wt(e_1)=a_1+a_2, \wt(e_2)=0, \wt(e_3)=-a_1-a_2
\]
\[
\wt(f_1)=a_2, \wt(f_2)=-a_2
\]
and we have
\[
\sigma(e_1)=e_3, \sigma(e_2)=-e_2,  \sigma(e_3)=e_1,
\]
\[
 \sigma(f_1)=f_2, \sigma(f_2)=-f_1.
\]
Thus $E\cong \tilde{\sO}(a_1+a_2)\oplus\tilde{\sO}^-(0)$, $F\cong \tilde{\sO}^-(a_2)$ (recall that $a_2$ is odd and $a_1+a_2$ is even).
\end{ex}

We use the following convention for ordered bases for a tensor product. If we have constructed  bases $u_1,\ldots, u_r$, $w_1,\ldots, w_s$ for vector spaces $U$ and $W$, the tensor product $U\otimes W$  inherits the basis $\{u_i\otimes w_j\}_{i,j}$, where we use the order in which $u_i\otimes w_j$ precedes $u_{i'}\otimes w_{j'}$ if $j<j'$ or if $j=j'$ and $i<i'$:
\[
u_1\otimes w_1, u_2\otimes w_1,\ldots, u_r\otimes w_1, u_1\otimes w_2,\ldots, u_r\otimes w_2, \ldots, u_1\otimes w_s,\ldots, u_r\otimes w_s.
\]
  This fixes the canonical isomorphism
\[
\det(U\otimes W)\cong \det(U)^{\otimes \dim W}\otimes\det(W)^{\otimes\dim U}
\]

To construct our relatively oriented basis for $\sE_{m,n}(y)$, we start with our choices of standard bases for $E$,  $F$ and $V$, and use our canonical bases for $\Sym^jV$, $j=m-3, m-2, m$. This gives us  canonical bases the terms in the resolutions \eqref{eq:resolution bundle} and \eqref{eq:resulution tangent space} From this we can extract  $N$-oriented bases for $\sE_{m,n}(y)$ and $\sT_{\Phi}(y)$, and then at the end add in a similarly chosen basis for $\sT_{\Gr(4, n+1})(\Phi(y))$ to give our $N$-oriented basis for $\sT_{H_n}(y)$. From our discussion above, these bases will automatically be compatible with our chosen relative orientation for $\sE_{m,n}$, possibly up to a sign that we can keep track of. In what follows we use {\em oriented} to refer to oriented bases as $N$-representation, and {\em relatively oriented} to refer to our globally chosen relative orientation. 
 
As $\det F$, $\det E$ occur in $\det \sE_{m,n}(y)$ and in $\det\sT_{H_n}(y)$ to the same powers modulo two, the choice of standard bases for $F$ and $E$ will not matter in our construction of a relatively oriented bases for $\sE_{m,n}(y)$, as long as we use the same bases for the construction for $\sE_{m,n}(y)$ and for $\sT_{H_n}(y)$.

\subsection{Oriented bases for $\sE_{m,n}(y_j)$}
Note that there are no zero weight monomials in $\Sym^rV$ when $r$ is odd.
Consequently, the canonical ordering of $\Sym^rV$ defined in \ref{def:CanonicalOrientation} comes in pairs $(g,g^*)$ where $g$ is positive and $g^*$ is negative. Also note
that $\sigma$ sends $g$ to $\pm g^*$ and hence the vector spaces with basis $(g,g^*)$ are invariant under the $N$-action. Note further that changing the order of the pairs yields an even permutation of a given basis, and thus is orientation preserving.

The $E\otimes \Sym^{m-2}V$-term breaks up into $6$-dimensional invariant subspaces with ordered bases
\[(e_1g,e_2g,e_3g,e_1g^*,e_2g^*,e_3g^*)\]
where $g$ is a monomial of degree $m-2$ and positive weight.
We rearrange this basis to
\[(e_1g,e_3g^*,e_2g,e_2g^*,e_3g,e_1g^*).\]
Note that this base change has determinant $+1$ and splits up into $3$ pairs of weight vectors, dual with respect to the $\sigma$-action, namely $(e_1g,e_3g^*)$, $(e_2g,e_2g^*)$ and $(e_3g,e_1g^*)$ where the first and the second one are oriented, that is,  the positive weight vector comes first. It might be that $\wt(e_3g)<\wt(e_1g^*)$ in which case we want to swap the order and keep track of an additional sign.

The $F\otimes \Sym^{m-3}V$-term breaks up into $4$-dimensional (when $g$ is positive) respectively $2$-dimensional (when $h$ is neutral) subspaces with ordered basis \[(f_1g,f_2g,f_1g^*,f_2g^*) \]
respectively
\[(f_1h,f_2h).\]
The latter is of the desired form.
We rearrange the first basis (change of basis with determinant $=+1$) to
\[(f_1g,f_2g^*,f_2g,f_1g^*)\]
which splits up into 2 $2$-dimensional invariant subspaces, with bases $(f_1g, f_2g^*)$ and $(f_2g,f_1g^*)$, both dual with respect to the $\sigma$-action.  The basis $(f_1g, f_2g^*)$ is an oriented basis as $N$-representation, while we might need to switch the order of $(f_2g,f_1g^*)$ to yield an oriented basis.  

We compute the number of total swaps to get our desired basis for $m$ odd.
\begin{proposition}
\label{prop:number of swaps}
\begin{enumerate}
\item Assume $m\equiv 1\mod 4$. Then the number of monomials $g$ of degree $m-2$ for which $\wt(ge_3)<\wt(g^*e_1)$ plus the number of monomials $g$ of degree $m-3$ for which $\wt(f_2g)<\wt(f_1g^*)$ is even. So in this case we have an even number of swaps.
\item Assume $m\equiv-1\mod 4$. Then the number of swaps is even for the first, second and fifth fixed point, while it is odd for the others.
\end{enumerate}
\end{proposition} 

\begin{lemma}
\label{lemma: cases}
Let $m$ be an odd, positive integer.
\\[5pt]
(1) $\#\{g:\text{$g$ is a monomial of degree $m-2$ and $0<\wt(g)<a_1+a_2$}\}=\frac{(m-1)(m+1)}{4}$\\[2pt]
(2)  $\#\{g:\text{$g$ is a monomial of degree $m-2$ and $0<\wt(g)<a_1-a_2$}\}=\frac{(m-1)^2}{4}$\\[2pt]
(3) $\#\{g: \text{$g$ is a monomial of degree $m-2$ and $0<\wt(g)<2a_2$}\}=\frac{m-1}{2}$\\[2pt]
(4) $\#\{g: \text{$g$ is a monomial of degree $m-2$ and $0<\wt(g)<2a_1$}\}=\frac{(m-1)^2}{2}$\\[2pt]
(5) $\#\{g: \text{$g$ is a monomial of degree $m-3$ and $0<\wt(g)<a_2$}\}=0$\\[2pt]
(6) $\#\{g: \text{$g$ is a monomial of degree $m-3$ and $0<\wt(g)<a_1$}\}=\frac{(m-1)(m-3)}{4}$
\end{lemma}

\begin{proof}
(1) Either $\wt(g)=j\cdot a_2$ for some $j\in \mathbb{Z}_{>0}$ or $\wt(g)=a_1-j\cdot a_2$ for $j\in \mathbb{Z}_{\ge 0}$. In the first case $g=x_0^{i}x_1^{i}x_2^px_3^q$ for some $i\ge 0$ and $p,q\ge 0$ such that $p>q$ and $p+q=m-2-2i$. For a fixed $i$ there are $\frac{m-1-2i}{2}$ such monomials, namely
\[
    x_0^ix_1^ix_2^{\frac{m-2-2i+1}{2}}x_3^{\frac{m-2-2i-1}{2}}, x_0^ix_1^ix_2^{\frac{m-2-2i+3}{2}}x_3^{\frac{m-2-2i-3}{2}},
    \ldots, x_0^ix_1^ix_2^{m-2-2i}.
\]
    So we get 
    \[\sum_{i=0}^{\frac{m-3}{2}}\frac{m-1-2i}{2}=1+2+\ldots +\frac{m-1}{2}=\binom{\frac{m+1}{2}}{2}=\frac{(m+1)(m-1)}{8}\]
    in total.
    
    In the second case we have $g=x_0^{i+1}x_1^ix_2^px_3^q$ for some $i\ge 0$, $p,q\ge 0$ with $p+q=m-2-2i-1$ and $q\ge p$. For a fixed $i$ there are $\frac{m-1-2i}{2}$ such $g$'s, namely
    \[x_0^{i+1}x_1x_2^{\frac{m-2-2i-1}{2}}x_3^{\frac{m-2-2i-1}{2}},x_0^{i+1}x_1x_2^{\frac{m-2-2i-3}{2}}x_3^{\frac{m-2-2i+1}{2}},\ldots,x_0^{i+1}x_1^ix_3^{m-2-2i-1}\]
    and in total we get 
    \[\sum_{i=0}^{\frac{m-3}{2}}\frac{m-1-2i}{2}=1+2+\ldots +\frac{m-1}{2}=\binom{\frac{m+1}{2}}{2}=\frac{(m+1)(m-1)}{8}.\]
    
    So there are $\frac{(m+1)(m-1)}{4}$ monomials $g$ of degree $m-2$ and $0<\wt(g)<a_1+a_2$.\\[5pt]
  (2) Either $\wt(g)=j\cdot a_2$ for $j\in \mathbb{Z}_{>0}$ or $\wt(g)=a_1-j\cdot a_2$ for some $j\in \mathbb{Z}_{>1}$. We already know that there are $\frac{(m-1)(m+1)}{8}$ monomials with weight $j\cdot a_2$. 
    The monomials of weight $a_1-j\cdot a_2$ with $j>1$ are of the form $x_0^{i+1}x_1^ix_2^px_3^q$ for some $i\ge 0$, $p,q\ge 0$ with $p+q=m-2-2i-1$ and $q> p+1$. For a fixed $i$ there are $\frac{m-3-2i}{2}$ such monomials, namely
    \[x_0^{i+1}x_1x_2^{\frac{m-2-2i-3}{2}}x_3^{\frac{m-2-2i+1}{2}},x_0^{i+1}x_1x_2^{\frac{m-2-2i-5}{2}}x_3^{\frac{m-2-2i+3}{2}},\ldots, x_0^{i+1}x_1^ix_3^{m-2-2i-1}.\]
    This makes in total $\binom{\frac{m-1}{2}}{2}=\frac{(m-1)(m-3)}{8}$.
    Adding the $2$ numbers we get
    \[\frac{(m-1)(m+1)}{8}+\frac{(m-1)(m-3)}{8}=\frac{m^2-1+m^2-4m+3}{8}=\frac{(m-1)^2}{4}\]
\ \\[2pt]
(3) We are looking for monomials of degree $m-2$ and weight $a_2$. These are of the form $x_0^ix_1^ix_2^{\frac{m-2-2i+1}{2}}x_3^{\frac{m-2-2i-1}{2}}$ and one checks that there are $\frac{m-1}{2}$ of these.\\[2pt]
(4) $g$ can either have weight $j\cdot a_2$ for some $j\in \mathbb{Z}_{>0}$ or $a_1+j\cdot a_2$ for $j\in \mathbb{Z}$ or $2a_1-j\cdot a_2$ for $j>0$. In the first case $g$ is of the form $g=x_0^ix_1^ix_2^px_3^q$ for $i,p,q\ge 0$ with $p>q$ and $m-2-2i=p+q$. For a fixed $i$ these are the following
    \[x_0^ix_1^ix_2^{\frac{m-2-2i+1}{2}}x_3^{\frac{m-2-2i-1}{2}},x_0^ix_1^ix_2^{\frac{m-2-2i+3}{2}}x_3^{\frac{m-2-2i-3}{2}},\ldots, x_0^ix_1^ix_2^{m-2-2i}\]
    One counts that the number of such monomials is $\frac{m-1-2i}{2}$.
    In total we get $\binom{\frac{d+1}{2}}{2}=\frac{(m-1)(m+1)}{8}$.
    
    In the second case the monomials are of the form $g=x_0^{i+1}x_1^ix_2^px_3^q$ for $i,p,q\ge 0$ with $m-2-2i=p+q$. For a fixed $i$ there are $m-2-(2i+1)+1$ such monomials and in total we get
\begin{multline*}
\sum_{i=0}^{\frac{m-3}{2}}m-2-2i=\frac{m-1}{2}\cdot (m-2)-2\cdot (1+2+\ldots+\frac{m-3}{2})\\=\frac{m-1}{2}\cdot(m-2)-2\cdot\binom{\frac{m-1}{2}}{2}=\frac{(m-1)^2}{4}
\end{multline*}    
    For the third weight class $2a_1-j\cdot a_2$ the monomials look as follows $x_0^{i+2}x_1^ix_2^px_3^q$ for $i,p,q\ge 0$ with $p+q=m-2-(2i+2)$ and $q>p$. For a fixed $i$ these are the following ones
    \[x_0^{i+2}x_1^ix_2^{\frac{m-2-(2i+2)-1}{2}}x_3^{\frac{m-2-(2i+2)+1}{2}},\ldots, x_0^{i+2}x_1^ix_3^{m-2-(2i+2)}\]
    which are $\frac{m-1-(2i+2)}{2}$ many. Summing up over all possible $i$, we get $\binom{\frac{m-1}{2}}{2}=\frac{(m-1)(m-3)}{8}$ many.
    
    So in total we get 
    \[\frac{(m-1)(m+1)}{8}+\frac{(m-1)^2}{4}+\frac{(m-1)(m-3)}{8}=\frac{(m-1)^2}{2}\]
\ \\
(5) No monomial has weight in this range. \\[2pt]
(6) These are all the monomials of degree $m-3$ of weight either $j\cdot a_2$ for some $j>0$ or $a_1-j\cdot a_2$ for $j>0$. In the first case these are of the form $x_0^ix_1^ix_2^px_3^q$ for $i,q,p\ge 0$ with $p+q=m-3-2i$ and $p>q$. For a fixed $i$ these are the following 
    \[x_0^ix_1^ix_2^{\frac{m-3-2i+2}{2}}x_3^{\frac{m-3-2i-2}{2}},x_0^ix_1^ix_2^{\frac{m-3-2i+4}{2}}x_3^{\frac{m-3-2i-4}{2}},\ldots,x_0^ix_1^ix_2^{m-3-2i}\]
    which makes $\frac{m-3-2i}{2}$ in total. So we get
    $\binom{\frac{m-1}{2}}{2}=\frac{(m-1)(m-3)}{8}$
    monomials of positive weight a multiple of $a_2$.
    
    The monomials of weight $a_1-j\cdot a_2$ are of the form $x_0^{i+1}x_1^ix_2^px_3^q$ for $p,q,i\ge 0$ with $q>p$ and $m-3-(2i+1)=p+q$. For a fixed $i$ we get the following monomials
    \[x_0^{i+1}x_1^ix_2^{\frac{m-4-2i-1}{2}}x_3^{\frac{m-4-2i+1}{2}},x_0^{i+1}x_1^ix_2^{\frac{m-4-2i-3}{2}}x_3^{\frac{m-4-2i+3}{2}},\ldots ,x_0^{i+1}x_1^ix_3^{m-4-2i}\]
    which makes $\frac{m-3-2i}{2}$ in total.
    Summing up over all possible $i$ we get $\binom{\frac{m-1}{2}}{2}=\frac{(m-1)(m-3)}{8}$.
    
    In total we get $\frac{(m-1)(m-3)}{4}$.
\end{proof}
  
\begin{proof}[Proof of Proposition \ref{prop:number of swaps}]
We need to swap $e_3g$ and $e_1g^*$ whenever
\[\wt(ge_3)<0\Leftrightarrow 0<\wt(g)<\wt(e_1)\]
and we need to swap $f_2g$ and $f_1g^*$ whenever
\[\wt(f_2g)<0\Leftrightarrow0<\wt(g)<\wt(f_1)\]
\ \\
(1) First fixed point: $\wt(e_1)=a_1+a_2$ and $\wt(f_1)=a_2$. By Lemma \ref{lemma: cases} we need to swap $\frac{(m-1)(m+1)}{4}+0=\frac{(m-1)(m+1)}{4}$ times which is always even since $m$ is odd.\\[2pt]
(2) Second fixed point: $\wt(e_1)=a_1+a_2$ and $\wt(f_1)=a_1$. By Lemma \ref{lemma: cases} we need to swap $\frac{(m-1)(m+1)}{4}+\frac{(m-1)(m-3)}{4}$ times which is always even (both summands are even) since $m$ is odd.\\[2pt]
(3)  Third fixed point: $\wt(e_1)=a_1-a_2$ and $\wt(f_1)=a_2$. By Lemma \ref{lemma: cases} we need to swap $\frac{(m-1)^2}{4}+0=\frac{(m-1)^2}{4}$ times which is always even for $m\equiv 1\mod 4$ and odd for $m\equiv -1\mod 4$.
\\[2pt]
(4) Fourth fixed point: $\wt(e_1)=a_1-a_2$ and $\wt(f_1)=a_1$. By Lemma \ref{lemma: cases} we need to swap $\frac{(m-1)^2}{4}+\frac{(m-1)(m-3)}{4}$ times which is an even number for $m\equiv 1\mod 4$ and odd for $m\equiv -1 \mod 4$.\\[2pt]
Fifth fixed point: $\wt(e_1)=2a_1$ and $\wt(f_1)=a_1$. By Lemma \ref{lemma: cases} we need to swap $\frac{(m-1)^2}{2}+\frac{(m-1)(m-3)}{4}$ times which is even for $m$ odd.
\\[2pt]
(6) Sixth fixed point: $\wt(e_1)=2a_2$ and $\wt(f_1)=a_2$. We get an even number of swaps for $m\equiv 1\mod 4$ namely $\frac{m-1}{2}$. We get an odd number for $m\equiv-1\mod 4$.
\end{proof}

We note that for each $j$, the map $E\to \Sym^2(V)$ sends each $e_i$ to a (degree 2) monomial. Thus, $\sE_{m,n}(y_j)$ has as basis the images of the monomials in $\Sym^mV$ that are not in the image of $E\otimes \Sym^{m-2}V\to \Sym^mV$. Using the canonical bases for $\Sym^{m-3}V$, $\Sym^{m-2}V$ and $\Sym^mV$, together with our choice of standard bases for $E$ and $F$ and the resulting bases for $F\otimes \Sym^{m-3}V$ and $E\otimes \Sym^{m-2}V$, the resolution \eqref{eq:resolution bundle} gives us the canonical isomorphism
\begin{equation}\label{eqn:CanDetIso}
\delta(y_j): \det\sE_{m,n}(y_j)\xrightarrow{\sim} \det\Sym^mV\otimes(\det E\otimes \Sym^{m-2}V)^{-1}\otimes \det F\otimes \Sym^{m-3}V\cong k
\end{equation}
\begin{definition} Let $m\ge3$ be an odd integer. We call a basis $g_1, g_2,\ldots, g_D$ of $\sE_{m,n}(y_j)$  {\em canonically oriented} if $g_1\wedge\ldots\wedge g_D$ maps to $1\in k$ under the isomorphism \eqref{eqn:CanDetIso}. We call such a basis 
{\em anti-canonically oriented} if  if $g_1\wedge\ldots\wedge g_D$ maps to $-1\in k$ under this isomorphism.

We extend these notions to the vector spaces $\sE_{m_1,\ldots, m_r,n}(y_j)$, with the  $m_j\ge3$ odd integers, in the evident way.
\end{definition}

\begin{proposition}\label{prop:SignsInOrientedBases} Let $m\ge3$ be an odd integer, fix $j\in\{1,\ldots, 6\}$, let $\{M_1,\ldots, M_d\}$ be the set of positive monomials $M_i$ in $\Sym^mV$ that are not in the image of $E\otimes \Sym^{m-2}V\to \Sym^mV$ (in any order), and let $m_i, m_i^*$ be the respective images of $M_i, M_i^*$ in  $\sE_{m,n}(y_j)$. Then the basis $m_1, m_1^*, m_2, m_2^*, \ldots, m_d, m_d^*$ of $\sE_{m,n}(y_j)$ is  canonically oriented in case $j=1,2,5$, or $j=3,4,6$ and $m\equiv 1\mod 4$, and is 
anti-canonically oriented  in case $j=3,4,6$ and $m\equiv -1\mod 4$.
\end{proposition}

\begin{proof} We first give the proof in case $j=1$. Let $W\subset E\otimes\Sym^{m-2}V$ be the subspace  $e_2\otimes \Sym^{m-2}V$. We have the maps $L(y_1):F\otimes\Sym^{m-3}V\to E\otimes\Sym^{m-2}V$, $\bigwedge^2L(y)^t:E\otimes\Sym^{m-2}V\to\Sym^mV$ and $\pi:\Sym^mV\to \sE_{m,n}(y_j)$ from the complex \eqref{eq:resolution bundle}.

We replace $L(y_1)$ with $\alpha=-L(y_1)$ and replace $\bigwedge^2L(y)^t$ with $\beta=-\bigwedge^2L(y)^t$ in our complex;  since this change is equivalent to the change of basis of $E$, $e_i\mapsto -e_i$, and since $\Sym^{m-2}V$ has even dimension, this change arises from  a change of basis of $E\otimes\Sym^{m-2}V$ of determinant $+1$. 

Since the restriction of the map $\beta:E\otimes\Sym^{m-2}V\to \Sym^mV$ to $W$ is just the injective map $\times x_0x_1:\Sym^{m-2}V\to \Sym^mV$, $W\cap\alpha(F\otimes \Sym^{m-3}V)=\{0\}$. 

Note that the map $\alpha:F\otimes \Sym^{m-3}V\to E\otimes\Sym^{m-2}V$ sends $f_1\otimes m$ to $e_1\otimes x_1m-e_2\otimes x_2m$ and sends $f_2\otimes m$ to $e_3\otimes x_0m-e_2\otimes x_3m$. Form a new basis for $E\otimes\Sym^{m-2}V$ by starting with the canonical basis, and replacing for each monomial $m\in  \Sym^{m-3}V$ the basis element $e_1\otimes x_1m$
with $e_1\otimes x_1m-e_2\otimes x_2m$ and the basis element $e_3\otimes x_0m$ with 
$e_3\otimes x_0m-e_2\otimes x_3m$. Then since $W\cap  \im(F\otimes \Sym^{m-3}V)=\{0\}$ and
$W\cap e_1\otimes\Sym^{m-2}V\oplus e_3\otimes \Sym^{m-2}V=\{0\}$ this is indeed a basis for 
$E\otimes\Sym^{m-2}V$ and the  change of basis matrix has determinant one.

Note also that 
\[
(e_1\otimes x_1m)^*=e_3\otimes x_0m^*, (e_2\otimes x_2m)^*=e_2\otimes x_3m^*
\]
and 
\[
\wt(e_1\otimes x_1m)=\wt(e_2\otimes x_2m),  \wt(e_3\otimes x_0m)=\wt(e_2\otimes x_3m),
\]
so under this change of basis, we send the $\sigma$-dual pair of weight vectors
 $(e_1\otimes x_1m, e_3\otimes x_0m^*)$ to the $\sigma$-dual pair of weight vectors
$(\alpha(f_1\otimes m), \alpha(f_2\otimes m^*))$.  

We pass from the canonical basis of $E\otimes\Sym^{m-2}V$ to the basis of oriented pairs 
$(g, g^*)$ by applying a permutation. We apply the same permutation after substituting $\alpha(f_1\otimes m)$ for $e_1\otimes x_1m$ and $\alpha(f_2\otimes m)$ for $e_3\otimes x_0m$.   We order the oriented pairs $(g,g^*)$ by placing those of the form $(\alpha(f_1\otimes m), \alpha(f_2\otimes m^*))$ for $(f_1\otimes m, f_2\otimes m^*)$ an oriented pair in $F\otimes \Sym^{m-3}V$ first, followed by the oriented pairs of the form $(e_1\otimes m, e_3\otimes m^*)$, $(e_3\otimes m, e_1\otimes m^*)$, or $(e_2\otimes m, e_2\otimes m^*)$. The first type gives a basis for $\alpha(F\otimes \Sym^{m-3}V)$ and the image in $\Sym^mV$ of the second type gives a basis of oriented pairs of monomials $(m, m^*)$ for $\beta(E\otimes\Sym^{m-2}V)\subset \Sym^mV$. We then complete the basis of
$\Sym^mV$ by using the  oriented pairs of monomials $(m, m^*)$ with  $m$ not in 
$\beta(E\otimes\Sym^{m-2}V)$. 

Altogether, this gives us a  basis $B_F=(b_1^F,\ldots, b_{d_F}^F)$ for $F\otimes\Sym^{m-3}V$, a  basis $B_E=(b_1^E:=\alpha(b_1^F),\ldots, b_{d_F}^E:=\alpha(b_{d_F}^F), b_{1}^E,\ldots, b_{d_E}^E)$ for $E\otimes\Sym^{m-2}V$, a  basis $B_V:=(\beta(b_1^E),\ldots, \beta(b_{d_E}^E), b_1^V,\ldots, b_{d_V}^V)$ for $\Sym^mV$, and the basis $B_{\sE} :=(\pi(b_1^V),\ldots, \pi(b_{d_V}^V))$ for $\sE_{m,n}(y_1)$. In particular, the basis elements $b_\ell^V$ are each an oriented pair of monomials $(M, M^*)$, and so 
$B_{\sE} $ is a basis for $\sE_{m,n}(y_1)$ of images of oriented pairs of monomials. 

The bases $B_F$, $B_E$,  $B_V$ and $B_\sE$ define generators for $\det F\otimes\Sym^{m-3}V$, $\det E\otimes \Sym^{m-2}V$, $\det\Sym^mV$ and $\det\sE_{m,n}(y_1)$, respectively, giving isomorphisms
\[
\phi_1(y_1):\det F\otimes\Sym^{m-3}V\otimes (\det E\otimes \Sym^{m-2}V)^{-1}\otimes \det\Sym^mV\xrightarrow{\sim} k
\]
\[
\phi_2(y_1):\det\sE_{m,n}(y_1)\xrightarrow{\sim} k
\]
and the composition $\phi_2(y_1)^{-1}\phi_1(y_1)$ is exactly the isomorphism $\delta(y_1)$ determined by the exact sequence \eqref{eq:resulution tangent space}.  By Proposition~\ref{prop:number of swaps},  the  isomorphism  $\phi_1(y_1)$ is the same as the one determined by the canonical bases,  so the basis $B_\sE$ is canonically oriented. This property is preserved by any reordering of oriented pairs, so we have proven the Proposition for $y_1$.

The proof in case $j=2,5$, or $j=3,4,6$ and $m\equiv1\mod 4$ is essentially the same. In case $j=3,4,6$ and $m\equiv-1\mod4$, Proposition~\ref{prop:number of swaps} says that the  isomorphism 
\[
\phi_1(y_j):\det F\otimes\Sym^{m-3}V\otimes (\det E\otimes \Sym^{m-2}V)^{-1}\otimes \det\Sym^mV\xrightarrow{\sim} k
\]
 constructed as for $\phi_1(y_1)$ will differ from the one using the canonical bases by a factor of $-1$, so the basis we construct will be anti-canonically oriented.
\end{proof}

\begin{remark}\label{rem:OrientedBases} We consider the direct sum $\sE_{m_1,\ldots, m_r;n}:=\oplus_{i=1}^r\sE_{m_i,n}$, with all $m_i$ odd (see Remark~\ref{rem:Cong2}). The case of interest for us is when $\sum_i3m_i+1=4n$ and $\sE_{m_1,\ldots, m_r;n}$ is relatively oriented, which by Proposition~\ref{prop:relOrientationObservations} implies that the number of $i$ such that $m_i\equiv-1\mod4$ is even. By Proposition~\ref{prop:SignsInOrientedBases}, this says that a basis of $\sE_{m_1,\ldots, m_r;n}$ consisting of oriented pairs of monomials $(g,g^*)$ for each of the $\sE_{m_i,n}$  is canonically oriented.
\end{remark}

\subsection{The equivariant Euler class}

Define $\sigma_m(a_1,a_2)(y_j)\in \{\pm1\}$, $j=1,\ldots, 6$, by
 \[
\sigma_m(a_1,a_2)(y_j)=\begin{cases}
\epsilon(a_1a_2)&m\equiv1\mod 4\\
-\epsilon(a_1)&m\equiv-1\mod 4,\ j=2,4,5\\
-\epsilon(a_2)&m\equiv-1\mod 4,\ j=1,3,6
\end{cases}
\]
 We call $\sigma_m(a_1,a_2)(y_j)$ the {\em sign for $e^N(\sE_{m,n})$ at $y_j$}.

As in \S\ref{sec:Det}, \eqref{eqn:BundleNotation}, we set $\sE_{m_1,\ldots,m_r;n}:=\oplus_{i=1}^r\sE_{m_i,n}$. Thus
\[
e^N(\sE_{m_1,\ldots,m_r;n}(y_j))=\prod_{i=1}^re^N(\sE_{m_i,n}(y_j))
\]
and we call
\[
\sigma_{m_1,\ldots, m_r}(a_1, a_2,)(y_j):=\prod_{i=1}^r\sigma_{m_i}(a_1,a_2)(y_j), 
\]
 the {\em sign for $\sE_{m_1,\ldots,m_r;n}$ at $y_j$}.
 
 This terminology for $r=1$ is justified by the following proposition; the justification in general  
 follows from the identity
 \[
 e^N(\sE_{m_1,\ldots,m_r;n}(y_j))=\prod_{i=1}^re^N(\sE_{m_i,n}(y_j)).
 \]
 
\begin{proposition}
\label{Prop:EulerClassBundle}  The Euler classes $e^N(\sE_{m,n}(y_j))$, $j=1,\ldots,6$, are as follows:
\scalebox{.95}{
\vbox{
\begin{align*}
&e^N(\sE_{m,n}(y_1))=\sigma_m(a_1,a_2)(y_1) \cdot \left[\prod_{i=0}^{m-1}((m-i)a_1-ia_2)\right]  \cdot m!!\cdot  a_2^{\frac{m+1}{2}} \cdot e^{\frac{3m+1}{2}}\\
&e^N(\sE_{m,n}(y_2))=\sigma_m(a_1,a_2)(y_2) \cdot \left[ \prod_{i=1}^{m-1}((m-i)a_1-ia_2)\right]  \cdot ma_2\cdot m!!\cdot a_1^{\frac{m+1}{2}}\cdot e^{\frac{3m+1}{2}}\\
 &e^N(\sE_{m,n}(y_3))=\sigma_m(a_1,a_2)(y_3) \cdot \left[ \prod_{i=0}^{m-1}((m-i)a_1+ia_2)\right]  \cdot m!!\cdot a_2^{\frac{m+1}{2}}\cdot e^{\frac{3m+1}{2}}\\
 &e^N(\sE_{m,n}(y_4))=\sigma_m(a_1,a_2)(y_4) \cdot  \left[\prod_{i=1}^{m-1}((m-i)a_1+ia_2)\right]  \cdot ma_2\cdot m!!\cdot a_1^{\frac{m+1}{2}}\cdot e^{\frac{3m+1}{2}}\\
 &e^N(\sE_{m,n}(y_5))=\sigma_m(a_1,a_2)(y_5) \cdot \left[\prod_{i=0}^{m-1}(a_1+(m-1-2i)a_2)\right]  \cdot m!!\cdot a_2^{\frac{m+1}{2}}\cdot e^{\frac{3m+1}{2}}\\
 &e^N(\sE_{m,n}(y_6))=\sigma_m(a_1,a_2)(y_6)\cdot\left[\prod_{i=0}^{\frac{m-3}{2}} ((m-1-2i)^2a_1^2-a_2^2)\right] \cdot a_2\cdot m!!\cdot a_1^{\frac{m+1}{2}}\cdot e^{\frac{3m+1}{2}}
 \end{align*} }
 }
\end{proposition}

Recall that  $\dim H_n=4n$ and $\rnk(\sE_{m_1,\ldots,m_r;n})=\sum_{i=1}^r(3m_i+1)$. 

\begin{proposition}\label{prop:Signs}
Assume  $\rnk(\sE_{m_1,\ldots,m_r;n}) =\dim H_n$ and that all the $m_i$ are odd.
Then the sign of $e^N(\sE_{m_1,\ldots,m_r;n})$ at $y_j$  is given by
\[\sigma_{m_1,\ldots, m_r}(a_1, a_2,)(y_j)=\begin{cases}\epsilon(a_1a_2)&\text{for $r$ odd,}\\
1&\text{for $r$ even.}\end{cases}\]
\end{proposition}
\begin{proof}
For fixed $a_1, a_2$, and $y_j$,   $\sigma_m(a_1,a_2)(y_j)$ depends only the residue of $m$ mod 4.  By Proposition~\ref{prop:relOrientationObservations} 
the number of $m_i$ with $m_i\equiv-1\mod 4$ is even. Thus 
\[
\sigma_{m_1,\ldots, m_r}(a_1, a_2)(y_j)=\prod_{m_i\equiv 1\mod 4}\sigma_{m_i}(a_1,a_2)(y_j).
\]
If $r$ is even then the number of the $m_i$  with $m_i\equiv1\mod 4$ is even,  so in this case
$\sigma_{m_1,\ldots, m_r}(a_1, a_2)(y_j)=1$.   If $r$ is odd, we have an odd number of $m_i$'s with $m_i\equiv 1\mod 4$ and for such $m_i$ we have $\sigma_{m_i}(a_1,a_2)(y_j)=\epsilon(a_1a_2)$. Thus  if $r$ is odd,  $\sigma_{m_1,\ldots, m_r}(a_1, a_2)(y_j)=\epsilon(a_1a_2)$ 
\end{proof}

In fact, the sign
$\sigma_{m_1,\ldots, m_r}(a_1, a_2)(y_j)$ is the same as the corresponding sign for the equivariant Euler class of $ \sT_{H_n}(y_j)$ (see Proposition~\ref{prop:Signs} and  Proposition~\ref{prop:EulerClassGrassmannian}) So all signs depending on the weights $a_1$ and $a_2$ cancel in the formula for $(\pi_{H_n})_*e(\sE_{m_1,\ldots,m_r;n})$ in \eqref{eq: Bott formula for twisted cubics}, as they should.

\begin{proof}[Proof of Proposition \ref{Prop:EulerClassBundle}]
For the fixed point $y_1:=(x_0x_2,x_0x_1,x_1x_3)$, the pairs of the basis of $\sE_{m,n}(y_1)$ are of the form
\begin{itemize}
    \item $(x_0^{m-i}x_3^i,x_1^{m-i}x_2^i)$, $i=0,\ldots,m-1$
    \item or $(x_2^{m-i}x_3^i,x_2^ix_3^{m-i})$, $i=0,\ldots,\frac{m-1}{2}$.
\end{itemize}
In the first case this is the representation $\rho^{(-1)^i}_{(m-i)a_1-ia_2}$ which has equivariant Euler class $$(-1)^i\epsilon((m-i)a_1-ia_2)\cdot ((m-i)a_1-ia_2)e.$$
In the second case this is the representation $\rho^{(-1)^i}_{(m-2i)a_2}$ which has Euler class $$(-1)^i\epsilon(m-2i)\epsilon(a_2)\cdot(m-2i)a_2e.$$

For the sign, if $m\equiv 1\mod 4$ we have $\epsilon(m-2i)=-1$ iff $i$ is odd. It follows that $(-1)^i\epsilon(m-2i)\epsilon(a_2)=\epsilon(a_2)$ and we get the following sign
\begin{align*}
  \prod_{i=0}^{m-1}(-1)^i\epsilon&((m-i)a_1-ia_2)\cdot \prod_{i=0}^{\frac{m-1}{2}}\epsilon(a_2)\\
    =&\epsilon(a_2)\cdot \epsilon(a_1)\cdot (-1)^{\frac{m-1}{2}}\cdot \prod_{i=1}^{m-1}\epsilon((m-i)a_1-ia_2)\\
    =&\epsilon(a_1)\epsilon(a_2)\cdot \left(\epsilon\left(\prod_{i=1}^4((m-i)a_1-ia_2)\right)\right)^{\frac{m-1}{4}}\\
    =&\epsilon(a_1)\epsilon(a_2)\cdot \left(\epsilon((4a_1-a_2)(3a_1-2a_2)(2a_1-3a_2)(a_1-4a_2))\right)^{\frac{m-1}{4}}\\
    =&\epsilon(a_1)\epsilon(a_2)\cdot\epsilon(-a_1a_2(6a_1^2-13a_1a_2+6a_2^2)^{\frac{m-1}{4}}\\
    =&\epsilon(a_1)\epsilon(a_2)\cdot \epsilon(-a_1a_2(6\cdot 1-a_1a_2+6\cdot 1))^{\frac{m-1}{4}}\\
    =& \epsilon(a_1)\epsilon(a_2)\cdot \epsilon((-a_1a_2)(-a_1a_2))^{\frac{m-1}{4}}\\
    =&\epsilon(a_1)\epsilon(a_2)\\
    =& \sigma_m(a_1,a_2)(y_1).
\end{align*}
When $m\equiv-1\mod 4$, $\epsilon(m-2i)=-1$ iff $i$ is odd.
Hence, $(-1)^i\epsilon(m-2i)\epsilon(a_2)=-\epsilon(a_2)$
In total we get the  sign 
\begin{align*}
 \prod_{i=0}^{m-1}(-1)^i\epsilon&((m-i)a_1-ia_2)\cdot \prod_{i=0}^{\frac{m-1}{2}}(-1)\epsilon(a_2)\\
    =&\epsilon(a_1)\cdot (-1)^{\frac{m-1}{2}}\cdot\prod_{i=1}^{m-1}\epsilon((m-i)a_1-ia_2)\\
=&-\epsilon(a_1)\cdot\left(\epsilon\left(\prod_{i=1}^4((m-i)a_1-ia_2)\right)\right)^{\frac{m-3}{4}}\\&\cdot (-1)\epsilon((m-(m-2))a_1-(m-2)a_2)\epsilon((m-(m-1))a_1-(m-1)a_2)\\
=&\epsilon(a_1)\cdot \epsilon(2a_1-a_2)\epsilon(a_1-2a_2)\\
=&-\epsilon(2a_1-a_2)=-\epsilon(a_2)\\
=& \sigma_m(a_1,a_2)(y_1).
\end{align*}

For the second fixed point $y_2:=(x_0x_2,x_2x_3,x_1x_3)$, the pairs of the basis of $\sE_{m,n}(y_2)$ are of the form
\begin{itemize}
    \item $(x_0^{m-i}x_3^i,x_1^{m-i}x_2^i)$, $i=1,\ldots,m-1$
    \item or $(x_2^m,x_3^m)$
    \item or $(x_0^{m-i}x_1^i,x_0^ix_1^{m-i})$, $i=0,\ldots,\frac{m-1}{2}$.
\end{itemize}
In the first case this is the representation $\rho^{(-1)^i}_{(m-i)a_1-ia_2}$ which has Euler class $$(-1)^i\epsilon((m-i)a_1-ia_2)\cdot ((m-i)a_1-ia_2)e.$$
In the second case this is $\rho_{ma_2}$ which has Euler class 
$$\epsilon(a_2)\cdot ma_2e$$
and in the last case this is the representation $\rho^{(-1)^i}_{(m-2i)a_1}$ which has Euler class $$(-1)^i\epsilon(m-2i)\epsilon(a_1)\cdot(m-2i)a_1e.$$

For the sign calculation, when $m\equiv1\mod 4$ we get
\begin{align*}
\prod_{i=1}^{m-1}(-1)^i\epsilon&((m-i)a_1-ia_2) \cdot \epsilon(ma_2)\cdot \prod_{i=0}^{\frac{m-1}{2}}(-1)^i\epsilon((m-2i)a_1)\\
    &= \epsilon(a_1)\epsilon(a_2)\cdot \prod_{i=1}^{m-1}(-1)^i\epsilon((m-i)a_1-ia_2)\\
    &= \epsilon(a_1)\epsilon(a_2)\cdot \left(\prod_{i=1}^{4}(-1)^i\epsilon((m-i)a_1-ia_2)\right)^{\frac{m-1}{4}}\\
    &= \sigma_m(a_1,a_2)(y_2).
\end{align*}

When $m\equiv-1\mod 4 $, the sign is given by
\begin{align*}
 \prod_{i=1}^{m-1}(-1)^i\epsilon&((m-i)a_1-ia_2) \cdot \epsilon(ma_2)\cdot \prod_{i=0}^{\frac{m-1}{2}}(-1)^i\epsilon((m-2i)a_1)\\
    =& -\epsilon(a_2)\cdot \prod_{i=1}^{m-1}(-1)^i\epsilon((m-i)a_1-ia_2)\\
    =& \epsilon(a_2)\cdot \left(\prod_{i=1}^{4}(-1)^i\epsilon((m-i)a_1-ia_2)\right)^{\frac{m-3}{4}}\\&\cdot \epsilon((m-(m-2))a_1-(m-2)a_2)\epsilon((m-(m-1))a_1-(m-1)a_2)\\
    =& \epsilon(a_2)\cdot  \epsilon(2a_1-a_2)\epsilon(a_1+2a_2)\\
    =& -\epsilon(a_1)\\
   =&  \sigma_m(a_1,a_2)(y_2)
\end{align*}

For the third fixed point $y_3:=(x_0x_3,x_0x_1,x_1x_2)$, the pairs of the basis of $\sE_m(y_3)$ are of the form
\begin{itemize}
    \item $(x_0^{m-i}x_2^i,x_1^{m-i}x_3^i)$, $i=0,\ldots,m-1$
    \item or $(x_2^{m-i}x_3^i,x_2^ix_3^{m-i})$, $i=0,\ldots,\frac{m-1}{2}$.
\end{itemize}
In the first case this is the representation $\rho_{(m-i)a_1+ia_2}$ which has Euler class $$\epsilon((m-i)a_1+ia_2)\cdot ((m-i)a_1+ia_2)e.$$
In the second case this is the representation $\rho^{(-1)^i}_{(m-2i)a_2}$ which has Euler class $$(-1)^i\epsilon(m-2i)\epsilon(a_2)\cdot(m-2i)a_2e.$$

For the sign,  if $m\equiv 1\mod 4$, then
\begin{align*}
\prod_{i=0}^{m-1}\epsilon&((m-i)a_1+ia_2)\cdot \prod_{i=0}^{\frac{m-1}{2}}\epsilon(a_2)\\
    &=\epsilon(a_1)\epsilon(a_2)\cdot \prod_{i=1}^{m-1}\epsilon((m-i)a_1+ia_2)\\
    &=\epsilon(a_1)\epsilon(a_2)\cdot\epsilon(\prod_{i=1}^4(m-i)a_1+ia_2)^{\frac{m-1}{4}}\\
    &=\epsilon(a_1)\epsilon(a_2)\cdot \epsilon(a_2(3a_1+2a_2)(2a_1+3a_2)a_1)^{\frac{m-1}{4}}\\
    &= \epsilon(a_1)\epsilon(a_2)\cdot \epsilon(a_1a_2(12\cdot 1+13a_1a_2))^{\frac{m-1}{4}}\\
    &= \epsilon(a_1)\epsilon(a_2)\\
    &=\sigma_m(a_1,a_2)(y_3).
\end{align*}
If $m\equiv-1\mod4$, then
\begin{align*}
\prod_{i=0}^{m-1}\epsilon&((m-i)a_1+ia_2)\cdot\prod_{i=1}^{\frac{m-1}{2}}(-1)^i\epsilon((m-2i)a_2)\\
=&-\epsilon(a_1)\prod_{i=1}^{m-1}\epsilon((m-i)a_1+ia_2)\\
=&-\epsilon(a_1) \epsilon(\prod_{i=1}^4(m-i)a_1+ia_2)^{\frac{m-3}{4}}\\&\cdot\epsilon((m-(m-2))a_1+(m-2)a_2)\epsilon((m-(m-1))a_1+(m-2)a_2)\\
=&-\epsilon(a_1)\epsilon(2a_1+a_2)\epsilon(a_1+2a_2)\\
=&-\epsilon(a_2)\\
    =&\sigma_m(a_1,a_2)(y_3).
\end{align*}

For the fourth fixed point $y_4:=(x_0x_3,x_2x_3,x_1x_2)$, the pairs of the basis of $\sE_{m,n}(y_4)$ are of the form
\begin{itemize}
    \item $(x_0^{m-i}x_2^i,x_1^{m-i}x_3^i)$, $i=1,\ldots,m-1$
    \item or $(x_2^m,x_3^m)$
    \item or $(x_0^{m-i}x_1^i,x_0^ix_1^{m-i})$, $i=0,\ldots,\frac{m-1}{2}$.
\end{itemize}
In the first case this is the representation $\rho_{(m-i)a_1+ia_2}$ which has Euler class $$\epsilon((m-i)a_1+ia_2)\cdot ((m-i)a_1+ia_2)e.$$
The second case is $\rho_{ma_2}$ which has Euler class $\epsilon(ma_2)\cdot ma_2e$.
In the last case this is the representation $\rho^{(-1)^i}_{(m-2i)a_1}$ which has Euler class $$(-1)^i\epsilon(m-2i)\epsilon(a_1)\cdot(m-2i)a_1e.$$

For the sign, when $m\equiv1\mod 4$ we get
\begin{align*}
\prod_{i=1}^{m-1}\epsilon&((m-ia_1+ia_2))\cdot \epsilon(ma_2)\cdot \prod_{i=0}^{\frac{m-1}{2}}(-1)^i\epsilon((m-2i)a_1)\\
    &= \epsilon(a_1)\epsilon(a_2)\cdot \prod_{i_1}^{m-1}\epsilon((m-ia_1+ia_2))\\
    &=\epsilon(a_1)\epsilon(a_2)\\
    &=\sigma_m(a_1,a_2)(y_4).
\end{align*}

When $m\equiv-1\mod 4$ we get
\begin{align*}
\prod_{i_1}^{m-1}\epsilon&((m-ia_1+ia_2))\cdot \epsilon(ma_2)\cdot \prod_{i=0}^{\frac{m-1}{2}}(-1)^i\epsilon((m-2i)a_1)\\
    =& -\epsilon(a_2)\cdot \prod_{i_1}^{m-1}\epsilon((m-ia_1+ia_2))\\
    =&-\epsilon(a_2)\cdot \left(\prod_{i=1}^4\epsilon((m-ia_1+ia_2))\right)^{\frac{m-3}{2}}\\&\cdot \epsilon(2a_1+a_2)\epsilon(a_1+2a_2)\\
    =& -\epsilon(a_1)\\
    =&\sigma_m(a_1,a_2)(y_4).
\end{align*}

For the fifth fixed point $y_5:=(x_0^2,x_0x_1,x_1^2)$, the pairs are 
\begin{itemize}
    \item  $(x_0x_2^{m-1-i}x_3^i,x_1x_2^ix_3^{m-1-i})$, $i=0,\ldots ,m-1$
    \item $(x_2^{m-i}x_3^i,x_2^ix_3^{m-i})$, $i=0,\ldots,\frac{m-1}{2}$.
\end{itemize}
In the first case this is $\rho^{(-1)^i}_{a_1+(m-1-2i)a_2}$ which has Euler class $$(-1)^i\epsilon(a_1+(m-1-2i)a_2)\cdot (a_1+(m-1-2i)a_2)e.$$ In this second case this is $\rho^{(-1)^i}_{(m-2i)a_2}$ which has Euler class $$(-1)^i\epsilon((m-2i)a_2)\cdot (m-2i)a_2e.$$

For the sign, if $m\equiv1\mod4$, then  
\begin{align*}
     \epsilon(a_2)\cdot\prod_{i=0}^{m-1}(-1)^i&\epsilon(a_1+(m-1-2i)a_2)\\
    =& \epsilon(a_2)\cdot (-1)^{\frac{m}{2}}\cdot \prod_{i=0}^{m-1}\epsilon(a_1+(m-1-2i)a_2)\\
    =& \epsilon(a_2)\cdot \epsilon(a_1)\\&\cdot \epsilon((a_1+4a_2)(a_1+2a_2)(a_1-2a_2)(a_1-4a_2))^{\frac{m-1}{4}}\\
    =&\epsilon(a_1)\epsilon(a_2)\cdot \epsilon(a_1^2(a_1^2-4a_2))^{\frac{m-1}{4}}\\
    =&\epsilon(a_1)\epsilon(a_2)\\
    =&\sigma_m(a_1,a_2)(y_5).
\end{align*}

If  $m\equiv-1\mod4$, then
\begin{align*}
     \prod_{i=0}^{m-1}(-1)^i\epsilon&(a_1+(m-1-2i)a_2)\cdot \prod_{i=0}^{\frac{m-1}{2}}(-1)^i\epsilon((m-2i)a_2)\\
    =& (-1)^{\frac{m-1}{2}}\cdot \prod_{i=0}^{m-1}\epsilon(a_1+(m-1-2i)a_2)\\
    =&  -\epsilon(a_1+2a_1)^2\epsilon(a_1)\cdot \epsilon((a_1+4a_2)(a_1+2a_2)(a_1-2a_2)(a_1-4a_2))^{\frac{m-3}{4}}\\
    =&-\epsilon(a_1)\\
    =&\sigma_m(a_1,a_2)(y_5).
\end{align*}

For the sixth fixed point $y_6:=(x_0x_3,x_0x_1,x_1x_2)$, the pairs are of the form
\begin{itemize}
    \item $(x_0^{m-1-i}x_1^i)x_2,x_0^ix_1^{m-1-i}x_3)$, $i=0,\ldots,\frac{m-1}{2}$
    \item $(x_0^{m-1-i}x_1^i)x_3,x_0^ix_1^{m-1-i}x_2)$, $i=0,\ldots,\frac{m-3}{2}$
    \item $(x_0^{m-i}x_1^i,x_0^ix_1^{m-i})$, $i=0,\ldots,\frac{m-1}{2}$.
\end{itemize}
The first case is $\rho^{(-1)^i}_{(m-1-2i)a_1+a_2}$ which has Euler class $$(-1)^i\epsilon((m-1-2i)a_1+a_2)\cdot ((m-1-2i)a_1+a_2)e.$$ The second case is $\rho^{(-1)^{i+1}}_{(m-1-2i)a_1-a_2}$ which has  Euler class $$(-1)^{i+1}\epsilon((m-1-2i)a_1-a_2)\cdot ((m-1-2i)a_1-a_2)e$$ and the last one is $\rho_{(m-2i)a_1}^{(-1)^i}$ which has Euler class $$(-1)^i\epsilon((m-2i)a_1)\cdot (m-2i)a_1e.$$

To compute the sign, if $m\equiv1\mod4$, then
\begin{align*}
     \epsilon(a_1)\prod_{i=0}^{\frac{m-1}{2}}(-1)^i&\epsilon((m-1-2i)a_1+a_2)\cdot \prod_{i=0}^{\frac{m-3}{2}}(-1)^{i+1}\epsilon((m-1-2i)a_1-a_2)\\
    =& \epsilon(a_1)\epsilon(a_2)\cdot \prod_{i=0}^{m-1}(-1)^i\epsilon(a_2-(m-1-2i)a_1)\\
    =& \epsilon(a_1)\epsilon(a_2)\cdot (-1)^{\frac{m-1}{2}}\cdot \prod_{i=0}^{m-1}\epsilon(a_2-(m-1-2i)a_1)\\
    =&\epsilon(a_1)\epsilon(a_2)\\
    =&\sigma_m(a_1,a_2)(y_6).
\end{align*}
The last step follows from the calculation for the fifth fixed point, one just needs to swap the roles of $a_1$ and $a_2$.

When $m\equiv-1\mod4$ the sign is 
\begin{align*}
    \prod_{i=0}^{\frac{m-1}{2}}(-1)^i&\epsilon((m-1-2i)a_1+a_2)\\&\hskip30pt\cdot \prod_{i=0}^{\frac{m-3}{2}}(-1)^{i+1}\epsilon((m-1-2i)a_1-a_2))
  \cdot \prod_{i=0}^{\frac{m-1}{2}}(-1)^i\epsilon((m-2i)a_1)\\
    =&\prod_{i=0}^{\frac{m-1}{2}}(-1)^i\epsilon((m-1-2i)a_1+a_2)\cdot \prod_{i=0}^{\frac{m-3}{2}}(-1)^{i+1}\epsilon((m-1-2i)a_1-a_2))\\
    =&\prod_{i=0}^{m-1}(-1)^i(a_2-(m-1-2i)a_1)\\
    =& -\epsilon(a_2)\\    
    =& \sigma_m(a_1,a_2)(y_6).
\end{align*}
Again, the last step follows from the calculation for the fifth fixed point, one just needs to swap the roles of $a_1$ and $a_2$.
\end{proof}

\subsection{Oriented bases for $ \sT_{H_3}(y)$}

For each of the six fixed points $y_1,\ldots,y_6$ we have chosen standard bases $(f_1,f_2)$, $(e_1,e_2, e_3)$ and $(x_0, x_1, x_2, x_3)$ for $F$, $E$ and $V$. These define bases of each term of the resolution \eqref{eq:resulution tangent space}, using our convention for the ordered basis of a tensor product. 
We get the basis 
\[(f_2^{\vee}f_1,f_1^{\vee}f_1,f_2^{\vee}f_2,f_1^{\vee}f_2)\]
of $\End(F):=F^\vee\otimes F$, which we rearrange to 
\[ (f_2^{\vee}f_1,f_1^{\vee}f_2)\oplus (f_1^{\vee}f_1)\oplus(f_2^{\vee}f_2)\]
using an even permutation.
We get the basis 
\[(e_3^{\vee}e_1,e_2^{\vee}e_1,e_1^{\vee}e_1,e_3^{\vee}e_2,e_2^{\vee}e_2,e_1^{\vee}e_2,e_3^{\vee}e_3,e_2^{\vee}e_3,e_1^{\vee}e_3)\]
of $\End(E):=E^\vee\otimes E$.
The following is an even permutation of this basis
\begin{align*}
 (e_3^{\vee}e_1,e_1^{\vee}e_3)\oplus (e_2^{\vee}e_1,e_2^{\vee}e_3)\oplus (e_3^{\vee}e_2,e_1^{\vee}e_2)\oplus(e_3^{\vee}e_3)\oplus (e_2^{\vee}e_2)\oplus (e_1^{\vee}e_1).
\end{align*}
Note that the direct summands in the bases for $\End(F)$ and $\End(E)$ either come in pairs which are dual with respect to the $\sigma$-action or have weight $0$.

We also get a basis of $\Hom(F,E):=F^\vee\otimes E$, which we rearrange using an even permutation
\[(f_2^{\vee}e_1,f_1^{\vee}e_1,f_2^{\vee}e_2,f_1^{\vee}e_2,f_2^{\vee}e_3,f_1^{\vee}e_3)\mapsto (f_2^{\vee}e_1,f_1^{\vee}e_3)\oplus (f_1^{\vee}e_1,f_2^{\vee}e_3)\oplus (f_2^{\vee}e_2,f_1^{\vee}e_2). \]

The first basis element of the three pairs has weight given by the following table.
\begin{equation}\label{eqn:SwapTable}
\vbox{
\begin{center}
\begin{tabular}{ |c|c|c|c|c| } 
\hline
& $\wt(f_2^{\vee}e_1)$& $\wt(f_1^{\vee}e_1)$& $\wt(f_2^{\vee}e_2)$ & $\#$swaps\\
 \hline
$y_1$& $a_1+2a_2$ & $a_1$ & $a_2$&1\\ 
$y_2$& $2a_1+a_2$ & $a_2$ & $a_1$ &1\\ 
$y_3$& $a_1$ & $a_1-2a_2$ & $a_2$ &2\\ 
$y_4$& $2a_1-a_2$ & $-a_2$ & $a_1$&2\\ 
$y_5$& $3a_1$ & $a_1$ & $a_1$ &0\\ 
$y_6$& $3a_2$ & $a_2$ & $a_2$ &3\\ 
 \hline
\end{tabular}
\end{center}
}
\end{equation}

Let $g=f_2^{\vee}e_1,f_1^{\vee}e_1,f_2^{\vee}e_2$. Then $(g,g^*)\otimes (x_0,x_1,x_2,x_3)$ has basis
\[(gx_0,g^*x_0,gx_1,g^*x_1,gx_2,g^*x_2,gx_3,g^*x_3)\] which we rearrange to
\[(gx_0,g^*x_1)\oplus (gx_1,g^*x_0)\oplus (gx_2,g^*x_3)\oplus (gx_3,g^*x_1).\]

We aim to have pairs where the first basis element is the one with weight greater or equal zero, which is always the case for $(gx_0,g^*x_1)$ and $(gx_2,g^*x_3)$.
In the pair $(x_1g,x_0g^*)$ the positive weight comes first iff
\[\wt(x_1g)\ge 0\Leftrightarrow \wt(g)\ge a_1\]
and in the pair $(x_3g,x_2g^*)$ the positive weight comes first iff
\[\wt(x_3g)\ge 0\Leftrightarrow \wt(g)\ge a_2.\]
The number of swaps we need to achieve this is listed in the table \eqref{eqn:SwapTable}.

We call an element of $\sT_{H_3}(y_j)$ that is the image of an element of $\Hom(F,E)\otimes V$ of the form $f^\vee_ie_jx_\ell$ a monomial in $\sT_{H_3}(y_j)$. As for $\sE_{m,n}$, we call a basis $b_1,\ldots, b_{12}$ of $\sT_{H_3}(y_j)$ canonically oriented, resp. anti-canonically oriented if $b_1\wedge\ldots\wedge b_{12}$ maps to $1$, resp. $-1$ in $k$ under the isomorphism $\det\sT_{H_3}(y_j)$ induced by the resolution \eqref{eq:resulution tangent space} and our choice of canonical bases for $\Hom(E,E)$, $\Hom(F,F)$, $\Hom(F, E)\otimes V$ and $1\in k$. 

\begin{proposition} For each $j=1,\ldots, 6$, each basis of $\sT_{H_3}(y_j)$ consisting of $\sigma$-dual pairs of monomial weight vectors is canonically oriented if the number of swaps in \eqref{eqn:SwapTable} is even, and is anti-canonically oriented if if the number of swaps in \eqref{eqn:SwapTable} is odd.
\end{proposition}

\begin{proof} The issue is the same as for Proposition~\ref{prop:SignsInOrientedBases}, namely, the map $\Hom(E,E)\oplus\Hom(F,F)\to\Hom(F, E)\otimes V$ does not send standard basis vectors to standard basis vectors. Just as in the proof of  
Proposition~\ref{prop:SignsInOrientedBases}, one can choose monomials $f^\vee_ie_jx_\ell$  in 
$\Hom(F, E)\otimes V$, closed under $f^\vee_ie_jx_\ell\mapsto (f^\vee_ie_jx_\ell)^*$,  that yield a basis in $\sT_{H_3}(y_j)$. One then changes the monomial basis in $\Hom(F, E)\otimes V$ by a matrix of determinant one that is the identity on this chosen set of monomials,   and that maps the complementary set of monomials to a basis of the image of 
$\Hom(E,E)\oplus\Hom(F,F)$ in $\Hom(F, E)\otimes V$ (the corresponding construction is easy to accomplish for the map $k\to \Hom(E,E)\oplus\Hom(F,F)$). The argument used in the proof of Proposition~\ref{prop:SignsInOrientedBases} then gives the result.
\end{proof}

Suppose we have odd integers $m_1,\ldots, m_r\ge 3$ and $n\ge4$ such that $\sum_i3m_i+1=4n$ and with $\sE_{m_1,\ldots, m_r;n}$ relatively oriented. For the fixed points $y_3$, $y_4$, $y_5$, we have an even number of swaps, and so a basis of
$\sT_{H_3}(y_j)$ constructed in this way will be relatively oriented with respect to our similarly constructed basis of $\sE_{m_1,\ldots, m_r;n}$ (see Remark~\ref{rem:OrientedBases}). For  $y_1$, $y_2$ and $y_6$, we have an odd number of swaps, so our basis of 
$\sT_{H_3}(y_j)$ will yield the ``opposite'' relative orientation with respect to our basis of $\sE_{m_1,\ldots, m_r;n}$, and we will correct this by  including an extra $(-1)$ factor in the Euler class computation for $\sT_{H_3}(y_j)$.

\subsection{Equivariant Euler class of $ \sT_{H_3}(y_i)$ for $i=1,\ldots,6$}
Recall from Theorem \ref{thm:EulerClasses} that for $a$ even 
\[e(\widetilde{\mathcal{O}}(a))=\begin{cases}\frac{a}{2}\widetilde{e}& \text{if }a\equiv2\mod 4\\ -\frac{a}{2}\widetilde{e}& \text{if }a \equiv 0 \mod 4\end{cases}\]
We say that the sign of $e(\widetilde{\mathcal{O}}(a))$ is $+1$ for $a\equiv2\mod 4$ and $-1$ for $\equiv0\mod4$.

Since $a_1$ and $a_2$ are odd we have $e(\widetilde{\mathcal{O}}(2a_1))=a_1\widetilde{e}$ and $e(\widetilde{\mathcal{O}}(2a_2))=a_2\widetilde{e}$. Furthermore one checks that
\begin{align*}&\text{sign}(e(\widetilde{\mathcal{O}}(a_1-a_2)))=-\text{sign}(e(\widetilde{\mathcal{O}}(a_1+3a_2)))=-\text{sign}(e(\widetilde{\mathcal{O}}(3a_1+a_2))).\\
=&-\text{sign}(e(\widetilde{\mathcal{O}}(a_1+a_2)))=\text{sign}(e(\widetilde{\mathcal{O}}(a_1-3a_2)))=\text{sign}(e(\widetilde{\mathcal{O}}(3a_1-a_2)))\end{align*}
In the following computations we use $\pm$ for the sign of the Euler classes of $\widetilde{\mathcal{O}}(a_1-a_2)$, $\widetilde{\mathcal{O}}(a_1+3a_2)$ and $\widetilde{\mathcal{O}}(3a_1+a_2)$
and $\mp$ for the sign of $\widetilde{\mathcal{O}}(a_1+a_2)$, $\widetilde{\mathcal{O}}(a_1-3a_2)$ and $\widetilde{\mathcal{O}}(3a_1-a_2)$ to indicate that the signs are opposite.

\subsubsection{First fixed point $y_1=(x_0x_2,x_0x_1,x_1x_3)$}
The six oriented monomial pairs in $ \sT_{H_3}(y_1)$ yield the following $N$-representations
\begin{itemize}
    \item $(f_1^{\vee}e_1x_0,f_2^{\vee}e_3x_1)=\rho^-_{2a_1}$
    \item $(f_1^{\vee}e_2x_0,f_2^{\vee}e_2x_1)=\rho_{a_1-a_2}$
    \item $(f_2^{\vee}e_1x_2,f_1^{\vee}e_3x_3)=\rho_{3a_2+a_1}$
    \item $(f_2^{\vee}e_2x_2,f_1^{\vee}e_2x_3)=\rho^-_{2a_2}$
    \item $(f_2^{\vee}e_1x_3,f_1^{\vee}e_3x_2)=\rho^-_{a_1+a_2}$
    \item $(f_1^{\vee}e_1x_3,f_2^{\vee}e_3x_2)=\rho_{a_1-a_2}$
\end{itemize}
and we get the following equivariant Euler class
\begin{align*}
\begin{split}
    e^N( \sT_{H_3}(y_1))=& (-1)\cdot(-a_1\widetilde{e})\cdot(\pm \frac{a_1-a_2}{2}\widetilde{e})\cdot  (\pm \frac{a_1+3a_2}{2}\widetilde{e})\\&\hskip50pt\cdot(-a_2 \widetilde{e})\cdot  (\mp (-1) \frac{a_1+a_2}{2}\widetilde{e})\cdot (\pm \frac{a_1-a_2}{2}\widetilde{e})\\
    =& -4a_1a_2(a_1-a_2)^2(a_1+3a_2)(a_1+a_2)e^6.
\end{split}
\end{align*}
Here we use the identity $\widetilde{e}^2=4e^2$ (Theorem~\ref{thm:EulerClasses}) for the last equality. The first $(-1)$ in this computation is needed, because we used an odd permutation to get a basis of desired form.

\subsubsection{Second fixed point $y_2=(x_0x_2,x_2x_3,x_1x_3)$}
The six oriented monomial pairs in $ \sT_{H_3}(y_2)$  yield the following $N$-representations
\begin{itemize}
    \item $(f_2^{\vee}e_2x_3,f_1^{\vee}e_2x_2)=\rho_{a_1-a_2}$
    \item $(f_1^{\vee}e_1x_2,f_2^{\vee}e_3x_3)=\rho^-_{2a_2}$
    \item $(f_2^{\vee}e_2x_0,f_1^{\vee}e_2x_1)=\rho^-_{2a_1}$
    \item $(f_2^{\vee}e_1x_0,f_1^{\vee}e_3x_1)=\rho_{3a_1+a_2}$
    \item $(f_2^{\vee}e_1x_1,f_1^{\vee}e_3x_0)=\rho^-_{a_1+a_2}$
    \item $(f_2^{\vee}e_3x_0,f_1^{\vee}e_1x_1)=\rho_{a_1-a_2}$
\end{itemize}
and thus we get the following equivariant Euler class
\begin{align*}
\begin{split}
    e^N( \sT_{H_3}(y_2))=&(-1)\cdot (\pm \frac{a_1-a_2}{2}\widetilde{e})\cdot (-a_2\widetilde{e})\\
    &\hskip40pt\cdot (-a_1 \widetilde{e})\cdot (\pm \frac{3a_1+a_2}{2}\widetilde{e})\cdot (\mp(-1) \frac{a_1+a_2}{2}\widetilde{e})\cdot (\pm \frac{a_1-a_2}{2}\widetilde{e})\\
    =& -4a_1a_2(a_1-a_2)^2(3a_1+a_2)(a_1+a_2)e^6
\end{split}
\end{align*}
Again we need an additional sign because we used an odd permutation to get the basis that splits up at the direct sum of summands listed above.

\subsubsection{Third fixed point $y_3=(x_0x_3,x_0x_1,x_1x_2)$}
The six oriented monomial pairs in $ \sT_{H_3}(y_3)$  yield the following $N$-representations
\begin{itemize}
    \item $(f_2^{\vee}e_2x_1,f_1^{\vee}e_2x_0)=\rho^-_{a_1+a_2}$
    \item $(f_2^{\vee}e_3x_1,f_1^{\vee}e_1x_0)=\rho^-_{2a_1}$
    \item $(f_2^{\vee}e_2x_2,f_1^{\vee}e_2x_3)=\rho^-_{2a_2}$
    \item $(f_1^{\vee}e_3x_3,f_2^{\vee}e_1x_2)=\rho^-_{a_1-3a_2}$
    \item $(f_1^{\vee}e_3x_2,f_2^{\vee}e_1x_3)=\rho_{a_1-a_2}$
    \item $(f_2^{\vee}e_3x_2,f_1^{\vee}e_1x_3)=\rho^-_{a_1+a_2}$
\end{itemize}
and we get the following equivariant Euler class
\begin{align*}
\begin{split}
    e^N( \sT_{H_3}(y_3))=&(\pm \frac{a_1+a_2}{2}\widetilde{e})\cdot (-a_1\widetilde{e})\cdot (-a_2 \widetilde{e})\\&\hskip40pt\cdot (\pm \frac{a_1-3a_2}{2}\widetilde{e})\cdot (\pm \frac{a_1-a_2}{2}\widetilde{e})\cdot (\pm \frac{a_1+a_2}{2}\widetilde{e})\\
    =&4a_1a_2(a_1+a_2)^2(a_1-a_2)(a_1-3a_2)e^6
\end{split}
\end{align*}

\subsubsection{Fourth fixed point $y_4=(x_0x_3,x_2x_3,x_1x_2)$}
The six oriented monomial pairs in $\sT_{H_3}(y_4)$  yield the following $N$-representations
\begin{itemize}
    \item $(f_2^{\vee}e_2x_3,f_1^{\vee}e_2x_2)=\rho^-_{a_1+a_2}$
    \item $(f_2^{\vee}e_3x_3,f_1^{\vee}e_1x_2)=\rho^-_{2a_2}$
    \item $(f_2^{\vee}e_2x_0,f_1^{\vee}e_2x_1)=\rho^-_{2a_1}$
    \item $(f_2^{\vee}e_1x_0,f_1^{\vee}e_3x_1)=\rho^-_{3a_1-a_2}$
    \item $(f_2^{\vee}e_1x_1,f_1^{\vee}e_3x_0)=\rho_{a_1-a_2}$
    \item $(f_2^{\vee}e_3x_0,f_1^{\vee}e_1x_1)=\rho^-_{a_1+a_2}$.
\end{itemize}
The equivariant Euler class of $ \sT_{H_3}(y_4)$ equals
\begin{align*}
\begin{split}
    e^N( \sT_{H_3}(y_4))=&(\pm\frac{a_1+a_2}{2}\widetilde{e})\cdot (-a_2\widetilde{e})\cdot (-a_1 \widetilde{e})\\&\hskip40pt\cdot (\pm \frac{3a_1-a_2}{2}\widetilde{e})\cdot (\pm \frac{a_1-a_2}{2}\widetilde{e})\cdot (\pm \frac{a_1+a_2}{2}\widetilde{e})\\
    =&4a_1a_2(a_1+a_2)^2(a_1-a_2)(3a_1-a_2)e^6
\end{split}
\end{align*}

\subsubsection{Fifth fixed point $y_5=(x_0^2,x_0x_1,x_1^2)$}
The six oriented monomial pairs in $\sT_{H_3}(y_5)$  yield the following $N$-representations
\begin{itemize}
    \item $(f_2^{\vee}e_2x_3,f_1^{\vee}e_2x_2)=\rho_{a_1-a_2}$
    \item $(f_1^{\vee}e_1x_2,f_2^{\vee}e_3x_3)=\rho^-_{a_1+a_2}$
    \item $(f_2^{\vee}e_2x_2,f_1^{\vee}e_2x_3)=\rho^-_{a_1+a_2}$
    \item $(f_2^{\vee}e_1x_2,f_1^{\vee}e_3x_3)=\rho_{3a_1+a_2}$
    \item $(f_2^{\vee}e_1x_3,f_1^{\vee}e_3x_2)=\rho^-_{3a_1-a_2}$
    \item $(f_1^{\vee}e_1x_3,f_2^{\vee}e_3x_2)=\rho_{a_1-a_2}$.
\end{itemize}
So the equivariant Euler class of $\sT_{H_3}(y_5)$ equals
\begin{align*}
\begin{split}
    e^N( \sT_{H_3}(y_5))=&(\pm\frac{a_1-a_2}{2}\widetilde{e})\cdot (\pm\frac{a_1+a_2}{2}\widetilde{e})\cdot (\pm\frac{a_1+a_2}{2} \widetilde{e})\\&\hskip50pt\cdot (\pm \frac{3a_1+a_2}{2}\widetilde{e})\cdot (\pm \frac{3a_1-a_2}{2}\widetilde{e})\cdot (\pm \frac{a_1-a_2}{2}\widetilde{e})\\
    =&(a_1+a_2)^2(a_1-a_2)^2(3a_1+a_2)(3a_1-a_2)e^6
\end{split}
\end{align*}

\subsubsection{Sixth fixed point $y_6=(x_2^2,x_2x_3,x_3^2)$}
The six oriented monomial pairs in $\sT_{H_3}(y_6)$ yield the following $N$-representations
\begin{itemize}
    \item $(f_1^{\vee}e_2x_0,f_2^{\vee}e_2x_1)=\rho_{a_1-a_2}$
    \item $(f_1^{\vee}e_1x_0,f_2^{\vee}e_3x_1)=\rho^-_{a_1+a_2}$
    \item $(f_2^{\vee}e_2x_0,f_1^{\vee}e_2x_1)=\rho^-_{a_1+a_2}$
    \item $(f_2^{\vee}e_1x_0,f_1^{\vee}e_3x_1)=\rho_{3a_2+a_1}$
    \item $(f_1^{\vee}e_3x_0,f_2^{\vee}e_1x_1)=\rho^-_{a_1-3a_2}$
    \item $(f_2^{\vee}e_3x_0,f_1^{\vee}e_1x_1)=\rho_{a_1-a_2}$.
\end{itemize}
So the equivariant Euler class of $ \sT_{H_3}(y_6)$ is equal to
\begin{align*}
\begin{split}
    e^N( \sT_{H_3}(y_6))=& (-1)\cdot(\pm\frac{a_1-a_2}{2}\widetilde{e})\cdot (\pm\frac{a_1+a_2}{2}\widetilde{e})\cdot (\pm\frac{a_1+a_2}{2} \widetilde{e})\\&\hskip50pt\cdot (\pm \frac{3a_2+a_1}{2}\widetilde{e})\cdot (\pm \frac{a_1-3a_2}{2}\widetilde{e})\cdot (\pm \frac{a_1-a_2}{2}\widetilde{e})\\ 
    =&-(a_1+a_2)^2(a_1-a_2)^2(3a_2+a_1)(a_1-3a_2)e^6
\end{split}
\end{align*}
Note that again we need an additional $(-1)$.

\subsection{Equivariant Euler class of $\sT_{\Gr(4,n+1)}(y_j)$}
To find $e^N(\sT_{H_n}(y_j))$ for $j=1,\ldots,6$ it remains to compute $e^N(\sT_{\Gr(4,n+1)}(y_j))$, which is the same for all 6 fixed points.

\begin{proposition}\label{prop:EulerClassGrassmannian}

Let $s=\floor{\frac{n+1}{2}}$. If $n$ is odd, then
\[e^N(\sT_{\Gr(4,n+1)}(y_j))= \prod_{i=3}^s (a_1^2-a_i^2) \cdot \prod_{i=3}^s(a_2^2-a_i^2)\cdot e^{2s-6}.
\]
If  $n$ is even, then
\[
e^N(\sT_{\Gr(4,n+1)}(y_j))=
\epsilon(a_1a_2)\cdot a_1a_2\prod_{i=3}^s (a_i^2-a_1^2) \cdot \prod_{i=3}^s(a_i^2-a_2^2)\cdot e^{2s-6}.\]
\end{proposition}
\begin{proof}

Since $\sT_{\Gr(4,n+1)}\cong E_4^{\vee}\otimes \mathcal{Q}$ where $E_4$ and $\mathcal{Q}$ are the tautological bundle and the quotient bundle respectively,
$\sT_{\Gr(4,n+1)}(y_j)=(\rho_{a_1}\oplus \rho_{a_2})\otimes (\rho_{a_3}\oplus\ldots \oplus\rho_{a_s})$ if $n$ is odd and $\sT_{\Gr(4,n+1)}(y_j)=(\rho_{a_1}\oplus \rho_{a_2})\otimes (\rho_{a_3}\oplus\ldots \oplus\rho_{a_s}\oplus k)$ if $n$ is even.

So we need to compute the equivariant Euler class of $\rho_{a_i}\otimes \rho_{a_j}$ for $i\neq j$. 
Let $u_1,u_2$ be an oriented basis of the N-representation $\rho_{a_i}$ and $v_1,v_2$ an oriented basis of $\rho_{a_j}$. Then we get the following basis for $\rho_{a_i}\otimes \rho_{a_j}$
\[u_1\otimes v_1,u_2\otimes v_1,u_1\otimes v_2,u_2\otimes v_2\]
which we rearrange (base change determinant 1) to 
\[u_1\otimes v_1,u_2\otimes v_2,u_2\otimes v_1,u_1\otimes v_2.\]
This splits up as the direct sum of two irreducible $N$-representations, namely
one with basis $u_1\otimes v_1,u_2\otimes v_2$ and one with basis $u_2\otimes v_1,u_1\otimes v_2$. The first one is $\rho_{a_i+a_j}$ and the basis is oriented. If $a_i>a_j$ the second one is $\rho_{a_i-a_j}^-$ and the basis above is oriented and if $a_j>a_i$ it is $\rho_{a_j-a_i}^-$ and the basis is not oriented. In both cases we get

\begin{align*}
e^N(\rho_{a_i}\otimes\rho_{a_j})=&(\pm \cdot\frac{a_i+a_j}{2}\cdot \widetilde{e})\cdot(\mp\cdot(-1)\cdot \frac{a_i-a_j}{2}\cdot \widetilde{e})\\
=& (a_i^2-a_j^2)\cdot e^2.
\end{align*}
When $n$ is even, we get the additional summand $\rho_{a_1}\oplus \rho_{a_2}$ which gives the additional factor $\epsilon(a_1)a_1e\cdot \epsilon(a_2)a_2e=\epsilon(a_1a_2)a_1a_2e^2$.
\end{proof}

\appendix
\section{Macaulay2 code and examples}

For $n\le 12$ the following cases are relatively orientable:
\begin{itemize}
\item $n=4$ and $m=5$, $n=10$ and $m=13$
\item $n=5$ and $m=(3,3)$, $n=11$ and $m=(3,11)$, $n=11$ and $m=(3,9)$, $n=11$ and $m=(7,7)$
\item $n=12$ and $m=(3,3,9)$, $n=12$ and $m=(3,5,7)$, $n=12$ and $m=(5,5,5)$
\end{itemize}

To compute \eqref{eq: Bott formula for twisted cubics} we wrote a Macaulay2 function.
The input is an integer $n$ which is the $n$ from $\mathbb{P}^n$ and a list $degs$ of degrees.

\begin{verbatim}
CubicCount=(n,degs)->(
R=QQ[a,b];
Ed=d->(
df=1;
i=1; while i<d+1 do(df=df*i;i=i+2);
P1=1;
i=1; while i<d do(P1=P1*(i*a-(d-i)*b);i=i+1);
F1=df*b^((d+1)//2)*d*a*P1;
P2=P1;
F2=df*a^((d+1)//2)*d*b*P2;
P3=1;
i=1; while i<d do(P3=P3*(i*a+(d-i)*b);i=i+1);
F3=df*b^((d+1)//2)*d*a*P3;
P4=P3;
F4=df*a^((d+1)//2)*d*b*P4;
P5=1;
i=-(d-1)//2; while i<(d-1)//2+1 do (P5=P5*(2*i*b+a);i=i+1);
F5=df*b^((d+1)//2)*P5;
P6=1;
i=-(d-1)//2; while i<(d-1)//2+1 do (P6=P6*(2*i*a+b);i=i+1);
F6=df*a^((d+1)//2)*P6;
{F1,F2,F3,F4,F5,F6});
j=0;E1=1;E2=1;E3=1;E4=1;E5=1;E6=1;
while j<length degs do
(E1=E1*(Ed(degs_j))_0;E2=E2*(Ed(degs_j))_1;E3=E3*(Ed(degs_j))_2;
E4=E4*(Ed(degs_j))_3;E5=E5*(Ed(degs_j))_4;E6=E6*(Ed(degs_j))_5;j=j+1);
T1=-2*a*(a+b)*(b-a)^2*(3*b+a)*2*b;
T2=-2*a*(a+b)*(b-a)^2*(3*a+b)*2*b;
T3=2*a*(b-a)*(b+a)^2*(3*b-a)*2*b;
T4=-2*a*(b-a)*(b+a)^2*(3*a-b)*2*b;
T5=(b+a)*(3*a+b)*(b-a)^2*(b+a)*(3*a-b);
T6=(a+b)^2*(3*b-a)*(3*b+a)*(b-a)^2;
if odd n==true then m=(n+1)//2 else m=n//2;
S=QQ[x_1..x_m];
if odd n==true then TG=1 else TG=x_1*x_2;
i=3; while i<(m+1) do (TG=TG*(x_1^2-x_i^2)*(x_2^2-x_i^2);i=i+1);
ma=map(S,R,{x_1,x_2});
E=((ma E1)/(ma T1)+(ma E2)/(ma T2)+(ma E3)/(ma T3)
+(ma E4)/(ma T4)+(ma E5)/(ma T5)+(ma E6)/(ma T6))/TG;
NumE= numerator E;
DenE= denominator E;
i=1;M={}; while i<(m+1) do (M=M|{x_i};i=i+1);
Sub=subsets(M,2);
l2=length Sub; T=E;
i=1; while i<l2 do
(Mc=toList(set M-set Sub_i);Per=Sub_i|Mc;m=map(S,S,Per); T=T+(m NumE)/(m DenE); i=i+1);
T)
\end{verbatim}
We computed the count in all relatively oriented cases where $n\le 12$ using our Macaulay2 function.

\begin{verbatim}
i2 : CubicCount(4,{5})

o2 = 765

o2 : frac(S)

i3 : CubicCount(10,{13})

o3 = 768328170191602020

o3 : frac(S)

i4 : CubicCount(5,{3,3})

o4 = 90

o4 : frac(S)

i5 : CubicCount(11,{3,11})

o5 = 4407109540744680

o5 : frac(S)

i6 : CubicCount(11,{5,9})

o6 = 313563865853700

o6 : frac(S)

i7 : CubicCount(11,{7,7})

o7 = 136498002303600

o7 : frac(S)

i8 : CubicCount(12,{3,3,9})

o8 = 43033957366680

o8 : frac(S)

i9 : CubicCount(12,{3,5,7})

o9 = 5860412510400

o9 : frac(S)

i10 : CubicCount(12,{5,5,5})

o10 = 1833366298500

o10 : frac(S)
\end{verbatim}

\end{document}